%% file: ms.tex
\pdfoutput=1

\documentclass[a4paper, 11pt, abstracton, titlepage, oneside]{scrartcl}
\makeatletter

\input{resources/packages.tex}
\input{resources/configuration.tex}
\input{resources/commands.tex}
\begin{document}
\emergencystretch 3em
\input{./resources/abbreviations}
\input{./resources/notation}
\input{./resources/titlepage.tex}
\input{./contents/introduction.tex}
\input{./contents/problemsetting.tex}
\input{./contents/pipeline.tex}
\input{./contents/experimentaldesign.tex}
\input{./contents/results.tex}
\input{./contents/conclusion.tex}
\input{./contents/appendix.tex}


%
\singlespacing{
\bibliographystyle{model5-names}
\bibliography{References}} 
\newpage
\onehalfspacing
\begin{appendices}
	\normalsize
\end{appendices}
\end{document}

%% file: resources/packages.tex
\usepackage[T1]{fontenc}
\usepackage[utf8]{inputenc}
\usepackage[onehalfspacing]{setspace}
\usepackage[headsepline]{scrlayer-scrpage}
\clearpairofpagestyles
\usepackage{layout}
\usepackage{geometry}
\usepackage{adjustbox}
\usepackage{authblk}
\usepackage[titletoc,title]{appendix}

\usepackage[hyphens]{url}
\usepackage{color}

\usepackage{amsmath,amssymb,amsfonts,amsthm, mathtools, bm}
\usepackage{thmtools, thm-restate}
\usepackage{nicefrac}
\usepackage{array}

\usepackage[normal]{threeparttable}
\usepackage{threeparttablex}
\usepackage{tabularx}
\usepackage{longtable}
\usepackage{booktabs}
\usepackage{colortbl}
\usepackage{multirow}
\usepackage{makecell}

\usepackage[ruled, noend, linesnumbered]{algorithm2e}
\SetKw{Continue}{continue}

\usepackage[acronym]{glossaries}

\usepackage{caption}
\usepackage{subcaption}
\usepackage{float}
\usepackage{graphicx}
\usepackage{graphics}
\usepackage{placeins}

\usepackage{tikz}
\usetikzlibrary{arrows.meta}
\usetikzlibrary{math}
\usetikzlibrary{shapes}
\usepackage{pgfplots}
\usepackage{svg}
\usepgfplotslibrary{groupplots}
\usepgfplotslibrary{dateplot}
\usepackage{tikzscale}

\usepackage{enumerate}
\usepackage{multicol}
\usepackage[author=Patrick]{pdfcomment}
\usepackage{xparse}
\usepackage{twoopt}
\usepackage{etoolbox}

\usepackage{eurosym}
\usepackage{ellipsis}
\usepackage{natbib}
\usepackage{paralist}
\usepackage{ulem}

\usepackage{import}
\usepackage{xcolor}
\usepackage{bbm}
\usepackage{rotating}
\usepackage{enumitem}

%% file: resources/configuration.tex
\AtBeginDocument{
	\newgeometry{
		textheight=9in,
		textwidth=6.5in,
		top=1in,
		headheight=14pt,
		headsep=25pt,
		footskip=30pt
	}
}
\newsavebox{\abstractbox}
\renewenvironment{abstract}
{\begin{lrbox}{0}\begin{minipage}{\textwidth}
			\begin{center}\normalfont\sectfont\abstractname\end{center}\quotation}
		{\endquotation\end{minipage}\end{lrbox}%
	\global\setbox\abstractbox=\box0 }

\makeatletter
\expandafter\patchcmd\csname\string\maketitle\endcsname
{\vskip\z@\@plus3fill}
{\vskip\z@\@plus2fill\box\abstractbox\vskip\z@\@plus1fill}
{}{}
\makeatother

\addtokomafont{captionlabel}{\footnotesize\bfseries} 
\addtokomafont{caption}{\footnotesize\bfseries} 

\bibpunct[, ]{(}{)}{,}{a}{}{,}%
\newtheorem{result}{Result}

\DeclareMathOperator*{\argmax}{arg\,max}

\counterwithin*{equation}{section}

\newenvironment{myfont}{\fontfamily{ccr}\selectfont}{\par}
\DeclareTextFontCommand{\textmyfont}{\myfont}

\AtBeginDocument{%
	\abovedisplayskip=7.5pt plus 0pt minus 0pt
	\abovedisplayshortskip=7.5pt plus 0pt minus 0pt
	\belowdisplayskip=7.5pt plus 0pt minus 0pt
	\belowdisplayshortskip=7.5pt plus 0pt minus 0pt
}

\newcolumntype{L}[1]{>{\raggedright\let\newline\\\arraybackslash\hspace{0pt}}p{#1}}
\newcolumntype{C}[1]{>{\centering\let\newline\\\arraybackslash\hspace{0pt}}p{#1}}
\newcolumntype{R}[1]{>{\raggedleft\let\newline\\\arraybackslash\hspace{0pt}}p{#1}}
\renewcommand{\emph}[1]{\textit{#1}}

%% file: resources/commands.tex
\newcommand{\resultAvgSBGreedy}{6.3}

\newcommand{\resultVehBASECBGreedy}{3.8}
\newcommand{\resultVehBASESBGreedy}{4.6}
\newcommand{\resultVehUPTOSBGreedy}{14.8}
\newcommand{\resultVehFiveHundertSBGreedy}{14.8}

\newcommand{\resultReqMAXSBGreedy}{17.1}

\newcommand{\resultReqHIGHERSBGreedy}{13.1}
\newcommand{\resultReqHIGHERSamplingGreedy}{10.3}
\newcommand{\resultVehUPTOSBGreedyCustomerSatisfaction}{16.8}

\newcommand{\resultOverallMAXSBGreedy}{17.1}
\newcommand{\resultOverallMAXCBGreedy}{9.1}

\newcommand{\instanceExtractionEvery}{four} 
\newcommand{\numTrainingInstances}{80}
\newcommand{\numPerturbations}{50}

\newcommand{\distanceSparsification}{1.5} 
\newcommand{\temporalSparsification}{720.0} 
\newcommand{\predictionHorizonSampling}{four} 
\newcommand{\reductionFactorSampling}{0.2}
\newcommand{\predictionHorizonSB}{four} 
\newcommand{\predictionHorizonCB}{four} 
\newcommand{\maxCapacityCB}{one}

\newcommand{\systemTimePeriod}{one} 
\newcommand{\lookUpX}{200} 
\newcommand{\lookUpY}{200} 
\newcommand{\averageDrivingSpeed}{20} 
\newcommand{\rebalancingX}{333} 
\newcommand{\rebalancingY}{333} 

%% file: resources/abbreviations.tex
\newacronym{abk:amod}{AMoD}{autonomous mobility-on-demand}
\newacronym{abk:k-dSPP}{$k$-dSPP}{$k$-disjoint shortest paths problem}
\newacronym{abk:ML}{ML}{machine learning}
\newacronym{abk:CO}{CO}{combinatorial optimization}
\newacronym{abk:SL}{SL}{structured learning}
\newacronym{abk:SB}{SB}{sample-based}
\newacronym{abk:CB}{CB}{cell-based}
\newacronym{abk:DRL}{DRL}{deep reinforcement learning}
\newacronym{abk:BFGS}{BFGS}{Broyden–Fletcher–Goldfarb–Shanno}
\newacronym{abk:SGD}{SGD}{stochastic gradient descent}
\newacronym{abk:FNN}{FNN}{feedforward neural network}

%% file: resources/titlepage.tex
\title{\large Learning-based Online Optimization for Autonomous Mobility-on-Demand Fleet Control}

\author[1]{\normalsize Kai Jungel}
\author[2]{\normalsize Axel Parmentier}
\author[3]{\normalsize Maximilian Schiffer}
\author[4]{\normalsize Thibaut Vidal}
\affil{\small 
	TUM School of Management, Technical University of Munich, 80333 Munich, Germany
	
	\scriptsize kai.jungel@tum.de

    \small
	\textsuperscript{2}CERMICS, \'{E}cole des Ponts, Marne-la-Vall\'{e}e, France
	
	\scriptsize axel.parmentier@enpc.fr

	\small
	\textsuperscript{3}TUM School of Management \& Munich Data Science Institute,
	
	Technical University of Munich, 80333 Munich, Germany
	
	\scriptsize schiffer@tum.de
 
    \small
	\textsuperscript{4}CIRRELT \& SCALE-AI Chair in Data-Driven Supply Chains, Department of Mathematics and Industrial Engineering, \'Ecole Polytechnique de Montr\'eal, Montr\'eal, Canada,
 
    Department of Computer Science, Pontifical Catholic University of Rio de Janeiro, Rio de Janeiro, Brazil
	
	\scriptsize thibaut.vidal@cirrelt.ca
 }

\date{}

\lehead{\pagemark}
\rohead{\pagemark}

\begin{abstract}
\begin{singlespace}
{\small\noindent \input{./contents/abstract.tex}
\smallskip}
\end{singlespace}
\end{abstract}

\maketitle

%% file: contents/abstract.tex
\Acrlong*{abk:amod} systems are a viable alternative to mitigate many transportation-related externalities in cities, such as rising vehicle volumes in urban areas and transportation-related pollution.
However, the success of these systems heavily depends on efficient and effective fleet control strategies.
In this context, we study online control algorithms for \acrlong*{abk:amod} systems and develop a novel hybrid 
\acrlong*{abk:CO} enriched \acrlong*{abk:ML} pipeline which learns online dispatching and rebalancing policies from optimal full-information solutions.
We test our hybrid pipeline on large-scale real-world scenarios with different vehicle fleet sizes and various request densities.
We show that our approach outperforms state-of-the-art greedy, and model-predictive control approaches with respect to various KPIs, e.g., by up to~\resultReqMAXSBGreedy\% and on average by~\resultAvgSBGreedy\% in terms of realized profit.

%% file: contents/introduction.tex
\section{Introduction}
\label{sec:Intro}

Growing urbanization leads to rising traffic volumes in urban areas and challenges the sustainability and economic efficiency of today's urban transportation systems. In fact, transportation systems' externalities, such as congestion and local emissions, can significantly damage public health and pose a growing threat to our environment \citep{Levy2010, MobilityReport2019}. 
Moreover, private cars currently require considerable space for parking and driving in urban areas.

Recent progress in autonomous driving and 5G technology has led to crucial advances towards enabling \gls*{abk:amod} systems~\citep{White2020}, which can contribute to mitigating many of the aforementioned challenges in the coming years.
An \gls*{abk:amod} system consists of a centrally-controlled fleet of self-driving vehicles serving on-demand ride requests~\citep{Pavone2015}. 
To manage the fleet, a central controller dispatches ride requests to vehicles in order to transport customers between their origin and destination. The central controller may further relocate idle vehicles to better match future demand, e.g., in high-demand areas.
An \gls*{abk:amod} system has many advantages over current conventional ride-hailing systems. Among others, it permits centralized control and therefore enables system-optimal dispatching decisions. Moreover, it allows for a high utilization of vehicles which results from rebalancing actions. However, the utilization of an \gls*{abk:amod} system depends on its efficient operation, i.e., system optimal online dispatching and rebalancing of the vehicle fleet.
Designing such an efficient control algorithm for \gls*{abk:amod} systems faces two major challenges. First, the underlying planning problem is an online control problem, i.e., the central controller must take dispatching and rebalancing decisions under uncertainty without knowing future ride requests.
Second, the scalability of the control mechanism for large-scale applications can be computationally challenging as dozens of thousands of requests will need to be satisfied every hour in densely populated areas.

Against this background, we develop a novel family of scalable online-control policies that receive an online \gls*{abk:amod} system state as input and return dispatching and rebalancing actions for the \gls*{abk:amod} fleet.
Specifically, we encode these online control policies via a hybrid \gls*{abk:CO} enriched \gls*{abk:ML} pipeline.
The rationale of this hybrid \gls*{abk:CO}-enriched \gls*{abk:ML} pipeline is to transform the online \gls*{abk:amod} system state into a combinatorial dispatching problem, and to learn the parameterization of this dispatching problem from optimal full-information solutions, based on historical data.
As a result, the solution of the dispatching problem imitates optimal full-information solutions and leads to anticipative online dispatching and rebalancing actions.

\subsection{Related Work}
\label{RelatedWork}

Our work contributes to two different research streams. From an application perspective, it connects to vehicle routing and dispatching problems for \gls*{abk:amod} systems. From a methodological perspective, it connects to the field of \gls*{abk:CO}-enriched \gls*{abk:ML}.
In this section, we first review related works in the field of \gls*{abk:amod} systems and refer to \cite{Vidal2020} for a general overview of vehicle routing problems. Second, we review related works at the intersection of \gls*{abk:ML} and \gls*{abk:CO}, and refer to \cite{Bengio2021} and \cite{Kotary2021} for a general overview.

\noindent \textbf{AMoD Systems.}
\gls*{abk:amod} systems have been studied from a system perspective and from an operational perspective.
Studies that focused on the system perspective investigated offline problem settings with full information on all ride requests, addressing research questions related to strategic decision-making and mesoscopic analyses, and are only loosely related to our work. Therefore, we exclude these from our review and refer to \cite{Zardini2022} for a general overview.

This study relates to the operational perspective, specifically to online systems with ride requests entering the system sequentially over time. To solve such an online dispatching problem, classical \gls*{abk:CO} approaches which leverage the combinatorial structure of the problem, as well as \gls*{abk:ML} approaches which learn to account for the uncertain appearance of future ride requests exist.
In the context of classical \gls*{abk:CO}, \cite{Bertsimas2019} introduced an online re-optimization algorithm, which iteratively optimizes the dispatching of available vehicles to new ride requests when they enter the system in a rolling-horizon fashion.
Other works focused explicitly on relocating unassigned vehicles to areas of high future demand to reduce waiting times for future ride requests \citep{Pavone2012, Zhang2016, Iglesias2018, Yang2022}.
As the efficiency of dispatching and rebalancing depends on future ride requests, model-predictive control \citep{Alonso-Mora2017a} 
and stochastic model-predictive control algorithms
\citep{Tsao2018} have been developed to incorporate information about 
future ride requests. A complementary stream of research studied similar problems from a \gls*{abk:DRL} perspective, explicitly focusing on rebalancing \citep{Jiao2021, 
Gammelli2021, 
Skordilis2022, 
Liang2022} 
or non-myopic (anticipative) dispatching \citep{
Xu2018, 
Li2019, 
Tang2019, 
Zhou2019, 
Eshkevari2022, 
Enders2022}. 

As can be seen, non-myopic dispatching and rebalancing for \gls{abk:amod} systems has recently been vividly studied. Existing approaches either follow a \gls*{abk:DRL} paradigm to account for uncertain system dynamics, or focus on classical optimization-based (control) approaches that are often extended by a sequential prediction component. However, a truly integrated approach that aims to combine learning and optimization to leverage the advantages of both concepts simultaneously has not been studied so far.

\noindent
\textbf{Combinatorial Optimization enriched Machine Learning.}
In many applications, combinatorial decisions are subject to uncertainty.
In this context, various approaches that enrich \gls*{abk:CO} by \gls*{abk:ML} components to account for uncertainty have recently been studied. For example, \cite{Donti2017} proposed an end-to-end learning approach for probabilistic models in the course of stochastic optimization, and \cite{Elmachtoub2022} introduced a loss to learn a prediction model for travel times by utilizing the structure of the underlying vehicle routing problem. Closest to our paradigm are approaches to which we can generally refer as \gls*{abk:CO}-enriched \gls*{abk:ML} pipelines \citep[cf.][]{Dalle2022}.
In particular, these approaches revolve around the concept of \acrfull{abk:SL}, a supervised learning methodology with a combinatorial solution space \citep[cf.][]{Nowozin2010}, which can be used to enrich \gls*{abk:CO}-algorithms with a prediction component. 
\cite{Parmentier2021a}~used \gls*{abk:SL} 
to parametrize a computationally tractable surrogate problem such that its solution approximates the solution of a hard combinatorial problem.
Using this concept, \cite{Parmentier2022} solved a single-machine scheduling problem with release dates by transforming it into a computationally tractable scheduling problem without release dates.

These recent works point towards the efficiency of \gls*{abk:CO}-enriched \gls*{abk:ML} when determining approximate solutions for hard-to-solve combinatorial problems. However, it remains an open question whether \gls*{abk:CO}-enriched \gls*{abk:ML} can also successfully be used to derive efficient online control policies, especially in the realm of \gls*{abk:amod} systems.

\subsection{Contributions}
\label{sec:contributions}

In this paper, we introduce a new state-of-the-art algorithm for large-scale online dispatching and rebalancing of \gls*{abk:amod} systems.
To reach this goal, we close the research gaps outlined above by introducing a novel family of online control policies, encoded via a hybrid \gls*{abk:CO}-enriched \gls*{abk:ML} pipeline.
In a nutshell, this hybrid \gls*{abk:CO}-enriched \gls*{abk:ML} pipeline receives an \gls*{abk:amod} system state, transforms it into a combinatorial dispatching problem, parameterizes the dispatching problem by utilizing \gls*{abk:SL}, and solves it to return online dispatching and rebalancing actions. The main rationale of this hybrid \gls*{abk:CO}-enriched \gls*{abk:ML} pipeline is to learn the parametrization of the dispatching problem from optimal full-information solutions for historical data, in such a way that the resulting policy imitates fully-informed dispatching and rebalancing actions.

The contribution of our work is several-fold.
First, we present a novel \gls*{abk:CO}-enriched \gls*{abk:ML} pipeline, amenable to derive online control policies, and tailor it to online dispatching and rebalancing for \gls*{abk:amod} systems.
Second, to establish this pipeline, we show how to solve the underlying \gls*{abk:CO} problem, namely a dispatching and rebalancing problem, in polynomial time.
Third, we show how to learn a parametrization for the underlying \gls*{abk:CO} problem based on optimal full-information solutions from historical data using an \gls*{abk:SL}-approach. 
Fourth, we derive two online control policies, namely a \gls*{abk:SB} approach, and a \gls*{abk:CB} approach, both utilizing our hybrid \gls*{abk:CO}-enriched \gls*{abk:ML} pipeline.
Fifth, we present a numerical study based on real-world data to benchmark the encoded policies against current state-of-the-art benchmarks: a greedy policy and a sampling-based receding horizon approach \citep[cf.][]{Alonso-Mora2017a}. 

Our results show that our learning-based policies outperform a greedy policy with respect to realized profit by up to~\resultReqMAXSBGreedy\% and by~\resultAvgSBGreedy\% on average. Moreover, our learning-based policies robustly perform better than a greedy policy across various scenarios, whereas a sampling-based policy oftentimes performs worse than greedy. In general, a sampling-based policy shows an inferior performance compared to our learning-based policies, even in the scenarios where it improves upon greedy. Finally, we show that our learning-based policies yield a good trade-off between anticorrelated key performance indicators, e.g., when maximizing the total profit or the number of served customers, we still obtain a lower distance driven per vehicle compared to policies that show worse performance on the aforementioned quantities.

\subsection{Organization}
\label{Organization}

The remainder of this paper is structured as follows.
In Section~\ref{sec:ProblemSetting}, we detail our problem setting. In~Section~\ref{sec:Pipeline}, we introduce our hybrid \gls*{abk:CO}-enriched \gls*{abk:ML} pipeline and derive the \gls*{abk:SL} problem to learn an \gls*{abk:ML}-predictor which parameterizes the underlying \gls*{abk:CO} problem. We present a real-world case study in Section~\ref{sec:computationalExperiments} and its numerical results in Section~\ref{sec:results}. Section~\ref{Conclusion} finally concludes the paper.

%% file: contents/problemsetting.tex
\section{Problem setting}
\label{sec:ProblemSetting}
We study an online vehicle routing problem from the perspective of an \gls*{abk:amod} system's central controller, which steers a homogeneous fleet of self-driving vehicles. The central controller either dispatches vehicles to serve ride requests or rebalances idle vehicles to anticipate future ride requests. 
Here, a ride request can represent a solitary ride or pooled rides as we assume decisions on pooling requests to be taken at an upper planning level. Accordingly, each ride request demands a mobility service from an origin to a destination location at a certain point in time.
The central controller takes decisions according to a discrete time horizon~$\{1,\ldots,T\}$.
We denote the period between the actual decision time~$t$ and the subsequent decision time~$t+1$ as system time period~$[t, t+1)$.
We consider a set of vehicles~$\mathcal{V}$ that serve requests.
Each ride request~$r\in R$ is a quintuple~$r=(o_r, d_r, s_r, a_r, p_r)$ with a pick-up location~$o_r$, a drop-off location~$d_r$, a start time~$s_r$, an arrival time~$a_r$, and its reward~$p_r$.
Here, the start time~$s_r$, and the arrival-time~$a_r$ are the times when the vehicle reaches the origin~$o_r$ and destination~$d_r$ of the ride request respectively. 
Let~$\tau(l_1,l_2)$ denote the driving time of a vehicle between two locations~$l_1$ and~$l_2$.
The start time~$s_r$, the arrival time~$a_r$, and the driving time~$\tau(l_1,l_2)$ are not subject to the time discretization but are from a continuous time range in~$\mathbb{R}$.
If the central controller assigns a ride request to a vehicle, the vehicle has to fulfill this ride request. We refer to a sequence of ride requests that a vehicle fulfills as a \emph{trip}.

\paragraph{System state:}
At any decision time~$t$, the state~$x_{v,t}$ of a vehicle~$v$ is a pair~$\big(l_{v,t},S_{v,t}\big)$ composed of its current location~$l_{v,t}$ and a sequence~$S_{v,t}=(r_1, r_2, ..., r_m)$ of unfinished requests which have been assigned to~$v$ at the previous decision time~$t-1$. A request~$r$ is unfinished when~$a_r > t$. 
In our setting, the sequence of unfinished requests contains~$|S_{v,t}| \in \{0,1\}$ elements as the central controller could only assign requests which started before the decision time~$t$, while in a general setting when the central controller could assign requests which start after the decision time~$t$, the sequence contains~$|S_{v,t}| \in \mathbb{N}_0$ elements. The sequence is empty when the vehicle has already fulfilled all assigned requests, which means that the vehicle is idling.
Regarding the overall \gls*{abk:amod} system state, the state~$x_t = (R_t,(x_{v,t})_{v \in V})$ comprises the batch~$R_t$ of requests that entered the system between the previous decision time~$t-1$ and~$t$ with~$s_r \in [t, t+1)$, as well as the state of each vehicle~$v \in \mathcal{V}$.

\paragraph{Feasible decisions:}
Given the current state of the system~$x_t$, the central controller decides at decision time~$t$ which requests~$r \in R_t$ to accept and assigns them to vehicles.
The central controller can also rebalance a vehicle after the last request in its trip or when it is idling.
A feasible decision~$y_{v,t}$ for vehicle~$v$ given its current state~$x_{v,t} = (l_{v,t}, S_{v,t})$ is a \emph{trip}, i.e.,
a sequence of requests~$\tilde S_{v,t}=(r_1, r_2, \dots, r_m, r_{m+1},\ldots , r_n,l)$. This trip~$\tilde S_{v,t}$ includes the unfinished requests which have already been assigned to~$v$ at the previous decision time, i.e., from the sequence~$S_{v,t}$, which is then potentially followed by new requests~$r_{m+1},\ldots , r_n$ from batch~$R_t$, and may include a final rebalancing location~$l$. A trip satisfies the following \emph{constraints}:
\begin{itemize}
    \item [(i)] Each request in a trip does not belong to a trip of another vehicle, meaning~$y_v \cap y_{v'} \cap R_t = \emptyset \text{ for all } v \neq v'$ where~$y_v \cap y_{v'} \cap R$ denotes the set of requests in~$y_v$ and~$y_{v'}$. Note that two vehicles can have the same rebalancing location.
    \item [(ii)] For every pair~$(r_i,r_{i+1})$ of successive requests, the vehicle can reach the origin of~$r_{i+1}$ on time after~$r_i$, i.e.,~$a_{r_i} + \tau(d_{r_i}, o_{r_{i+1}}) \leq s_{r_{i+1}}$.
    \item [(iii)] If the first request~$r_1$ is new, i.e.,~$r_1 \in R_t$, then it must be reachable from the current location of the vehicle:~$t+\tau(l_{v,t},o_{r_1}) \leq s_{r_1}$.
\end{itemize}
We denote by~$\mathcal{Y}^{v}(x_{v,t},R_t)$ the set of trips available for~$v$ given the vehicle state~$x_{v,t}$ and the set of new ride requests~$R_t$.
Then, the set of feasible dispatching decisions given the current state of the system is
\begin{equation*}
    \mathcal{Y}(x_t) = \Big\{(y_{v,t})_{v \in V}\colon y_{v,t} \in \mathcal{Y}^{v}(x_{v,t}, R_t) \text{ for all } v \in \mathcal{V} \Big\}.
\end{equation*}

\paragraph{Evolution of the system:}
Once the central controller dispatches the decision~$y_t$ from decision time~$t$, the system evolves until~$t+1$. 
Each vehicle~$v$ drives along its trip.
If a vehicle completes its trip before~$t+1$, it idles. If the vehicle does not reach its rebalancing location~$l$ before~$t+1$, it starts idling at~$t+1$ and awaits the next operator's decision.
We can therefore deduce from~$y_{v,t}$ its position~$l_{v,t+1}$ at~$t+1$. 
The system removes all requests from~$R_t$ that were assigned at time~$t$.
Furthermore, requests from~$R_t$ which remained unassigned at~$t$ leave the system. New requests~$R_{t+1}$ enter the system between~$t$ and~$t+1$.

\paragraph{Policy:}
Let~$\mathcal{X}$ denote the set of possible \gls*{abk:amod} system states.
A (deterministic) policy~$\delta: \mathcal{X} \rightarrow \mathcal{Y}$ is a mapping that assigns to any \gls*{abk:amod} system state~$x$ a decision~$y \in \mathcal{Y}(x)$.

\paragraph{Objective:}
Let us denote by~$p_{v,t}$ the reward made by vehicle~$v$ at time~$t$.
The reward at time~$t$ is therefore
\begin{equation}
    P_t = \sum_{v \in V} p_{v,t}.
    \label{eq:control_objective}
\end{equation}
We may adjust the objective by changing the interpretation of~$p_{v,t}$; for example, we can maximize the total profit by accounting for trip-dependent rewards and costs.
Alternatively, we can maximize the number of served requests by setting~$p_{v,t}$ to the number of requests that vehicle $v$ serves at time~$t$.
Our objective is to find a policy~$\delta$ that maximizes the total reward over all vehicle trips at the end of the problem time horizon~$T$,
\begin{equation}
    \max_{\delta} \; \mathbb{E} \left( \sum_{t=1}^T P_t | \delta \right).
    \label{eq:ControlProblem}
\end{equation}

\paragraph{Full-information upper bound:}
We can compute a full-information upper bound by supposing that all requests~$\bigcup_t R_t$ are already in the system at time~$t=1$.

\paragraph{Scope and limitations of the model:}
A few comments on this modeling approach are in order. 
First, we do not permit late arrivals at pick-up locations. When a request enters the system, the central controller decides to either fulfill or reject this request. This assumption aligns with real-world practice: customers typically can choose between various mobility providers such that an unfulfilled request usually leads to a customer choosing another provider's service or resubmitting the request a bit later. Consequently, operators aim to take immediate action to provide direct feedback to customers. Notably, it is possible to adapt our setting for delayed pickups by including unfulfilled requests in the request batch that enters the system in the next system time period.
Second, we assume that ride-pooling decisions occur at an upper planning stage.
Analyzing whether a hierarchical or an integrated ride pooling approach is superior for \gls*{abk:CO}-enriched \gls*{abk:ML} online control of \gls*{abk:amod} systems remains an open research question for future work.
Third, we consider travel time uncertainty on a lower operational level. We can express time uncertainty on the duration of each ride request which has no impact on the presented model.

%% file: contents/pipeline.tex
\section{CO-enriched ML pipeline}
\label{sec:Pipeline}
\begin{figure}[b]
    \centering
    \fontsize{9}{10}\selectfont
    \def\svgwidth{1\columnwidth}
    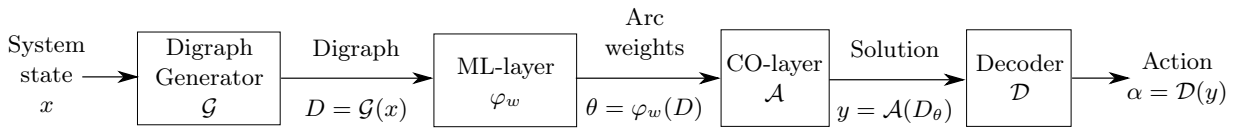
    \caption{The hybrid \gls*{abk:CO}-enriched \gls*{abk:ML} pipeline.}
    \label{fig:MLforOR}
\end{figure}

Recall that a policy maps an \gls*{abk:amod} system state to a feasible dispatching and rebalancing decision.
Since the set of feasible dispatching and rebalancing decisions~$\mathcal{Y}(x)$ is composed of collections of disjoint trips, which can be represented as a collection of disjoint paths in a digraph with requests as vertices, it is tempting to model the problem of choosing a decision as a \gls*{abk:k-dSPP}.
Unfortunately, there is no natural way of defining the parameters of this problem.
To mitigate this, we propose a learning pipeline illustrated in Figure~\ref{fig:MLforOR} to encode our policies.
It includes a \gls*{abk:k-dSPP} as a combinatorial optimization layer that selects the decisions. 
This layer captures the combinatorial structure of the decisions and can be efficiently solved using polynomial algorithms.

Our pipeline is as follows. In general, the input to the hybrid \gls*{abk:CO}-enriched \gls*{abk:ML} pipeline is a system state~$x$. In our specific application (cf. Section~\ref{sec:ProblemSetting}), we can interpret a system state~$x$ as an~\gls*{abk:amod} system state at time~$t$,~$x_t = (R_t,(x_{v,t})_{v \in V})$, which comprises the current ride request batch~$R_t$ and the vehicle states~$x_{v,t}$ of all vehicles~$v \in \mathcal{V}$.
The pipeline forwards this system state~$x$ to the digraph generator layer~$\mathcal{G}$, which transforms the system state into a digraph~$D$ so that a path in the digraph~$D$ is a feasible vehicle trip.
The second layer, the \gls*{abk:ML}-layer, receives the digraph~$D$ and uses an \gls*{abk:ML}-predictor~$\varphi_{w}$ to predict the respective arc weights~$\theta$ of digraph~$D$. Accordingly, the \gls*{abk:ML}-layer returns a weighted digraph~$D_{\theta}$, which we refer to as the \gls*{abk:CO}-layer-digraph.
The third layer, the \gls*{abk:CO}-layer, then consists of a~$k$-dSP-algorithm~$\mathcal{A}$ that solves a~\gls*{abk:k-dSPP} on the weighted \gls*{abk:CO}-layer-digraph~$D_{\theta}$ and returns a solution~$y$. 
Finally, a decoder~$\mathcal{D}$ transforms the vehicle trips~$y$ to a feasible dispatching and rebalancing action~$\alpha$.

The main rationale of this hybrid \gls*{abk:CO}-enriched \gls*{abk:ML} pipeline is to learn the weights~$\theta$ of the \gls*{abk:CO}-layer-digraph~$D_{\theta}$ based on optimal full-information solutions so that the \gls*{abk:k-dSPP} returns a good solution~$y \in \mathcal{Y}(x)$.

In the remainder of this section, we show in Section~\ref{sec:kdspp} how the digraph generator layer~$\mathcal{G}$ generally models a dispatching problem as a digraph so that we can set its weights in the \gls*{abk:ML}-layer and solve it with a \gls*{abk:k-dSPP} in the \gls*{abk:CO}-layer.
We then detail two variants of the \gls*{abk:CO}-enriched \gls*{abk:ML} pipeline in Section~\ref{sec:SP}, and Section~\ref{sec:CP} respectively. These two pipelines represent two different possibilities to enable rebalancing actions and only differ in their digraph generator layer~$\mathcal{G}$. In Section~\ref{sec:StructuredLearning}, we finally derive the \gls*{abk:SL} problem to learn an \gls*{abk:ML}-predictor for these pipelines.

\subsection{Dispatching Problem}
\label{sec:kdspp}
\begin{figure}[b]
    \centering
    \begin{subfigure}[t]{1\textwidth}
         \centering
         \def\svgwidth{6.3cm}
         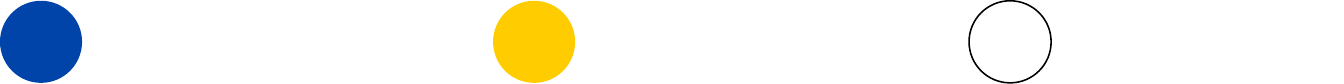
     \end{subfigure}
     \hfill
    \begin{subfigure}[t]{7.5cm}
         \centering
         \def\svgwidth{7.5cm}
         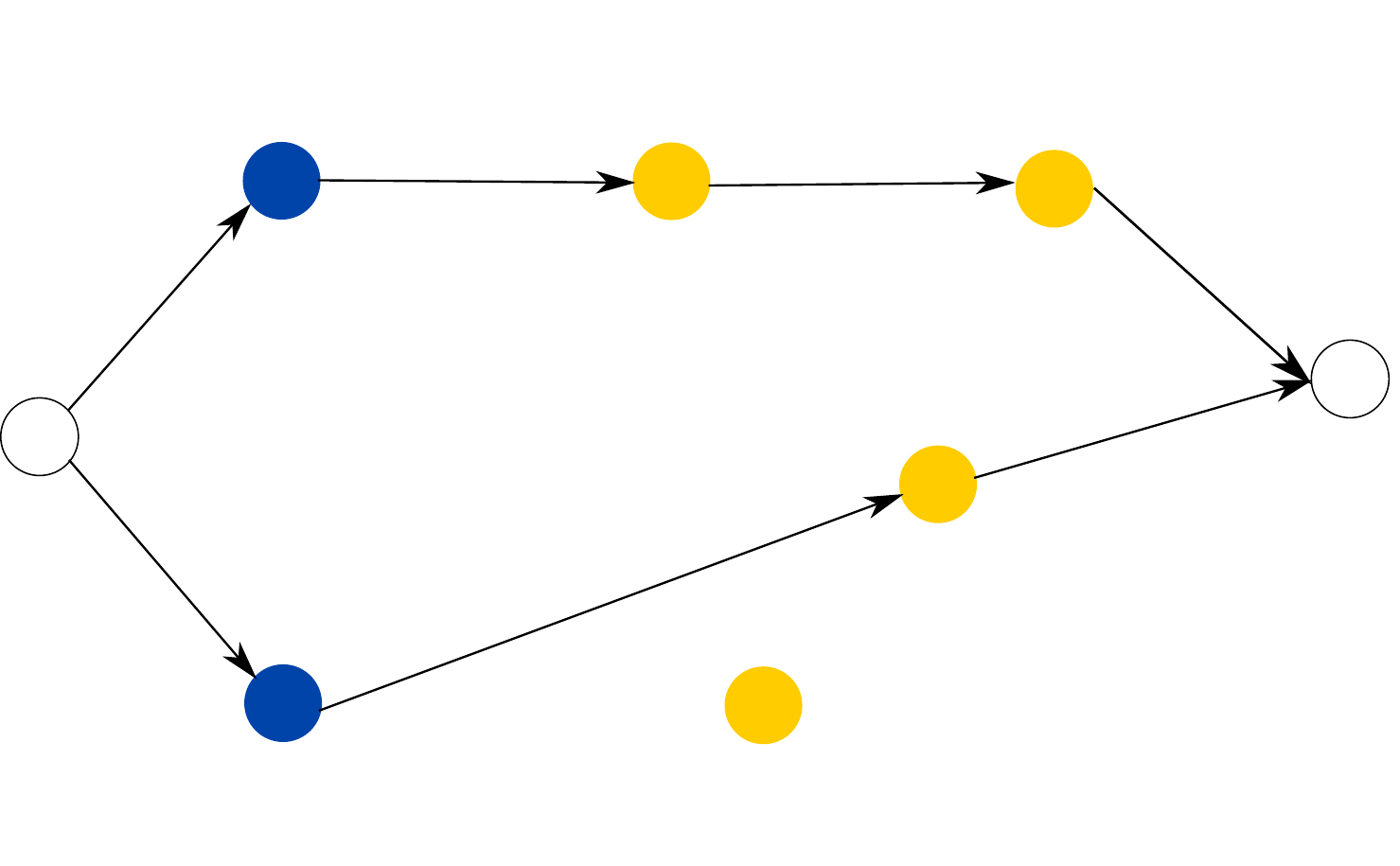
         \caption{Problem instance: digraph~$D$.}
         \label{fig:offlineProblem_Graph}
     \end{subfigure}
     \hfill
     \begin{subfigure}[t]{7.5cm}
         \centering
         \def\svgwidth{7.5cm}
         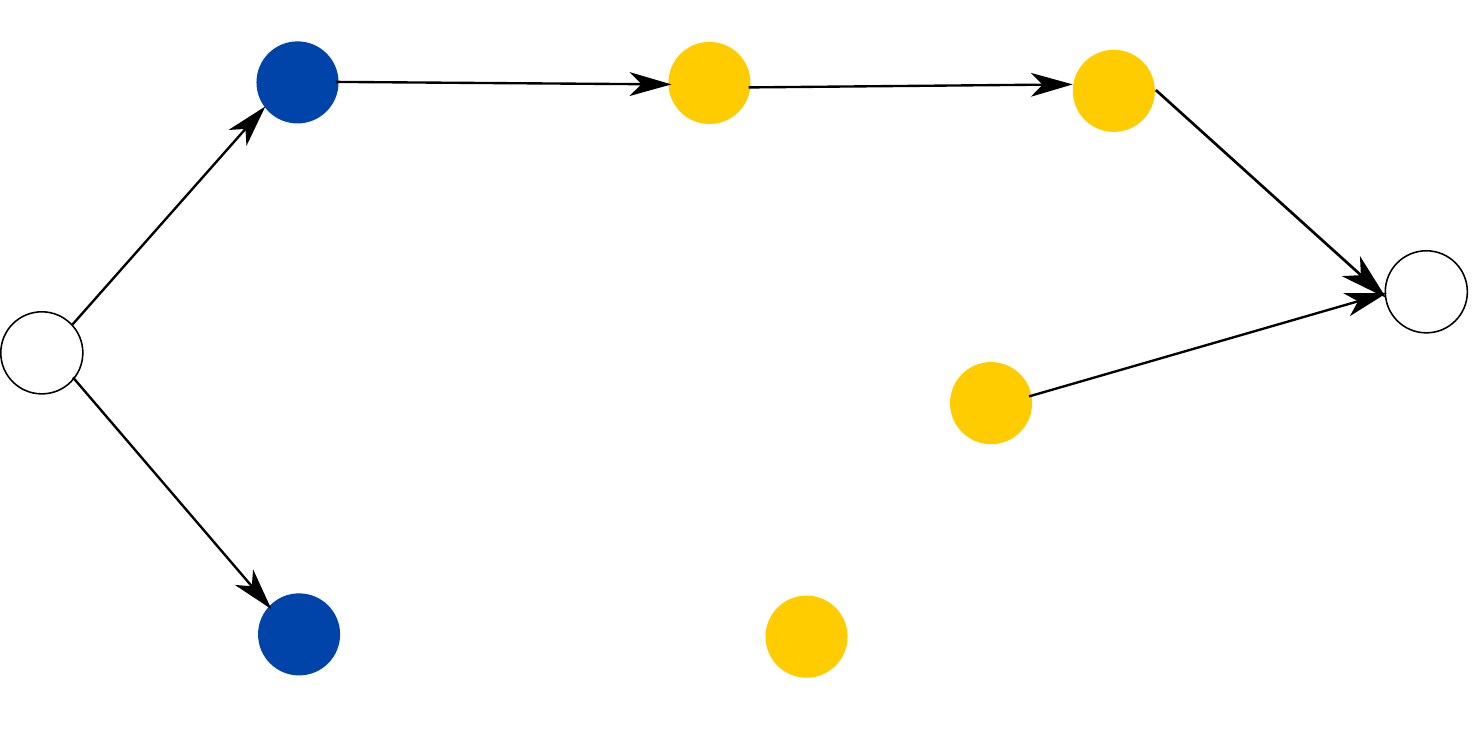
         \caption{Solution:~$k$-disjoint shortest paths.}
         \label{fig:offlineProblem_Solution}
     \end{subfigure}
    \caption{Example of the dispatching problem modeled as a digraph~(a) and a corresponding solution~(b).}
    \label{fig:GraphRepresentation}
\end{figure}

In the following, we show how to formulate our dispatching and rebalancing problem as a \gls*{abk:k-dSPP}, which yields the basis for the policies that we introduce in Section~\ref{sec:SP} and Section~\ref{sec:CP}.
Specifically, we define a dispatching graph as illustrated in Figure~\ref{fig:offlineProblem_Graph}.
We construct our dispatching graph as a weighted and directed graph~$D=(V, A, \theta)$ with a set of vertices~$V$, a set of arcs~$A$, and a vector of weights~$\theta$, which comprises a weight~$\theta_a$ for each arc~$a \in A$. 
The vertex set~$V= V^{\text{v}} \cup V^{\text{r}} \cup V^{\text{a}} \cup \{s^\text{o},s^\text{d}\}$ consists of a subset of vehicle vertices~$V^{\text{v}}$, a subset of request vertices~$V^{\text{r}}$, a subset of artificial vertices~$V^{\text{a}}$, a dummy source vertex~$s^\text{o}$, and a dummy sink vertex~$s^\text{d}$. Here, we associate each request vertex with a customer request, while we associate each vehicle vertex with a vehicle's initial position at the beginning of the system time period or the location where the vehicle will end its previous trip. 
Furthermore, we associate each artificial vertex with a rebalancing action. For now, we assume the set of artificial vertices~$V^{\text{a}}$ to be empty to simplify the example in Figure~\ref{fig:GraphRepresentation}. We then construct the arc set as follows: we create an arc between two request vertices if request-to-request connectivity is possible, i.e., Constraint~(ii) in Section~\ref{sec:ProblemSetting} holds.
Moreover, we create an arc between a vehicle vertex and a request vertex if vehicle-to-request connectivity is possible, i.e., Constraint~(iii) in Section~\ref{sec:ProblemSetting} holds. Finally, we connect the dummy source vertex with each vehicle vertex, and connect all vehicle, request, and artificial vertices to the dummy sink. 

Note that~$D$ is acyclic by construction. Then, a path in~$D$ that starts at the dummy source vertex and ends at the dummy sink vertex represents a feasible vehicle trip. This trip can unambiguously be associated with a certain vehicle as its second vertex is always a vehicle vertex.
To finally devise our dispatching problem, we assign a weight~$\theta_{a}$ to each arc~$a \in A$. We set the weights of arcs that enter the sink or leave the source to zero. Additionally, we set the weights of arcs that enter a request vertex depending on the use case:
in the \gls*{abk:CO}-enriched \gls*{abk:ML} pipeline our \gls*{abk:ML}-predictor~$\varphi_{\theta}$ predicts the weights.

In a full-information setting, we can set the weight of a request's ingoing arc according to its negative reward~$-p_r$. 
We can then find an optimal solution to the \gls*{abk:CO}-layer problem that maximizes the total reward of the central controller by solving a \gls*{abk:k-dSPP} on~$D_{\theta}$, with~$k$ being the number of vehicles of the \gls*{abk:amod} fleet. 
Figure~\ref{fig:offlineProblem_Solution} illustrates an example of a resulting solution in which each disjoint path defines a trip sequence for a specific vehicle. 
The resulting solution is feasible as the disjointness of the paths ensures Constraint~(i), while Constraints~(ii)--(iii) are ensured due to the structure of~$D_{\theta}$. The complexity of solving a \gls*{abk:k-dSPP} on~$D_{\theta}$ is in~$O(|A|(|V| + k) + k|V|\log|V|)$ \citep[cf.][]{Schiffer2021}, such that we obtain a polynomial algorithm to solve the problem of dispatching vehicles to requests. 

Three comments on this algorithm are in order. First, the resulting solution can contain empty trip sequences that indicate inactive vehicles if this benefits the overall reward maximization. In this case, some paths consist of a vehicle vertex being directly connected to the source and sink vertex. Second, one may change the objective of the central controller by modifying the graph's weights, e.g., to maximize the profit of the central controller, we correct a request's reward with arising costs for fulfilling this request. Third, one may further speed up the resulting algorithm by deleting connections of~$D$ to receive a sparse graph, e.g., by connecting only request vertices for which the associated requests are within a certain vicinity.
We discuss such a sparsification heuristic in Appendix~\ref{appendix:Sparsification}.

\subsection{Sample-Based (SB) Pipeline}
\label{sec:SP}
\begin{figure}[b]
    \centering
    \def\svgwidth{15.0cm}
    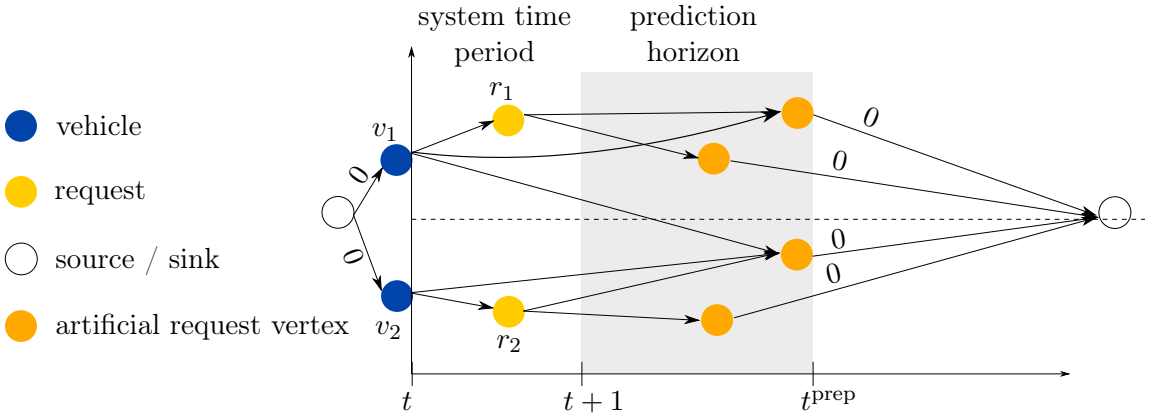
    \caption{Example of an extended digraph with predicted requests in the prediction horizon~$[t+1, t^{\text{prep}})$.}
    \label{fig:onlineGraphRebalancing}
\end{figure}
We recall that the main rationale of the hybrid \gls*{abk:CO}-enriched \gls*{abk:ML} pipeline is to use an \gls*{abk:ML}-predictor to predict the weights of the \gls*{abk:CO}-layer digraph, such that solving the \gls*{abk:k-dSPP} yields an anticipative dispatching and rebalancing solution.
A major drawback of the so far presented \gls*{abk:CO}-layer digraph~$D_{\theta}$ is that it only considers vehicle states~$x_{v,t}$ for all~$v \in V$ and the current request batch~$R_t$ to take a dispatching action in system time period~$[t, t+1$). Accordingly, decisions based on this policy are rebalancing-agnostic, and do not consider future ride requests or any possibility to incorporate rebalancing actions.
The \gls*{abk:SB} pipeline addresses both of these drawbacks by sampling future ride requests and incorporating them into our digraph.
In the following, we detail the components of the \gls*{abk:SB} pipeline.

\paragraph{Digraph Generator:} 
To extend the digraph, we first extend a system state by enlarging its system time period~$[t, t+1)$ with a prediction horizon~$[t+1, t^{\text{pred}})$. Then, the system time period and the prediction horizon form an extended horizon for which we predict a set of artificial requests. To extend the digraph~$D$, we create an artificial vertex~$r\in V^{\text{a}}$ with~$s_r \in [t+1, t^{\text{pred}})$ for each artificial request, and enlarge the arc set accordingly to obtain an extended digraph, see Figure~\ref{fig:onlineGraphRebalancing}.

In this paper, we predict artificial requests by sampling from a request distribution calibrated on historical data, which has proven to yield state-of-the-art results within related online algorithms \citep[cf.][]{Alonso-Mora2017a}. To do so, we partition the fleet's operating area into a rectangular grid, divided into equally-sized squares, and refer to each square as a rebalancing cell. Each rebalancing cell comprises its own request distribution.
We note that one can easily replace this sampling with a more sophisticated predictor. In this case, one could even avoid the usage of rebalancing cells. However, as the size of these cells is usually rather small, the resulting accuracy increase remains limited.

\paragraph{ML-layer:}
Given all the information available, we predict the weight~$\theta_a$ for each arc of the digraph using a generalized linear model,~$\theta_a = \langle w | \phi(a,x) \rangle,$ where~$\phi(a,x)$ is a features vector, which depends on arc~$a$ and system state~$x$. Here, the system state~$x$ can also incorporate historical information. Note that we detail the features used in Appendix~\ref{appendix:Features}.
Intuitively, a weight on an artificial request's ingoing arc represents the learned reward of rebalancing to the location of the artificial request.

\paragraph{CO-layer:}
The \gls*{abk:CO}-layer solves the \gls*{abk:k-dSPP} on the \gls*{abk:CO}-layer-digraph~$D_{\theta}$.

\paragraph{Decoder:}
The decoder~$\mathcal{D}$ receives~$k$~disjoint paths, each representing a vehicle trip, and returns a dispatching and rebalancing action for each vehicle. 
When a vehicle path goes through~$r$, we dispatch the vehicle to~$r$ if~$r$ is an actual request, i.e.,~$s_r \in [t, t+1)$, and we rebalance the vehicle to the location of~$r$ if~$r$ is an artificial request, i.e.,~$s_r \in [t+1, t^{\text{pred}})$.
We refer to the resulting \gls*{abk:SB} policy with~$\delta_{SB}$.

A major drawback of the \gls*{abk:SB} policy presented in Section~\ref{sec:SP}, and state-of-the-art benchmarks \cite[see, e.g.,][]{Alonso-Mora2017a}, is the need for additional distributional information on future ride requests.
To omit the need for this distributional information, we propose a \gls*{abk:CB} pipeline, which aims to extend the digraph with capacity vertices in each rebalancing cell such that it is possible to learn distributional information when learning the weights of the digraph. 
In the following, we detail the components of the \gls*{abk:CB} pipeline.

\subsection{Cell-Based (CB) Pipeline}
\label{sec:CP}
\begin{figure}[b]
    \centering
    \def\svgwidth{16.0cm}
    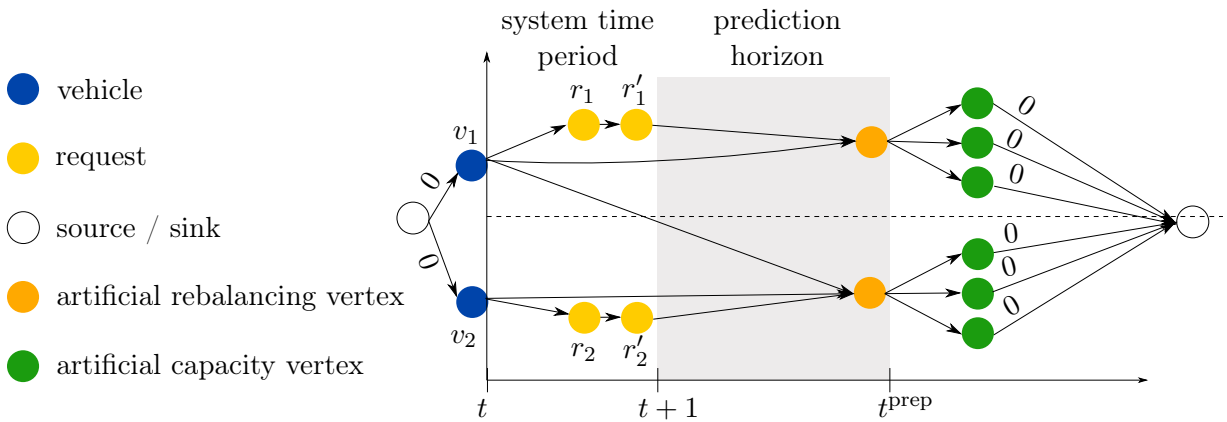
    \caption{Extended digraph with rebalancing vertices in prediction horizon~$[t+1, t^{\text{prep}})$ and capacity vertices.}
    \label{fig:onlineGraphCapacity}
\end{figure}

\paragraph{Digraph Generator:}
We partition the fleet's operating area into a rectangular grid, divided into equally-sized rebalancing cells. We add one rebalancing action as an artificial rebalancing vertex~$r\in V^{\text{a}}$ for each rebalancing cell and connect it to a set of artificial capacity vertices from~$V^{\text{a}}$ as shown in Figure~\ref{fig:onlineGraphCapacity}. 
We refer to the digraph with artificial vertices as the extended digraph.

\paragraph{ML-layer:}
Again, we predict the weights~$\theta_a$ using a generalized linear model~$\theta_a = \langle w | \phi(a,x) \rangle$.
The weight~$\theta_a$ of each rebalancing vertex's ingoing arc reflects the reward of a relocation to the respective rebalancing cell, and the weight of each capacity vertex's ingoing arc reflects the reward of each additional vehicle rebalancing to this cell.

\paragraph{CO-layer:}
Following this pipeline requires us to rethink our \gls*{abk:k-dSPP} approach as we face the following conflict: in the created graph structure, multiple paths need to share a rebalancing vertex to learn the value of multiple rebalancing actions to the same location correctly. On the other hand, we still require paths to be vertex-disjoint with respect to request vertices to ensure dedicated request-to-vehicle assignments. To resolve this contradiction, we adapt the \gls*{abk:k-dSPP} from finding vertex-disjoint paths to finding arc-disjoint paths, which allows paths to share a rebalancing vertex. To ensure disjoint request vertices between paths, we then split each request vertex into two vertices and introduce an arc with weight~$0$ to connect both vertices.
Then we can solve an arc-disjoint \gls*{abk:k-dSPP} on the \gls*{abk:CO}-layer-digraph.
The resulting paths may share rebalancing vertices but are disjunct with respect to request and capacity vertices.

\paragraph{Decoder:}
The decoder~$\mathcal{D}$ receives~$k$ disjoint paths and returns a dispatching and rebalancing action for each vehicle. 
When a vehicle path goes through~$r$, we dispatch the vehicle to~$r$ if~$r$ is an actual request, i.e.,~$s_r \in [t, t+1)$, and we rebalance the vehicle to the location of~$r$ if~$r$ is an artificial rebalancing vertex.
The advantage of the \gls*{abk:CB} pipeline is that it embeds distributional information within the weights of the \gls*{abk:CO}-layer-digraph and does not require any distributional information about future ride requests. 
We refer to this \gls*{abk:CB} policy with~$\delta_{CB}$.

\subsection{Structured Learning methodology}
\label{sec:StructuredLearning}
The policies resulting from the pipelines introduced in Sections~\ref{sec:SP} and~\ref{sec:CP} are clearly sensitive to the weights of the \gls*{abk:CO}-layer-digraph~$D_{\theta}$.
To effectively learn these weights~$\theta$, we leverage \gls*{abk:SL}, i.e., we aim at learning~$\theta$ in a supervised fashion from past full-information solutions. To do so, we must derive a direct relation between the weights~$\theta$ and the resulting policy decisions encoded via solution paths~$y$. Therefore, we introduce the following notation:
we denote a \gls*{abk:CO}-layer-digraph as~$x_i$ and embed the set of all feasible solutions $\mathcal{Y}(x_i)$ in $\mathbb{R}^{|A|}$ so that $y\in \mathcal{Y}(x_i)$ is a vector $(y_a)_{a\in A}$ with
\begin{equation}
    y_a =
\begin{cases}
1, \quad \text{if arc }a\text{ is in }k\text{-disjoint solution paths}\\
0, \quad \text{otherwise.}
\end{cases}
\end{equation}
Then, we can formulate the \gls*{abk:CO}-layer-problem, namely the \gls*{abk:k-dSPP}, as a maximization problem 
\begin{equation}
    \hat y_i (\theta) = \argmax_{y \in \mathcal{Y}(x_i)} \; \theta^T y.
    \label{eq:co_layer_problem}
\end{equation}
Here, our goal is to find a~$\theta$ such that~$\hat y_i(\theta)$ leads to an efficient policy for all~$x_i$ at our disposal.
To do so, 
we use a learning-by-imitation approach.
Practically, we build a training set~$(x_1,y_1^*),\ldots,(x_n,y_n^*)$ containing~$n$ instances~$x_i$ of the CO-layer problem~\eqref{eq:co_layer_problem} as well as their solution~$y_i^*$ that we want to imitate. 
We then formulate the learning problem as,
\begin{equation}
    \min_{\theta} \frac{1}{n}\sum_{i=1}^n L(\theta, x_i, y_i^*),
    \label{eq:learning_problem}
\end{equation}
where the loss function~$L(\theta,x,y^*)$ quantifies the non-optimality when using~$\hat y (\theta)$ instead of~$y^*$.
In the following, we explain how we build~$\{(x_1, y_1^*), ..., (x_n, y_n^*) \}$, which loss~$L$ we use, and how to solve~\eqref{eq:learning_problem}.

\begin{figure}[b]
    \centering
    \begin{subfigure}[t]{1\textwidth}
         \centering
         \def\svgwidth{9.7cm}
         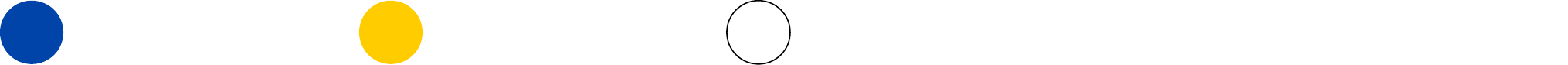
     \end{subfigure}
     \hfill
    \begin{subfigure}[t]{7.5cm}
         \centering
         \def\svgwidth{7.5cm}
         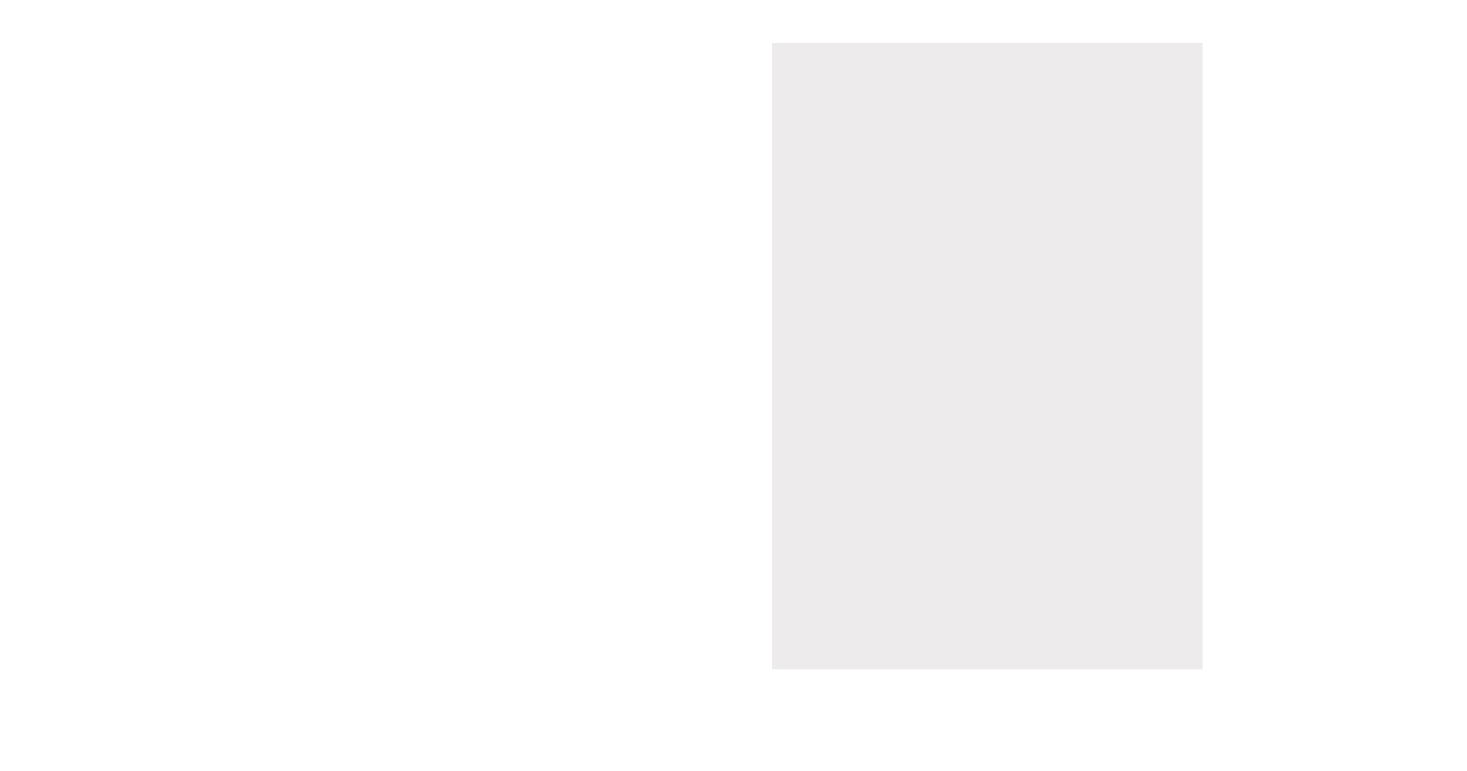
         \caption{Digraph of a system state~$x$.}
         \label{fig:training_instance}
     \end{subfigure}
     \hfill
     \begin{subfigure}[t]{7.5cm}
         \centering
         \def\svgwidth{7.5cm}
         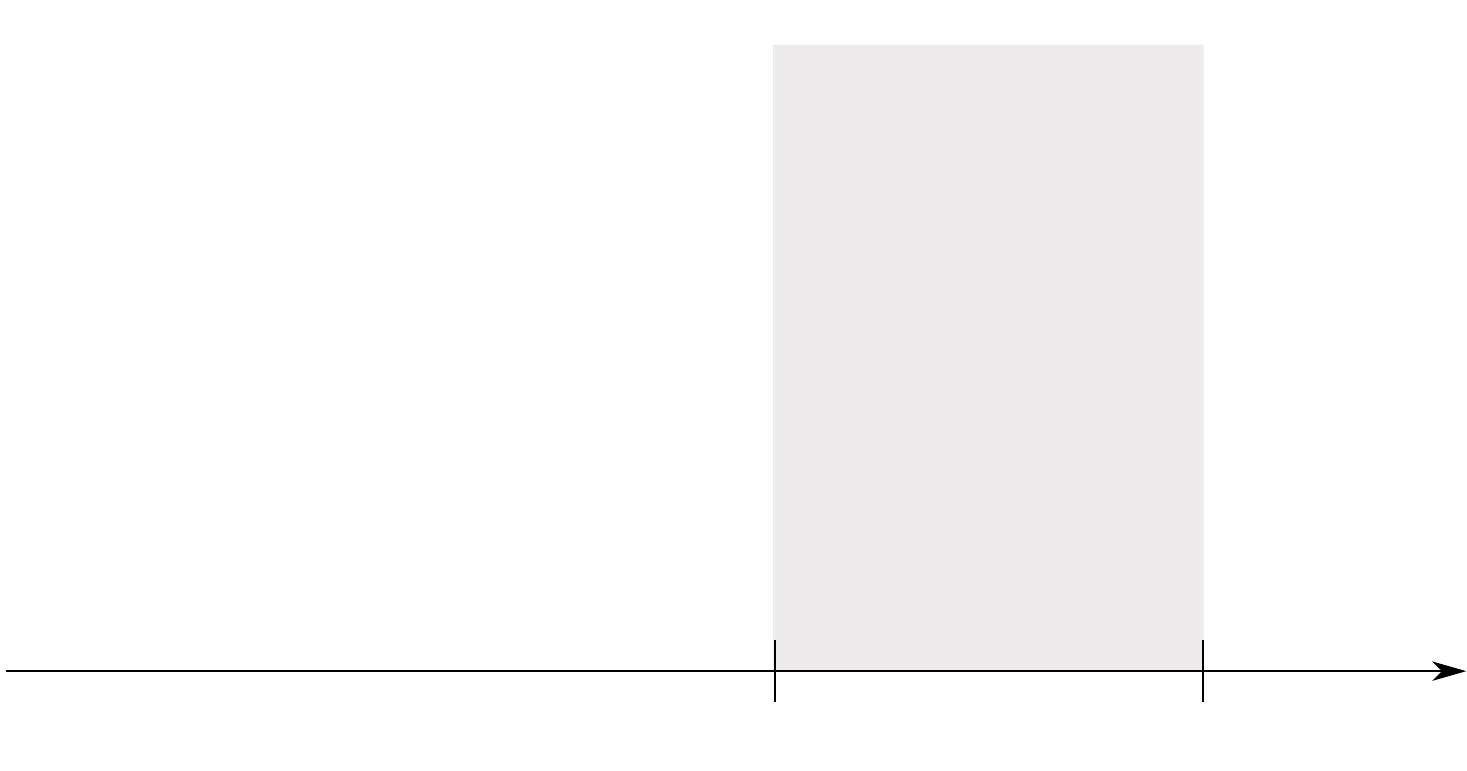
         \caption{Optimal full-information solution~$y^*$.}
         \label{fig:prediction_instance}
     \end{subfigure}
    \caption{Example digraph (left) and the corresponding optimal full-information solution (right). The solid black paths in the optimal full-information solution represent the dispatching solution as~$k$-disjoint shortest paths. The \gls*{abk:SL} approach learns the weights of the digraph such that the solution coincides with the solid black path of the optimal full-information solution.}
    \label{fig:LearningGraph}
\end{figure}

\paragraph{Building the training set: }
We aim to imitate an optimal full-information bound, which cannot be derived during online decision-making because it requires information about future requests.
Accordingly, we build a training set from historical data, which allows us to generate a training set based on full-information solutions.
Our training set contains a training instance~$(x_i, y_i^*)$ for each system time period of every problem instance at our disposal. We generate~$y_i^*$ as follows: we solve the full-information problem of the respective historical problem instance for the complete problem time horizon. To do so, we modify the \gls*{abk:CO}-enriched \gls*{abk:ML} pipeline. In the digraph generator layer, we generate a digraph with all requests of this full-information problem. In the \gls*{abk:ML}-layer, we set the weights of a request's ingoing arc according to the reward which the central controller receives for fulfilling the respective request. Subsequently, we solve the \gls*{abk:CO}-layer-problem~\eqref{eq:co_layer_problem} which retrieves the full-information solution.
From this full-information solution, we rebuild the digraph~$x_i$ of Problem~\eqref{eq:co_layer_problem} for the system time period of interest~(see Figure~\ref{fig:training_instance}) and assign the corresponding full-information solution~$y^*_i$ to it (see Figure~\ref{fig:prediction_instance}).

\paragraph{Fenchel-Young loss:}
Given a training point~$(x_i,y_i^*)$, we want to find a~$\theta$ such that the solution of problem~\eqref{eq:co_layer_problem}, namely~$\hat y_i(\theta)$ equals the optimal solution~$y_i^*$.
It is therefore natural to use as loss~$L(\theta,x_i,y_i^*)$ the non-optimality of~$\hat y_i(\theta) = \argmax_{y \in \mathcal{Y}(x_i)} \; \theta^T y$ as a full-information solution 
\begin{equation}
    L(\theta, x_i, y_i^*) = \underbrace{\max_{y \in \mathcal{Y}(x_i)}(\theta^T y) - \theta^T y_i^*,}_{\text{non-optimiality of~$\hat y_i$ as a full-information solution}}
    \label{eq:loss}
\end{equation}
i.e., we minimize the difference between the solution~$y$ induced by~$\theta$ and the full-information solution~$y_i^*$. Figure~\ref{fig:LearningGraph} visualizes the intuition of this loss function: the loss function minimizes the difference between the sum of weights from the \gls*{abk:CO}-layer solution (left side of~\eqref{eq:loss} and Figure~\ref{fig:training_instance}) and the sum of weights from the full-information solution (right side of~\eqref{eq:loss} and Figure~\ref{fig:prediction_instance}). When minimizing loss function~\eqref{eq:loss}, we aim to adapt the weights~$\theta$ such that the \gls*{abk:CO}-layer solution coincides with the full-information solution.

Unfortunately, the learning problem~\eqref{eq:learning_problem} combined with loss~\eqref{eq:loss} contains a trivial solution with~\mbox{$\theta=0$}. In this case, any feasible solution from~$\mathcal{Y}$ is an optimal solution of the \gls*{abk:CO}-layer-problem~\eqref{eq:co_layer_problem}, which results in a random policy.
We follow~\citet{Berthet2020} and fix this issue by adding a Gaussian perturbation~$Z$ to~$\theta$, which has a regularizing effect,
\begin{equation}
   L(\theta, x_i, y_i^*) = \underbrace{\mathbb{E}\Big[\max_{y \in \mathcal{Y}(x_i)}(\theta + Z)^T y \Big] - \theta^T y_i^*.}_{\text{non-optimality of~$\hat y_i$ as a solution of the perturbed CO-layer problem}}
   \label{eq:loss_perturbed}
\end{equation}
\paragraph{Solving the learning problem:}
An analysis based on Fenchel-Young convex duality enables us to show that this loss is smooth and convex \citep[cf.][]{Berthet2020} with the following gradient
\begin{equation}
    \nabla_{\theta} L(\theta,x_i, y_i^*):= \mathbb{E}\Big[\argmax_{y \in \mathcal{Y}(x_i)}(\theta + Z)^T y \Big] - y_i^*.
    \label{eq:gradient}
\end{equation}
For further discussion of smoothness and convexity properties and the gradient derivation, we refer to Appendix~\ref{appendix:fenchel-young-loss}.

Considering this gradient, we can minimize the learning problem~\eqref{eq:learning_problem}, with loss~\eqref{eq:loss_perturbed} with respect to~$\theta$. As~$\theta$ is a prediction of the \gls*{abk:ML}-predictor~$\varphi_w$, we need to backpropagate the gradient through~$\varphi_w$ to learn~$w$. To do so, the predictor~$\varphi_w$ must be differentiable. A straightforward differentiable predictor is a linear model
\begin{equation}
    \theta_a = \langle w | \phi(a, x_i) \rangle  \quad \forall a \in A,
    \label{eq:MLpredictor}
\end{equation}
and with this linear model,
$w \mapsto L(\varphi_w(x_i),x_i, y_i^*)$ remains convex.
Although we have now defined all items in~\eqref{eq:loss_perturbed} and~\eqref{eq:gradient}, we can still neither calculate the loss nor its gradient as computing the expectation remains intractable.
Therefore, we use a sample average approximation of~\eqref{eq:gradient} to calculate a gradient. With this gradient, we can then minimize a sample average approximation of~\eqref{eq:loss_perturbed} with a \gls*{abk:BFGS} algorithm. The \gls*{abk:BFGS} is a quasi-Newton algorithm that performed well on such problems \citep[cf.][]{Parmentier2021a, Parmentier2022}.

%% file: figures/StructuredLearningParmentier_8.pdf_tex
\begingroup%
  \makeatletter%
  \providecommand\color[2][]{%
    \errmessage{(Inkscape) Color is used for the text in Inkscape, but the package 'color.sty' is not loaded}%
    \renewcommand\color[2][]{}%
  }%
  \providecommand\transparent[1]{%
    \errmessage{(Inkscape) Transparency is used (non-zero) for the text in Inkscape, but the package 'transparent.sty' is not loaded}%
    \renewcommand\transparent[1]{}%
  }%
  \providecommand\rotatebox[2]{#2}%
  \newcommand*\fsize{\dimexpr\f@size pt\relax}%
  \newcommand*\lineheight[1]{\fontsize{\fsize}{#1\fsize}\selectfont}%
  \ifx\svgwidth\undefined%
    \setlength{\unitlength}{863.56149763bp}%
    \ifx\svgscale\undefined%
      \relax%
    \else%
      \setlength{\unitlength}{\unitlength * \real{\svgscale}}%
    \fi%
  \else%
    \setlength{\unitlength}{\svgwidth}%
  \fi%
  \global\let\svgwidth\undefined%
  \global\let\svgscale\undefined%
  \makeatother%
  \begin{picture}(1,0.09028234)%
    \lineheight{1}%
    \setlength\tabcolsep{0pt}%
    \put(0,0){\includegraphics[width=\unitlength,page=1]{StructuredLearningParmentier_8.pdf}}%
    \put(0.38676691,0.04286348){\color[rgb]{0,0,0}\makebox(0,0)[t]{\lineheight{1.25}\smash{\begin{tabular}[t]{c}ML-layer \\$\varphi_{w}$\end{tabular}}}}%
    \put(0.60078248,0.04476927){\color[rgb]{0,0,0}\makebox(0,0)[t]{\lineheight{1.25}\smash{\begin{tabular}[t]{c}CO-layer\\$\mathcal{A}$\end{tabular}}}}%
    \put(0,0){\includegraphics[width=\unitlength,page=2]{StructuredLearningParmentier_8.pdf}}%
    \put(0.49742424,0.08184051){\color[rgb]{0,0,0}\makebox(0,0)[t]{\lineheight{1.25}\smash{\begin{tabular}[t]{c}Arc \\weights\\\\$\theta = \varphi_{w} (D)$\end{tabular}}}}%
    \put(0.69686888,0.05590372){\color[rgb]{0,0,0}\makebox(0,0)[t]{\lineheight{1.25}\smash{\begin{tabular}[t]{c}Solution\\\\$y = \mathcal{A}(D_{\theta})$\end{tabular}}}}%
    \put(0.79647176,0.04542886){\color[rgb]{0,0,0}\makebox(0,0)[t]{\lineheight{1.25}\smash{\begin{tabular}[t]{c}Decoder\\$\mathcal{D}$\end{tabular}}}}%
    \put(0,0){\includegraphics[width=\unitlength,page=3]{StructuredLearningParmentier_8.pdf}}%
    \put(0.92187209,0.04682837){\color[rgb]{0,0,0}\makebox(0,0)[t]{\lineheight{1.25}\smash{\begin{tabular}[t]{c}Action\\$\alpha = \mathcal{D}(y)$\end{tabular}}}}%
    \put(0,0){\includegraphics[width=\unitlength,page=4]{StructuredLearningParmentier_8.pdf}}%
    \put(0.02082751,0.06022609){\color[rgb]{0,0,0}\makebox(0,0)[t]{\lineheight{1.25}\smash{\begin{tabular}[t]{c}System\\state\\ $x$\end{tabular}}}}%
    \put(0.15052353,0.05782333){\color[rgb]{0,0,0}\makebox(0,0)[t]{\lineheight{1.25}\smash{\begin{tabular}[t]{c}Digraph \\Generator\\$\mathcal{G}$\end{tabular}}}}%
    \put(0.26751571,0.05714105){\color[rgb]{0,0,0}\makebox(0,0)[t]{\lineheight{1.25}\smash{\begin{tabular}[t]{c}Digraph\\\\$D = \mathcal{G}(x)$\end{tabular}}}}%
    \put(0,0){\includegraphics[width=\unitlength,page=5]{StructuredLearningParmentier_8.pdf}}%
  \end{picture}%
\endgroup%

%% file: figures/Paper_GraphRepresentationOffline_v2_SubfigureLegend.pdf_tex
\begingroup%
  \makeatletter%
  \providecommand\color[2][]{%
    \errmessage{(Inkscape) Color is used for the text in Inkscape, but the package 'color.sty' is not loaded}%
    \renewcommand\color[2][]{}%
  }%
  \providecommand\transparent[1]{%
    \errmessage{(Inkscape) Transparency is used (non-zero) for the text in Inkscape, but the package 'transparent.sty' is not loaded}%
    \renewcommand\transparent[1]{}%
  }%
  \providecommand\rotatebox[2]{#2}%
  \newcommand*\fsize{\dimexpr\f@size pt\relax}%
  \newcommand*\lineheight[1]{\fontsize{\fsize}{#1\fsize}\selectfont}%
  \ifx\svgwidth\undefined%
    \setlength{\unitlength}{644.72098984bp}%
    \ifx\svgscale\undefined%
      \relax%
    \else%
      \setlength{\unitlength}{\unitlength * \real{\svgscale}}%
    \fi%
  \else%
    \setlength{\unitlength}{\svgwidth}%
  \fi%
  \global\let\svgwidth\undefined%
  \global\let\svgscale\undefined%
  \makeatother%
  \begin{picture}(1,0.06226872)%
    \lineheight{1}%
    \setlength\tabcolsep{0pt}%
    \put(0,0){\includegraphics[width=\unitlength,page=1]{Paper_GraphRepresentationOffline_v2_SubfigureLegend.pdf}}%
    \put(0.08800643,0.01846942){\color[rgb]{0,0,0}\makebox(0,0)[lt]{\lineheight{1.25}\smash{\begin{tabular}[t]{l}vehicle\end{tabular}}}}%
    \put(0.45096203,0.02147612){\color[rgb]{0,0,0}\makebox(0,0)[lt]{\lineheight{1.25}\smash{\begin{tabular}[t]{l}request\end{tabular}}}}%
    \put(0.80449695,0.01968875){\color[rgb]{0,0,0}\makebox(0,0)[lt]{\lineheight{1.25}\smash{\begin{tabular}[t]{l}source / sink\end{tabular}}}}%
  \end{picture}%
\endgroup%

%% file: figures/Paper_GraphRepresentationOffline_v2_SubfigureGraph.pdf_tex
\begingroup%
  \makeatletter%
  \providecommand\color[2][]{%
    \errmessage{(Inkscape) Color is used for the text in Inkscape, but the package 'color.sty' is not loaded}%
    \renewcommand\color[2][]{}%
  }%
  \providecommand\transparent[1]{%
    \errmessage{(Inkscape) Transparency is used (non-zero) for the text in Inkscape, but the package 'transparent.sty' is not loaded}%
    \renewcommand\transparent[1]{}%
  }%
  \providecommand\rotatebox[2]{#2}%
  \newcommand*\fsize{\dimexpr\f@size pt\relax}%
  \newcommand*\lineheight[1]{\fontsize{\fsize}{#1\fsize}\selectfont}%
  \ifx\svgwidth\undefined%
    \setlength{\unitlength}{709.73893734bp}%
    \ifx\svgscale\undefined%
      \relax%
    \else%
      \setlength{\unitlength}{\unitlength * \real{\svgscale}}%
    \fi%
  \else%
    \setlength{\unitlength}{\svgwidth}%
  \fi%
  \global\let\svgwidth\undefined%
  \global\let\svgscale\undefined%
  \makeatother%
  \begin{picture}(1,0.61878108)%
    \lineheight{1}%
    \setlength\tabcolsep{0pt}%
    \put(0,0){\includegraphics[width=\unitlength,page=1]{Paper_GraphRepresentationOffline_v2_SubfigureGraph.pdf}}%
    \put(0.18494183,0.53180388){\color[rgb]{0,0,0}\makebox(0,0)[lt]{\lineheight{1.25}\smash{\begin{tabular}[t]{l}$v_1$\end{tabular}}}}%
    \put(0.46625881,0.53464621){\color[rgb]{0,0,0}\makebox(0,0)[lt]{\lineheight{1.25}\smash{\begin{tabular}[t]{l}$r_1$\end{tabular}}}}%
    \put(0.73063299,0.52926066){\color[rgb]{0,0,0}\makebox(0,0)[lt]{\lineheight{1.25}\smash{\begin{tabular}[t]{l}$r_2$\end{tabular}}}}%
    \put(0.66887252,0.30704373){\color[rgb]{0,0,0}\makebox(0,0)[lt]{\lineheight{1.25}\smash{\begin{tabular}[t]{l}$r_3$\end{tabular}}}}%
    \put(0.52984776,0.05259432){\color[rgb]{0,0,0}\makebox(0,0)[lt]{\lineheight{1.25}\smash{\begin{tabular}[t]{l}$r_5$\end{tabular}}}}%
    \put(0,0){\includegraphics[width=\unitlength,page=2]{Paper_GraphRepresentationOffline_v2_SubfigureGraph.pdf}}%
    \put(0.12949081,0.35342565){\color[rgb]{0,0,0}\rotatebox{-0.31917685}{\makebox(0,0)[lt]{\lineheight{1.25}\smash{\begin{tabular}[t]{l}$0$\end{tabular}}}}}%
    \put(0.12698649,0.23023015){\color[rgb]{0,0,0}\rotatebox{-0.55357223}{\makebox(0,0)[lt]{\lineheight{1.25}\smash{\begin{tabular}[t]{l}$0$\end{tabular}}}}}%
    \put(0.17829996,0.16134674){\color[rgb]{0,0,0}\makebox(0,0)[lt]{\lineheight{1.25}\smash{\begin{tabular}[t]{l}$v_2$\end{tabular}}}}%
    \put(0,0){\includegraphics[width=\unitlength,page=3]{Paper_GraphRepresentationOffline_v2_SubfigureGraph.pdf}}%
    \put(0.85492526,0.42635877){\color[rgb]{0,0,0}\rotatebox{0.25543144}{\makebox(0,0)[lt]{\lineheight{1.25}\smash{\begin{tabular}[t]{l}$0$\end{tabular}}}}}%
    \put(0.7345934,0.30399574){\color[rgb]{0,0,0}\rotatebox{-2.3979291}{\makebox(0,0)[lt]{\lineheight{1.25}\smash{\begin{tabular}[t]{l}$0$\end{tabular}}}}}%
    \put(0.73772856,0.24369674){\color[rgb]{0,0,0}\rotatebox{0.73304266}{\makebox(0,0)[lt]{\lineheight{1.25}\smash{\begin{tabular}[t]{l}$0$\end{tabular}}}}}%
    \put(0.76119454,0.41305877){\color[rgb]{0,0,0}\rotatebox{0.12434524}{\makebox(0,0)[lt]{\lineheight{1.25}\smash{\begin{tabular}[t]{l}$0$\end{tabular}}}}}%
    \put(0.31659412,0.50970814){\color[rgb]{0,0,0}\rotatebox{0.02355137}{\makebox(0,0)[lt]{\lineheight{1.25}\smash{\begin{tabular}[t]{l}$\theta_1$\end{tabular}}}}}%
    \put(0.43790258,0.39837818){\color[rgb]{0,0,0}\rotatebox{-0.17029261}{\makebox(0,0)[lt]{\lineheight{1.25}\smash{\begin{tabular}[t]{l}$\theta_7$\end{tabular}}}}}%
    \put(0.57511211,0.50937908){\color[rgb]{0,0,0}\rotatebox{-0.17356466}{\makebox(0,0)[lt]{\lineheight{1.25}\smash{\begin{tabular}[t]{l}$\theta_8$\end{tabular}}}}}%
    \put(0.43170504,0.22109776){\color[rgb]{0,0,0}\rotatebox{0.30296919}{\makebox(0,0)[lt]{\lineheight{1.25}\smash{\begin{tabular}[t]{l}$\theta_6$\end{tabular}}}}}%
    \put(0.36389389,0.07296672){\color[rgb]{0,0,0}\rotatebox{-0.09055358}{\makebox(0,0)[lt]{\lineheight{1.25}\smash{\begin{tabular}[t]{l}$\theta_4$\end{tabular}}}}}%
    \put(0.5487229,0.16939182){\color[rgb]{0,0,0}\rotatebox{1.100875}{\makebox(0,0)[lt]{\lineheight{1.25}\smash{\begin{tabular}[t]{l}$\theta_5$\end{tabular}}}}}%
    \put(0,0){\includegraphics[width=\unitlength,page=4]{Paper_GraphRepresentationOffline_v2_SubfigureGraph.pdf}}%
    \put(0.55830931,0.36462064){\color[rgb]{0,0,0}\makebox(0,0)[lt]{\lineheight{1.25}\smash{\begin{tabular}[t]{l}$r_4$\end{tabular}}}}%
    \put(0.76750917,0.19617385){\color[rgb]{0,0,0}\makebox(0,0)[lt]{\lineheight{1.25}\smash{\begin{tabular}[t]{l}$r_6$\end{tabular}}}}%
    \put(0,0){\includegraphics[width=\unitlength,page=5]{Paper_GraphRepresentationOffline_v2_SubfigureGraph.pdf}}%
    \put(0.73101316,0.35385164){\color[rgb]{0,0,0}\rotatebox{0.12434521}{\makebox(0,0)[lt]{\lineheight{1.25}\smash{\begin{tabular}[t]{l}$0$\end{tabular}}}}}%
    \put(0.29268217,0.23365663){\color[rgb]{0,0,0}\rotatebox{0.30296919}{\makebox(0,0)[lt]{\lineheight{1.25}\smash{\begin{tabular}[t]{l}$\theta_3$\end{tabular}}}}}%
    \put(0.29419015,0.38415706){\color[rgb]{0,0,0}\rotatebox{0.30296919}{\makebox(0,0)[lt]{\lineheight{1.25}\smash{\begin{tabular}[t]{l}$\theta_2$\end{tabular}}}}}%
    \put(0.67761002,0.09507987){\color[rgb]{0,0,0}\rotatebox{0.30296919}{\makebox(0,0)[lt]{\lineheight{1.25}\smash{\begin{tabular}[t]{l}$\theta_9$\end{tabular}}}}}%
    \put(0,0){\includegraphics[width=\unitlength,page=6]{Paper_GraphRepresentationOffline_v2_SubfigureGraph.pdf}}%
    \put(0.62101493,0.03164944){\color[rgb]{0,0,0}\rotatebox{0.73304272}{\makebox(0,0)[lt]{\lineheight{1.25}\smash{\begin{tabular}[t]{l}$0$\end{tabular}}}}}%
    \put(0,0){\includegraphics[width=\unitlength,page=7]{Paper_GraphRepresentationOffline_v2_SubfigureGraph.pdf}}%
    \put(0.54651595,0.5646659){\color[rgb]{0,0,0}\rotatebox{0.73304272}{\makebox(0,0)[lt]{\lineheight{1.25}\smash{\begin{tabular}[t]{l}$0$\end{tabular}}}}}%
    \put(0.83721428,0.23655711){\color[rgb]{0,0,0}\rotatebox{0.73304272}{\makebox(0,0)[lt]{\lineheight{1.25}\smash{\begin{tabular}[t]{l}$0$\end{tabular}}}}}%
    \put(0.00815505,0.3518336){\color[rgb]{0,0,0}\makebox(0,0)[lt]{\lineheight{1.25}\smash{\begin{tabular}[t]{l}$s^o$\end{tabular}}}}%
    \put(0.95515104,0.39230246){\color[rgb]{0,0,0}\makebox(0,0)[lt]{\lineheight{1.25}\smash{\begin{tabular}[t]{l}$s^d$\end{tabular}}}}%
  \end{picture}%
\endgroup%

%% file: figures/Paper_GraphRepresentationOffline_v2_SubfigureSolution.pdf_tex
\begingroup%
  \makeatletter%
  \providecommand\color[2][]{%
    \errmessage{(Inkscape) Color is used for the text in Inkscape, but the package 'color.sty' is not loaded}%
    \renewcommand\color[2][]{}%
  }%
  \providecommand\transparent[1]{%
    \errmessage{(Inkscape) Transparency is used (non-zero) for the text in Inkscape, but the package 'transparent.sty' is not loaded}%
    \renewcommand\transparent[1]{}%
  }%
  \providecommand\rotatebox[2]{#2}%
  \newcommand*\fsize{\dimexpr\f@size pt\relax}%
  \newcommand*\lineheight[1]{\fontsize{\fsize}{#1\fsize}\selectfont}%
  \ifx\svgwidth\undefined%
    \setlength{\unitlength}{707.98476623bp}%
    \ifx\svgscale\undefined%
      \relax%
    \else%
      \setlength{\unitlength}{\unitlength * \real{\svgscale}}%
    \fi%
  \else%
    \setlength{\unitlength}{\svgwidth}%
  \fi%
  \global\let\svgwidth\undefined%
  \global\let\svgscale\undefined%
  \makeatother%
  \begin{picture}(1,0.49730883)%
    \lineheight{1}%
    \setlength\tabcolsep{0pt}%
    \put(0,0){\includegraphics[width=\unitlength,page=1]{Paper_GraphRepresentationOffline_v2_SubfigureSolution.pdf}}%
    \put(0.18540005,0.48372804){\color[rgb]{0,0,0}\makebox(0,0)[lt]{\lineheight{1.25}\smash{\begin{tabular}[t]{l}$v_1$\end{tabular}}}}%
    \put(0.46741395,0.48657744){\color[rgb]{0,0,0}\makebox(0,0)[lt]{\lineheight{1.25}\smash{\begin{tabular}[t]{l}$r_1$\end{tabular}}}}%
    \put(0.73244328,0.48117855){\color[rgb]{0,0,0}\makebox(0,0)[lt]{\lineheight{1.25}\smash{\begin{tabular}[t]{l}$r_2$\end{tabular}}}}%
    \put(0.6705298,0.258411){\color[rgb]{0,0,0}\makebox(0,0)[lt]{\lineheight{1.25}\smash{\begin{tabular}[t]{l}$r_3$\end{tabular}}}}%
    \put(0.53116046,0.00333114){\color[rgb]{0,0,0}\makebox(0,0)[lt]{\lineheight{1.25}\smash{\begin{tabular}[t]{l}$r_5$\end{tabular}}}}%
    \put(0,0){\includegraphics[width=\unitlength,page=2]{Paper_GraphRepresentationOffline_v2_SubfigureSolution.pdf}}%
    \put(0.17874171,0.11235296){\color[rgb]{0,0,0}\makebox(0,0)[lt]{\lineheight{1.25}\smash{\begin{tabular}[t]{l}$v_2$\end{tabular}}}}%
    \put(0,0){\includegraphics[width=\unitlength,page=3]{Paper_GraphRepresentationOffline_v2_SubfigureSolution.pdf}}%
    \put(0.55969252,0.31613058){\color[rgb]{0,0,0}\makebox(0,0)[lt]{\lineheight{1.25}\smash{\begin{tabular}[t]{l}$r_4$\end{tabular}}}}%
    \put(0.76729211,0.14514771){\color[rgb]{0,0,0}\makebox(0,0)[lt]{\lineheight{1.25}\smash{\begin{tabular}[t]{l}$r_6$\end{tabular}}}}%
    \put(0.00570262,0.29983172){\color[rgb]{0,0,0}\makebox(0,0)[lt]{\lineheight{1.25}\smash{\begin{tabular}[t]{l}$s^o$\end{tabular}}}}%
    \put(0.95503992,0.34040067){\color[rgb]{0,0,0}\makebox(0,0)[lt]{\lineheight{1.25}\smash{\begin{tabular}[t]{l}$s^d$\end{tabular}}}}%
  \end{picture}%
\endgroup%

%% file: figures/Paper_GraphRepresentationOnline2_2.pdf_tex
\begingroup%
  \makeatletter%
  \providecommand\color[2][]{%
    \errmessage{(Inkscape) Color is used for the text in Inkscape, but the package 'color.sty' is not loaded}%
    \renewcommand\color[2][]{}%
  }%
  \providecommand\transparent[1]{%
    \errmessage{(Inkscape) Transparency is used (non-zero) for the text in Inkscape, but the package 'transparent.sty' is not loaded}%
    \renewcommand\transparent[1]{}%
  }%
  \providecommand\rotatebox[2]{#2}%
  \newcommand*\fsize{\dimexpr\f@size pt\relax}%
  \newcommand*\lineheight[1]{\fontsize{\fsize}{#1\fsize}\selectfont}%
  \ifx\svgwidth\undefined%
    \setlength{\unitlength}{1405.13758799bp}%
    \ifx\svgscale\undefined%
      \relax%
    \else%
      \setlength{\unitlength}{\unitlength * \real{\svgscale}}%
    \fi%
  \else%
    \setlength{\unitlength}{\svgwidth}%
  \fi%
  \global\let\svgwidth\undefined%
  \global\let\svgscale\undefined%
  \makeatother%
  \begin{picture}(1,0.36126714)%
    \lineheight{1}%
    \setlength\tabcolsep{0pt}%
    \put(0,0){\includegraphics[width=\unitlength,page=1]{Paper_GraphRepresentationOnline2_2.pdf}}%
    \put(0.32164211,0.24631226){\color[rgb]{0,0,0}\makebox(0,0)[lt]{\lineheight{1.25}\smash{\begin{tabular}[t]{l}$v_1$\end{tabular}}}}%
    \put(0.43110586,0.05948311){\color[rgb]{0,0,0}\makebox(0,0)[lt]{\lineheight{1.25}\smash{\begin{tabular}[t]{l}$r_2$\end{tabular}}}}%
    \put(0,0){\includegraphics[width=\unitlength,page=2]{Paper_GraphRepresentationOnline2_2.pdf}}%
    \put(0.31345123,0.20281531){\color[rgb]{0,0,0}\rotatebox{55.168932}{\makebox(0,0)[lt]{\lineheight{1.25}\smash{\begin{tabular}[t]{l}$0$\end{tabular}}}}}%
    \put(0.29587648,0.14372695){\color[rgb]{0,0,0}\rotatebox{-69.128098}{\makebox(0,0)[lt]{\lineheight{1.25}\smash{\begin{tabular}[t]{l}$0$\end{tabular}}}}}%
    \put(0.32481093,0.07129138){\color[rgb]{0,0,0}\makebox(0,0)[lt]{\lineheight{1.25}\smash{\begin{tabular}[t]{l}$v_2$\end{tabular}}}}%
    \put(0,0){\includegraphics[width=\unitlength,page=3]{Paper_GraphRepresentationOnline2_2.pdf}}%
    \put(0.42478765,0.28255088){\color[rgb]{0,0,0}\makebox(0,0)[lt]{\lineheight{1.25}\smash{\begin{tabular}[t]{l}$r_1$\end{tabular}}}}%
    \put(0.75097035,0.25781274){\color[rgb]{0,0,0}\rotatebox{-22.536366}{\makebox(0,0)[lt]{\lineheight{1.25}\smash{\begin{tabular}[t]{l}$0$\end{tabular}}}}}%
    \put(0.72237319,0.1155131){\color[rgb]{0,0,0}\rotatebox{13.231632}{\makebox(0,0)[lt]{\lineheight{1.25}\smash{\begin{tabular}[t]{l}$0$\end{tabular}}}}}%
    \put(0,0){\includegraphics[width=\unitlength,page=4]{Paper_GraphRepresentationOnline2_2.pdf}}%
    \put(0.34805481,0.00444101){\color[rgb]{0,0,0}\makebox(0,0)[lt]{\lineheight{1.25}\smash{\begin{tabular}[t]{l}$t$\end{tabular}}}}%
    \put(0.49003376,0.00444101){\color[rgb]{0,0,0}\makebox(0,0)[lt]{\lineheight{1.25}\smash{\begin{tabular}[t]{l}$t+1$\end{tabular}}}}%
    \put(0.69819843,0.00444101){\color[rgb]{0,0,0}\makebox(0,0)[lt]{\lineheight{1.25}\smash{\begin{tabular}[t]{l}$t^{\text{prep}}$\end{tabular}}}}%
    \put(0,0){\includegraphics[width=\unitlength,page=5]{Paper_GraphRepresentationOnline2_2.pdf}}%
    \put(0.7244292,0.21840518){\color[rgb]{0,0,0}\rotatebox{-13.581751}{\makebox(0,0)[lt]{\lineheight{1.25}\smash{\begin{tabular}[t]{l}$0$\end{tabular}}}}}%
    \put(0.7257874,0.14711145){\color[rgb]{0,0,0}\rotatebox{8.2973467}{\makebox(0,0)[lt]{\lineheight{1.25}\smash{\begin{tabular}[t]{l}$0$\end{tabular}}}}}%
    \put(0,0){\includegraphics[width=\unitlength,page=6]{Paper_GraphRepresentationOnline2_2.pdf}}%
    \put(0.04480735,0.24768794){\color[rgb]{0,0,0}\makebox(0,0)[lt]{\lineheight{1.25}\smash{\begin{tabular}[t]{l}vehicle\end{tabular}}}}%
    \put(0,0){\includegraphics[width=\unitlength,page=7]{Paper_GraphRepresentationOnline2_2.pdf}}%
    \put(0.04350421,0.19105231){\color[rgb]{0,0,0}\makebox(0,0)[lt]{\lineheight{1.25}\smash{\begin{tabular}[t]{l}request\end{tabular}}}}%
    \put(0,0){\includegraphics[width=\unitlength,page=8]{Paper_GraphRepresentationOnline2_2.pdf}}%
    \put(0.04428608,0.12988704){\color[rgb]{0,0,0}\makebox(0,0)[lt]{\lineheight{1.25}\smash{\begin{tabular}[t]{l}source / sink\end{tabular}}}}%
    \put(0,0){\includegraphics[width=\unitlength,page=9]{Paper_GraphRepresentationOnline2_2.pdf}}%
    \put(0.04416098,0.07179704){\color[rgb]{0,0,0}\makebox(0,0)[lt]{\lineheight{1.25}\smash{\begin{tabular}[t]{l}artificial request vertex\end{tabular}}}}%
    \put(0.42957515,0.34504597){\color[rgb]{0,0,0}\makebox(0,0)[t]{\lineheight{1.25}\smash{\begin{tabular}[t]{c}system time\\period\end{tabular}}}}%
    \put(0.60402806,0.34504597){\color[rgb]{0,0,0}\makebox(0,0)[t]{\lineheight{1.25}\smash{\begin{tabular}[t]{c}prediction\\horizon\end{tabular}}}}%
  \end{picture}%
\endgroup%

%% file: figures/Paper_GraphRepresentationOnline_Capacity.pdf_tex
\begingroup%
  \makeatletter%
  \providecommand\color[2][]{%
    \errmessage{(Inkscape) Color is used for the text in Inkscape, but the package 'color.sty' is not loaded}%
    \renewcommand\color[2][]{}%
  }%
  \providecommand\transparent[1]{%
    \errmessage{(Inkscape) Transparency is used (non-zero) for the text in Inkscape, but the package 'transparent.sty' is not loaded}%
    \renewcommand\transparent[1]{}%
  }%
  \providecommand\rotatebox[2]{#2}%
  \newcommand*\fsize{\dimexpr\f@size pt\relax}%
  \newcommand*\lineheight[1]{\fontsize{\fsize}{#1\fsize}\selectfont}%
  \ifx\svgwidth\undefined%
    \setlength{\unitlength}{1493.31152185bp}%
    \ifx\svgscale\undefined%
      \relax%
    \else%
      \setlength{\unitlength}{\unitlength * \real{\svgscale}}%
    \fi%
  \else%
    \setlength{\unitlength}{\svgwidth}%
  \fi%
  \global\let\svgwidth\undefined%
  \global\let\svgscale\undefined%
  \makeatother%
  \begin{picture}(1,0.34094027)%
    \lineheight{1}%
    \setlength\tabcolsep{0pt}%
    \put(0,0){\includegraphics[width=\unitlength,page=1]{Paper_GraphRepresentationOnline_Capacity.pdf}}%
    \put(0.3667187,0.23377749){\color[rgb]{0,0,0}\makebox(0,0)[lt]{\lineheight{1.25}\smash{\begin{tabular}[t]{l}$v_1$\end{tabular}}}}%
    \put(0.46369225,0.05396194){\color[rgb]{0,0,0}\makebox(0,0)[lt]{\lineheight{1.25}\smash{\begin{tabular}[t]{l}$r_2$\end{tabular}}}}%
    \put(0.3509757,0.18983547){\color[rgb]{0,0,0}\rotatebox{55.168932}{\makebox(0,0)[lt]{\lineheight{1.25}\smash{\begin{tabular}[t]{l}$0$\end{tabular}}}}}%
    \put(0.33745209,0.13524046){\color[rgb]{0,0,0}\rotatebox{-69.128098}{\makebox(0,0)[lt]{\lineheight{1.25}\smash{\begin{tabular}[t]{l}$0$\end{tabular}}}}}%
    \put(0.36467808,0.06507299){\color[rgb]{0,0,0}\makebox(0,0)[lt]{\lineheight{1.25}\smash{\begin{tabular}[t]{l}$v_2$\end{tabular}}}}%
    \put(0,0){\includegraphics[width=\unitlength,page=2]{Paper_GraphRepresentationOnline_Capacity.pdf}}%
    \put(0.46276951,0.26787636){\color[rgb]{0,0,0}\makebox(0,0)[lt]{\lineheight{1.25}\smash{\begin{tabular}[t]{l}$r_1$\end{tabular}}}}%
    \put(0.82973792,0.25933902){\color[rgb]{0,0,0}\rotatebox{-34.494391}{\makebox(0,0)[lt]{\lineheight{1.25}\smash{\begin{tabular}[t]{l}$0$\end{tabular}}}}}%
    \put(0.82469271,0.08661073){\color[rgb]{0,0,0}\rotatebox{30.839152}{\makebox(0,0)[lt]{\lineheight{1.25}\smash{\begin{tabular}[t]{l}$0$\end{tabular}}}}}%
    \put(0,0){\includegraphics[width=\unitlength,page=3]{Paper_GraphRepresentationOnline_Capacity.pdf}}%
    \put(0.38654951,0.00417879){\color[rgb]{0,0,0}\makebox(0,0)[lt]{\lineheight{1.25}\smash{\begin{tabular}[t]{l}$t$\end{tabular}}}}%
    \put(0.5191407,0.00417879){\color[rgb]{0,0,0}\makebox(0,0)[lt]{\lineheight{1.25}\smash{\begin{tabular}[t]{l}$t+1$\end{tabular}}}}%
    \put(0.71601858,0.00417879){\color[rgb]{0,0,0}\makebox(0,0)[lt]{\lineheight{1.25}\smash{\begin{tabular}[t]{l}$t^{\text{prep}}$\end{tabular}}}}%
    \put(0,0){\includegraphics[width=\unitlength,page=4]{Paper_GraphRepresentationOnline_Capacity.pdf}}%
    \put(0.81991987,0.23038336){\color[rgb]{0,0,0}\rotatebox{-26.006494}{\makebox(0,0)[lt]{\lineheight{1.25}\smash{\begin{tabular}[t]{l}$0$\end{tabular}}}}}%
    \put(0.8205473,0.14765962){\color[rgb]{0,0,0}\rotatebox{6.2109191}{\makebox(0,0)[lt]{\lineheight{1.25}\smash{\begin{tabular}[t]{l}$0$\end{tabular}}}}}%
    \put(0,0){\includegraphics[width=\unitlength,page=5]{Paper_GraphRepresentationOnline_Capacity.pdf}}%
    \put(0.04122888,0.26726471){\color[rgb]{0,0,0}\makebox(0,0)[lt]{\lineheight{1.25}\smash{\begin{tabular}[t]{l}vehicle\end{tabular}}}}%
    \put(0,0){\includegraphics[width=\unitlength,page=6]{Paper_GraphRepresentationOnline_Capacity.pdf}}%
    \put(0.04000269,0.21237669){\color[rgb]{0,0,0}\makebox(0,0)[lt]{\lineheight{1.25}\smash{\begin{tabular}[t]{l}request\end{tabular}}}}%
    \put(0,0){\includegraphics[width=\unitlength,page=7]{Paper_GraphRepresentationOnline_Capacity.pdf}}%
    \put(0.04073839,0.15322648){\color[rgb]{0,0,0}\makebox(0,0)[lt]{\lineheight{1.25}\smash{\begin{tabular}[t]{l}source / sink\end{tabular}}}}%
    \put(0,0){\includegraphics[width=\unitlength,page=8]{Paper_GraphRepresentationOnline_Capacity.pdf}}%
    \put(0.04062068,0.09697009){\color[rgb]{0,0,0}\makebox(0,0)[lt]{\lineheight{1.25}\smash{\begin{tabular}[t]{l}artificial rebalancing vertex\end{tabular}}}}%
    \put(0,0){\includegraphics[width=\unitlength,page=9]{Paper_GraphRepresentationOnline_Capacity.pdf}}%
    \put(0.04062068,0.03955614){\color[rgb]{0,0,0}\makebox(0,0)[lt]{\lineheight{1.25}\smash{\begin{tabular}[t]{l}artificial capacity vertex\end{tabular}}}}%
    \put(0,0){\includegraphics[width=\unitlength,page=10]{Paper_GraphRepresentationOnline_Capacity.pdf}}%
    \put(0.82278865,0.1992997){\color[rgb]{0,0,0}\rotatebox{-14.900643}{\makebox(0,0)[lt]{\lineheight{1.25}\smash{\begin{tabular}[t]{l}$0$\end{tabular}}}}}%
    \put(0,0){\includegraphics[width=\unitlength,page=11]{Paper_GraphRepresentationOnline_Capacity.pdf}}%
    \put(0.82016709,0.11952978){\color[rgb]{0,0,0}\rotatebox{14.346759}{\makebox(0,0)[lt]{\lineheight{1.25}\smash{\begin{tabular}[t]{l}$0$\end{tabular}}}}}%
    \put(0,0){\includegraphics[width=\unitlength,page=12]{Paper_GraphRepresentationOnline_Capacity.pdf}}%
    \put(0.46946819,0.3256769){\color[rgb]{0,0,0}\makebox(0,0)[t]{\lineheight{1.25}\smash{\begin{tabular}[t]{c}system time\\period\end{tabular}}}}%
    \put(0.63362035,0.3256769){\color[rgb]{0,0,0}\makebox(0,0)[t]{\lineheight{1.25}\smash{\begin{tabular}[t]{c}prediction\\horizon\end{tabular}}}}%
    \put(0.50574683,0.2690967){\color[rgb]{0,0,0}\makebox(0,0)[lt]{\lineheight{1.25}\smash{\begin{tabular}[t]{l}$r'_1$\end{tabular}}}}%
    \put(0.50646465,0.05285671){\color[rgb]{0,0,0}\makebox(0,0)[lt]{\lineheight{1.25}\smash{\begin{tabular}[t]{l}$r'_2$\end{tabular}}}}%
    \put(0,0){\includegraphics[width=\unitlength,page=13]{Paper_GraphRepresentationOnline_Capacity.pdf}}%
  \end{picture}%
\endgroup%

%% file: figures/Learning2_SubfigureLegend.pdf_tex
\begingroup%
  \makeatletter%
  \providecommand\color[2][]{%
    \errmessage{(Inkscape) Color is used for the text in Inkscape, but the package 'color.sty' is not loaded}%
    \renewcommand\color[2][]{}%
  }%
  \providecommand\transparent[1]{%
    \errmessage{(Inkscape) Transparency is used (non-zero) for the text in Inkscape, but the package 'transparent.sty' is not loaded}%
    \renewcommand\transparent[1]{}%
  }%
  \providecommand\rotatebox[2]{#2}%
  \newcommand*\fsize{\dimexpr\f@size pt\relax}%
  \newcommand*\lineheight[1]{\fontsize{\fsize}{#1\fsize}\selectfont}%
  \ifx\svgwidth\undefined%
    \setlength{\unitlength}{974.34498325bp}%
    \ifx\svgscale\undefined%
      \relax%
    \else%
      \setlength{\unitlength}{\unitlength * \real{\svgscale}}%
    \fi%
  \else%
    \setlength{\unitlength}{\svgwidth}%
  \fi%
  \global\let\svgwidth\undefined%
  \global\let\svgscale\undefined%
  \makeatother%
  \begin{picture}(1,0.04120301)%
    \lineheight{1}%
    \setlength\tabcolsep{0pt}%
    \put(0,0){\includegraphics[width=\unitlength,page=1]{Learning2_SubfigureLegend.pdf}}%
    \put(0.05823358,0.0069325){\color[rgb]{0,0,0}\makebox(0,0)[lt]{\lineheight{1.25}\smash{\begin{tabular}[t]{l}vehicle\end{tabular}}}}%
    \put(0.28454472,0.01091154){\color[rgb]{0,0,0}\makebox(0,0)[lt]{\lineheight{1.25}\smash{\begin{tabular}[t]{l}request\end{tabular}}}}%
    \put(0.51847781,0.00854613){\color[rgb]{0,0,0}\makebox(0,0)[lt]{\lineheight{1.25}\smash{\begin{tabular}[t]{l}source / sink\end{tabular}}}}%
    \put(0,0){\includegraphics[width=\unitlength,page=2]{Learning2_SubfigureLegend.pdf}}%
    \put(0.84549918,0.01091154){\color[rgb]{0,0,0}\makebox(0,0)[lt]{\lineheight{1.25}\smash{\begin{tabular}[t]{l}artificial vertex\end{tabular}}}}%
  \end{picture}%
\endgroup%

%% file: figures/Learning2_Prediction.pdf_tex
\begingroup%
  \makeatletter%
  \providecommand\color[2][]{%
    \errmessage{(Inkscape) Color is used for the text in Inkscape, but the package 'color.sty' is not loaded}%
    \renewcommand\color[2][]{}%
  }%
  \providecommand\transparent[1]{%
    \errmessage{(Inkscape) Transparency is used (non-zero) for the text in Inkscape, but the package 'transparent.sty' is not loaded}%
    \renewcommand\transparent[1]{}%
  }%
  \providecommand\rotatebox[2]{#2}%
  \newcommand*\fsize{\dimexpr\f@size pt\relax}%
  \newcommand*\lineheight[1]{\fontsize{\fsize}{#1\fsize}\selectfont}%
  \ifx\svgwidth\undefined%
    \setlength{\unitlength}{703.80826188bp}%
    \ifx\svgscale\undefined%
      \relax%
    \else%
      \setlength{\unitlength}{\unitlength * \real{\svgscale}}%
    \fi%
  \else%
    \setlength{\unitlength}{\svgwidth}%
  \fi%
  \global\let\svgwidth\undefined%
  \global\let\svgscale\undefined%
  \makeatother%
  \begin{picture}(1,0.52523034)%
    \lineheight{1}%
    \setlength\tabcolsep{0pt}%
    \put(0,0){\includegraphics[width=\unitlength,page=1]{Learning2_Prediction.pdf}}%
    \put(0.27263122,0.45079632){\color[rgb]{0,0,0}\rotatebox{1.6968753}{\makebox(0,0)[lt]{\lineheight{1.25}\smash{\begin{tabular}[t]{l}$\theta_1$\end{tabular}}}}}%
    \put(0.28564718,0.34019532){\color[rgb]{0,0,0}\rotatebox{0.04262907}{\makebox(0,0)[lt]{\lineheight{1.25}\smash{\begin{tabular}[t]{l}$\theta_2$\end{tabular}}}}}%
    \put(0.37012741,0.36857046){\color[rgb]{0,0,0}\rotatebox{0.51862491}{\makebox(0,0)[lt]{\lineheight{1.25}\smash{\begin{tabular}[t]{l}$\theta_8$\end{tabular}}}}}%
    \put(0.27273873,0.26593211){\color[rgb]{0,0,0}\rotatebox{1.400879}{\makebox(0,0)[lt]{\lineheight{1.25}\smash{\begin{tabular}[t]{l}$\theta_3$\end{tabular}}}}}%
    \put(0.57284415,0.2297725){\color[rgb]{0,0,0}\rotatebox{0.08385081}{\makebox(0,0)[lt]{\lineheight{1.25}\smash{\begin{tabular}[t]{l}$\theta_6$\end{tabular}}}}}%
    \put(0.57083455,0.33423633){\color[rgb]{0,0,0}\rotatebox{-0.08206237}{\makebox(0,0)[lt]{\lineheight{1.25}\smash{\begin{tabular}[t]{l}$\theta_7$\end{tabular}}}}}%
    \put(0.46295499,0.17447961){\color[rgb]{0,0,0}\rotatebox{3.0226965}{\makebox(0,0)[lt]{\lineheight{1.25}\smash{\begin{tabular}[t]{l}$\theta_5$\end{tabular}}}}}%
    \put(0.26048757,0.14184553){\color[rgb]{0,0,0}\rotatebox{-1.107518}{\makebox(0,0)[lt]{\lineheight{1.25}\smash{\begin{tabular}[t]{l}$\theta_4$\end{tabular}}}}}%
    \put(0.46763823,0.38725995){\color[rgb]{0,0,0}\rotatebox{-0.31243466}{\makebox(0,0)[lt]{\lineheight{1.25}\smash{\begin{tabular}[t]{l}$\theta_9$\end{tabular}}}}}%
    \put(0,0){\includegraphics[width=\unitlength,page=2]{Learning2_Prediction.pdf}}%
    \put(0.22523055,0.00626475){\color[rgb]{0,0,0}\makebox(0,0)[lt]{\lineheight{1.25}\smash{\begin{tabular}[t]{l}$t$\end{tabular}}}}%
    \put(0.48356688,0.00839598){\color[rgb]{0,0,0}\makebox(0,0)[lt]{\lineheight{1.25}\smash{\begin{tabular}[t]{l}$t+1$\end{tabular}}}}%
    \put(0,0){\includegraphics[width=\unitlength,page=3]{Learning2_Prediction.pdf}}%
    \put(0.80567631,0.00929232){\color[rgb]{0,0,0}\makebox(0,0)[lt]{\lineheight{1.25}\smash{\begin{tabular}[t]{l}$t^{\text{prep}}$\end{tabular}}}}%
    \put(0,0){\includegraphics[width=\unitlength,page=4]{Learning2_Prediction.pdf}}%
  \end{picture}%
\endgroup%

%% file: figures/Learning2_Training.pdf_tex
\begingroup%
  \makeatletter%
  \providecommand\color[2][]{%
    \errmessage{(Inkscape) Color is used for the text in Inkscape, but the package 'color.sty' is not loaded}%
    \renewcommand\color[2][]{}%
  }%
  \providecommand\transparent[1]{%
    \errmessage{(Inkscape) Transparency is used (non-zero) for the text in Inkscape, but the package 'transparent.sty' is not loaded}%
    \renewcommand\transparent[1]{}%
  }%
  \providecommand\rotatebox[2]{#2}%
  \newcommand*\fsize{\dimexpr\f@size pt\relax}%
  \newcommand*\lineheight[1]{\fontsize{\fsize}{#1\fsize}\selectfont}%
  \ifx\svgwidth\undefined%
    \setlength{\unitlength}{703.80825789bp}%
    \ifx\svgscale\undefined%
      \relax%
    \else%
      \setlength{\unitlength}{\unitlength * \real{\svgscale}}%
    \fi%
  \else%
    \setlength{\unitlength}{\svgwidth}%
  \fi%
  \global\let\svgwidth\undefined%
  \global\let\svgscale\undefined%
  \makeatother%
  \begin{picture}(1,0.52523032)%
    \lineheight{1}%
    \setlength\tabcolsep{0pt}%
    \put(0,0){\includegraphics[width=\unitlength,page=1]{Learning2_Training.pdf}}%
    \put(0.22523055,0.00626475){\color[rgb]{0,0,0}\makebox(0,0)[lt]{\lineheight{1.25}\smash{\begin{tabular}[t]{l}$t$\end{tabular}}}}%
    \put(0.48356689,0.00839598){\color[rgb]{0,0,0}\makebox(0,0)[lt]{\lineheight{1.25}\smash{\begin{tabular}[t]{l}$t+1$\end{tabular}}}}%
    \put(0,0){\includegraphics[width=\unitlength,page=2]{Learning2_Training.pdf}}%
    \put(0.80567631,0.00929232){\color[rgb]{0,0,0}\makebox(0,0)[lt]{\lineheight{1.25}\smash{\begin{tabular}[t]{l}$t^{\text{prep}}$\end{tabular}}}}%
    \put(0,0){\includegraphics[width=\unitlength,page=3]{Learning2_Training.pdf}}%
  \end{picture}%
\endgroup%

%% file: contents/experimentaldesign.tex
\section{Experimental Design}
\label{sec:computationalExperiments}

This section details our experimental design by first introducing our real-world case study in Section~\ref{sec:casestudy} and afterward detailing our benchmark policies in Section~\ref{sec:benchmarks}.
\subsection{Case Study}
\label{sec:casestudy}
Our case study focuses on an online \gls*{abk:amod} control problem for the New York City Taxi data set \citep{TLC} from the year~2015.
We chose the year 2015 as it is the most recent year with exact recordings of starting and destination locations from ride requests, and restrict the simulation area to Manhattan.
We model vehicle movements with exact longitude and latitude values. To avoid a computational overload for on-the-fly computation of driving times, we use a look-up table to retrieve the distance and the time needed to drive from a starting location to a destination location. We base this look-up table on data from Uber Movement for the year 2020 (Data retrieved from Uber Movement, (c) 2022 Uber Technologies, Inc., https://movement.uber.com.), which accounts for traffic conditions accordingly.
Each location in this look-up table is specified by a rather small square cell of size \lookUpX~meters~$\times$ \lookUpY~meters.
If the starting and destination locations are in the same cell, we assume a line distance and an average driving speed of \averageDrivingSpeed~kilometers per hour.
The rebalancing cells are specified by larger square cells of size \rebalancingX~meters~$\times$ \rebalancingY~meters to represent rebalancing areas.

We only consider weekdays as we assume different travel behavior on weekends, and focus on an \gls*{abk:amod} system operating for one hour from~9~AM till~10~AM. 
We assume a system time period of~\systemTimePeriod~minute such that the central controller can apply dispatching and rebalancing actions every minute.
We split data into disjoint sets of training data, validation data, and testing data. The training data and the validation data comprise 5 working days~$\times$~1~hour of request data, and the testing data comprises~20 working days~$\times$~1~hour of request data. We stop the training after 70 training iterations or a maximum of 70~hours of training time. When discussing results, we report the averaged results over the testing data.

To further speed up the computation times of our \gls*{abk:k-dSPP} algorithm, we apply sparsification techniques to sparsify the \gls*{abk:CO}-layer-digraph. For a detailed description of this sparsification, we refer to Appendix~\ref{appendix:Sparsification}. For our specific case study, we finally only include an arc in the \gls*{abk:CO}-layer-digraph if the distance to the next request is below~\distanceSparsification~kilometers, and if the next request is reachable within~\temporalSparsification~seconds. These thresholds permit to significantly sparsen our digraph but have hardly any effect on the solution quality. Note that these thresholds are problem and instance specific and should be recalibrated for other use cases. 

We generate different scenarios with varying fleet sizes ranging from~$500$ to~$5000$ vehicles, and scenarios with different ride request densities, ranging from~$10\%$ to~$100\%$ of the original ride request volume. In each scenario, we optimize over profit and the number of satisfied customers. When we optimize over profit, we take the revenue for fulfilling the request extracted from \citep{TLC} minus operational costs of 0.45\$ per driven kilometer \citep{Bosch2018}. When we optimize over the number of satisfied customers, we take a revenue of 1 for fulfilling a ride request minus marginal costs of 0.00001 per driven kilometer, to retrieve a lexicographic objective that maximizes the number of satisfied ride requests but prefers sustainable solutions with fewer kilometers driven over less sustainable solutions. 
If not specified differently, we assume a base scenario that considers~$30\%$ of the original ride request volume and a vehicle fleet size of~$2000$ vehicles. Here, we chose a request share of~$30\%$ to reflect a typical market share of a ride-hailing fleet \citep[cf.][]{Statista}. We then set the vehicle fleet size such that it enables a~$100\%$ request fulfillment under an optimal full-information approach.

\paragraph{Discussion}
A few comments on the setup of this case study are in order.
First, we decided to evaluate the system on a one-hour period for several reasons: Learning a policy that operates longer than a one-hour period requires solving an offline dispatching problem which is computationally hard for large-scale scenarios. However, operating the system longer than a one-hour period is easily possible by operating individually trained policies in each hour. As the performance of our algorithms is similar over different hours of a day, we only report this one-hour period for better comparison.
Second, we experience that our \gls*{abk:CO}-enriched \gls*{abk:ML} pipeline leads to very general policies, shown by superior average performance over weekdays, although the request distribution might change from Monday to Friday. In this context, one may argue that demand patterns can change more drastically between different seasons e.g., from winter to summer, or weekday to weekend. For such cases, our analyses are still valid as one can easily either train separate models for drastically different scenarios or add context to the predictor to enhance its generalizability further. Third, we limit our studies to one case study and train our model for this respective case. This is common practice, as generalizing across different case studies can be an unnecessary tedious task. In practice as well as in research there is some consensus that applying a learning-based model to a different ground truth is better handled by retraining, which is also the case for our model: to apply our pipeline to a different city, it needs to be retrained.

\subsection{Benchmark policies}
\label{sec:benchmarks}
In our computational studies, we apply the following benchmark policies.
\begin{description}
\item[\bf Greedy Policy:] The greedy benchmark policy is an online policy that greedily re-optimizes dispatching decisions in a rolling horizon fashion. The greedy policy does not include any rebalancing vertices into the digraph and, instead of learning the weights of arcs in the digraph, sets the weight of a request's ingoing arc according to its reward. Solving a \gls{abk:k-dSPP} on this digraph leads to greedy dispatching actions which only incorporate requests that are within the actual system time period. 

\item[\bf Sampling Policy:] The sampling benchmark policy is an online policy that re-optimizes dispatching and rebalancing decisions in a rolling horizon fashion and samples future requests for the prediction horizon from historical data. We can interpret the sampling policy also as a variant of an \gls*{abk:CO}-enriched \gls*{abk:ML} pipeline with a digraph generator that incorporates sampled ride requests as artificial rebalancing vertices in the digraph. Instead of learning the weights of arcs in the digraph, it sets the weight of a request's ingoing arc according to its reward. To prefer real requests from the system time period over sampled requests from the prediction horizon, we multiply the reward of sampled ride requests with a discount factor of~\reductionFactorSampling~as we received the best results for this factor. We consider a prediction horizon of~\predictionHorizonSampling~minutes as it leads to best results (see Appendix~\ref{appendix:hyper_samplingBench}). 

\item[\bf Full-information bound (FI):] We compute a full-information bound with complete information about all ride requests. Solving the full-information problem optimally constitutes an upper bound for the online problem. To compute the full-information solution, we assume a system time period equal to the problem time horizon~$[1, T)$, such that all ride requests~$\bigcup_t R_t$ are available at the system at time~$t=1$.
Accordingly, the digraph at~$t=1$ comprises all requests of the complete problem time horizon. Therefore, we do not have to predict the arc weights in the digraph but can set the arc weights of an arc equal to the reward for the ingoing request vertex.

\item[\bf Sample-Based (\gls*{abk:SB}) Policy:] The \gls*{abk:SB} policy is an online policy that follows the variant of the \gls*{abk:CO}-enriched \gls*{abk:ML} pipeline presented in Section~\ref{sec:SP}. The policy re-optimizes dispatching and rebalancing actions in a rolling horizon fashion. For the \gls*{abk:SB} policy, we consider a prediction horizon of~\predictionHorizonSB~minutes, as we received the best results for this planning horizon (see Appendix~\ref{appendix:hyper_idea}).

\item[\bf Cell-Based (\gls*{abk:CB}) Policy:] The \gls*{abk:CB} policy is an online policy that follows the variant of the \gls*{abk:CO}-enriched \gls*{abk:ML} pipeline presented in Section~\ref{sec:CP}. The policy re-optimizes dispatching and rebalancing actions in a rolling horizon fashion. We set the number of capacity vertices for every rebalancing cell to~\maxCapacityCB, as we achieved the best results with this setting. For the \gls*{abk:CB} policy, we consider a prediction horizon of~\predictionHorizonCB~minutes, as we received the best results for this planning horizon (see Appendix~\ref{appendix:hyper_idea}).
\end{description}

One comment on the selection of our benchmark policies is in order. One may wonder why we do not explicitly consider any \gls*{abk:DRL}-based policy. In fact, all publications that present a \gls*{abk:DRL}-based policy benchmark this policy solely against a greedy baseline algorithm. Accordingly, by benchmarking our learning-based policies against the greedy benchmark suffices to discuss the performance of our policies against \gls*{abk:DRL}-based approaches and keeps the results section as concise as possible.

%% file: contents/results.tex
\section{Results}
\label{sec:results}
In the following, we first present the results of the real-world case study for different vehicle fleet sizes and request densities in Section~\ref{sec:AggregatedAnalysis}. We then provide an in-depth structural analysis for the \gls*{abk:SB} and \gls*{abk:CB} policies and compare them with the other benchmark policies with respect to external effects and computational time in Section~\ref{sec:AppliedDiscussion}.
\subsection{Performance Analysis}
\label{sec:AggregatedAnalysis}
To analyze the performance of our learning-based policies, we benchmark their performance for various scenarios with varying fleet sizes and request densities as well as for varying objectives, first focusing on profit maximization.
\paragraph{Varying fleet size:}
\begingroup
\renewcommand\baselinestretch{0.5}
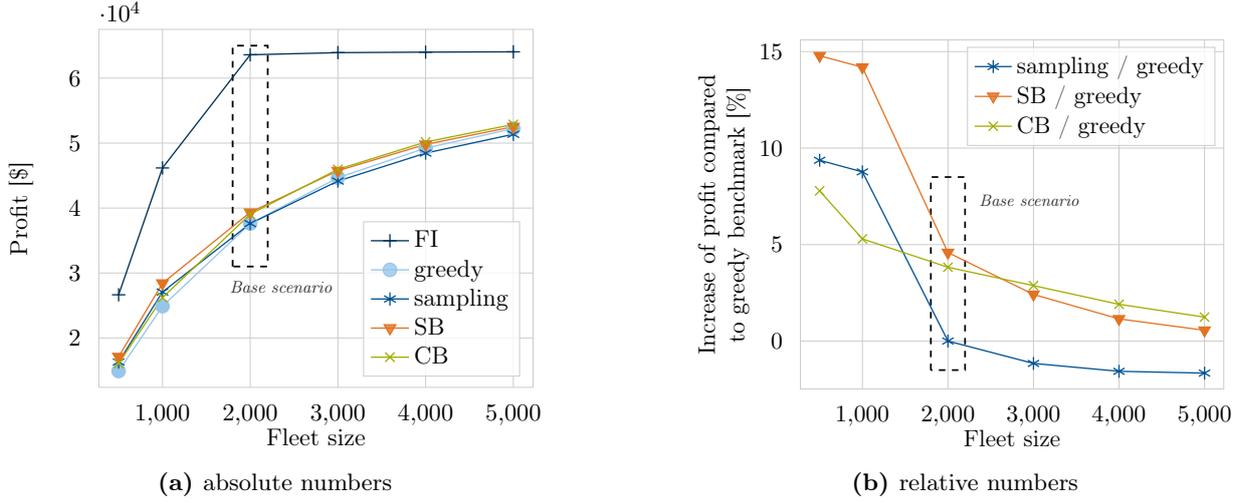
\begin{figure}[t]
\captionsetup[subfigure]{font=footnotesize}
\centering
    \begin{subfigure}[t]{0.45\textwidth}
      \centering
      \resizebox{\textwidth}{!}{ 
      \input{./figures/evaluate_vehicleTest_num_vehicles_profit_absolute.tex}}
      \caption{absolute numbers}\label{subfig:vehicleSizesProfitAbsolute}
    \end{subfigure}%
    \hfill
    \begin{subfigure}[t]{0.45\textwidth}
      \centering
      \resizebox{\textwidth}{!}{
      \input{./figures/evaluate_vehicleTest_num_vehicles_profit_ratio.tex}}
      \caption{relative numbers}\label{subfig:vehicleSizesProfitRelative}
    \end{subfigure}%
    \caption{Scenarios with different fleet sizes when optimizing the profit.}
    \label{fig:scenario:vehicleSizesProfit}
\end{figure}
\endgroup
Figure~\ref{fig:scenario:vehicleSizesProfit} shows the performance of the \gls*{abk:SB} policy and the \gls*{abk:CB} policy compared to the greedy and sampling policy, and to the full-information bound for different fleet sizes when maximizing the profit. Here, the data point with a fleet size of 2000 vehicles constitutes our base scenario. As can be seen, there is a significant gap between the full-information bound, which remains an upper bound, and all other policies, see Figure~\ref{subfig:vehicleSizesProfitAbsolute}. This gap reduces for increasing or decreasing fleet sizes as in the former case, the assignment of requests to vehicles gets less error sensitive, while in the latter case, the objective value of the full-information solution decreases overproportionally due to the reduced solution space.

To ease the comparison between all online policies, Figure~\ref{subfig:vehicleSizesProfitRelative} compares both the \gls*{abk:SB} and the \gls*{abk:CB} policy as well as the sampling policy against the greedy policy. In the base scenario, both the \gls*{abk:SB} and \gls*{abk:CB} policy outperform the greedy and the sampling policy: while the \gls*{abk:SB} policy improves upon the greedy policy by \resultVehBASESBGreedy\%, the \gls*{abk:CB} policy improves upon the greedy policy by \resultVehBASECBGreedy\%, and the sampling policy matches the performance of the greedy policy. For increasing fleet sizes, the performance difference between the \gls*{abk:SB} and \gls*{abk:CB} policies decreases. Remarkably, the sampling policy performs even worse than greedy in these cases. For decreasing fleet sizes, all policies show an increased performance gain compared to the greedy policy as the anticipation and prioritization of profitable requests is even more crucial in this setting. Here, the \gls*{abk:SB} policy still outperforms all other policies but the sampling policy shows a larger improvement compared to the \gls*{abk:CB} policy. These results indicate that the \gls*{abk:CO}-enriched \gls*{abk:ML} pipeline outperforms state-of-the-art matching and rebalancing algorithms, e.g., the sampling benchmark which is similar to the model predictive control approach of \cite{Alonso-Mora2017a}, as well as state-of-the-art \gls*{abk:DRL} algorithms which only outperform greedy policies by around 2-5\% \citep[cf.][]{Enders2022}.

\begin{result}
Both the \gls*{abk:SB} and \gls*{abk:CB} policies robustly improve over the greedy policy, while the sampling policy shows a fragile performance, often performing worse than the greedy policy. The \gls*{abk:SB} policy outperforms all other policies with improvements over the greedy policy of up to \resultVehUPTOSBGreedy\%. The \gls*{abk:CB} policy outperforms the sampling policy in most scenarios although the \gls*{abk:CB} policy does not require distributional information on future ride requests.
\end{result}

\paragraph{Varying request density:}
\begingroup
\renewcommand\baselinestretch{0.5}
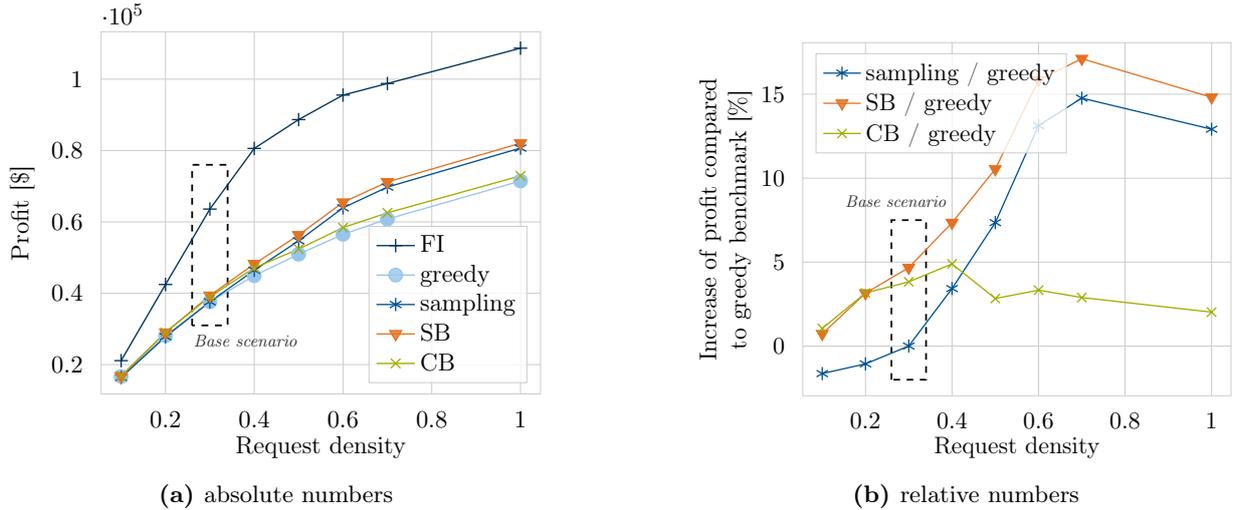
\begin{figure}[t]
\captionsetup[subfigure]{font=footnotesize}
\centering
    \begin{subfigure}[t]{0.45\textwidth}
      \centering
      \resizebox{\textwidth}{!}{
      \input{./figures/evaluate_sparsityTest_sparsity_factor_profit_absolute.tex}}
      \caption{absolute numbers}
    \end{subfigure}
    \hfill
    \begin{subfigure}[t]{0.45\textwidth}
      \centering
      \resizebox{\textwidth}{!}{
      \input{./figures/evaluate_sparsityTest_sparsity_factor_profit_ratio.tex}}
      \caption{relative numbers}
    \end{subfigure}
    \caption{Scenarios with different request densities when optimizing the profit.}
    \label{fig:scenario:requestDensitiesProfit}
\end{figure}
\endgroup
Figure~\ref{fig:scenario:requestDensitiesProfit} extends our performance analysis by showing the performance of all different policies for different request densities when fixing the number of vehicles to 2000. For scenarios with a higher request density than in our base scenario, the \gls*{abk:SB} policy outperforms all other online policies with average improvements of \resultReqHIGHERSBGreedy\% compared to a greedy policy. The sampling policy shows average improvements of \resultReqHIGHERSamplingGreedy\% for increasing request densities.
For decreasing request densities, the sampling policy performs even worse than the greedy policy.

Unsurprisingly, the sampling policy does not perform as well on low request densities as it does on high request densities. Indeed, it relies on a single sampled scenario to take decisions. When there are very few requests, a single scenario can be seen as mostly noise and does not provide accurate information on the distribution of requests.
Addressing this issue with a sampling-based approach would require sampling several scenarios, which would lead to a two-stage stochastic formulation. Such a formulation would make the policy much more computationally intensive, and difficult to apply in practice.
On the contrary, our \gls*{abk:CB} and \gls*{abk:SB} policies robustly perform better than greedy, also for low request densities. This holds even for the \gls*{abk:SB} policy although it inherits the same weakness as the sampling policy as it samples rebalancing locations.
This shows the robustness of our learning-based policies compared to the sampling policy, and further indicates their superiority in terms of solution quality.

\begin{result}
    The \gls*{abk:CB} and \gls*{abk:SB} policies yield a stable performance and improve upon greedy across all request densities. Contrarily, the sampling policy performs well on high request densities but badly on low request densities.
\end{result}

\paragraph{Maximizing number of satisfied customers:}
Figure~\ref{fig:scenario:vehicleSizesCustomers} extends our analysis to a different objective and shows the performance of all policies for a varying fleet size when maximizing the number of satisfied ride requests. As can be seen, our \gls*{abk:SB} policy outperforms all other policies, improving by up to \resultVehUPTOSBGreedyCustomerSatisfaction\% compared to a greedy policy for fleet sizes between 1000 and 5000 vehicles. 
\begingroup
\renewcommand\baselinestretch{0.5}
\begin{figure}[tbp]
\captionsetup[subfigure]{font=footnotesize}
\centering
    \begin{subfigure}[t]{0.45\textwidth}
      \centering
      \resizebox{\textwidth}{!}{
      \input{./figures/evaluate_vehicleTest_num_vehicles_amount_satisfied_customers_absolute.tex}}
      \caption{absolute numbers}
    \end{subfigure}
    \hfill
    \begin{subfigure}[t]{0.45\textwidth}
      \centering
      \resizebox{\textwidth}{!}{
      \input{./figures/evaluate_vehicleTest_num_vehicles_amount_satisfied_customers_ratio.tex}}
      \caption{relative numbers}
    \end{subfigure}
    \hfill \vfill
    \caption{Scenarios with different vehicle fleet sizes when optimizing the number of satisfied customers.}
    \label{fig:scenario:vehicleSizesCustomers}
\end{figure}
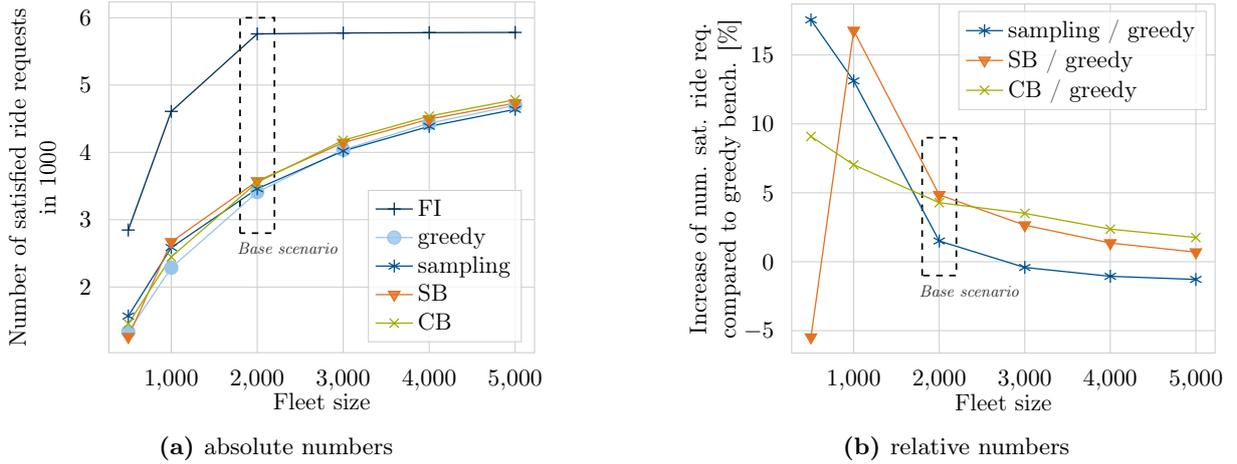
\endgroup
Again, we observe that the difference between a greedy policy and our learning-based policies decreases for increasing fleet sizes, while the sampling policy performs even worse than a greedy policy in such cases. However, we also observe one outlier: while our \gls*{abk:SB} policy performs worse than the greedy policy in the artificial case of 500 vehicles, the sampling policy outperforms the greedy policy in this case. This is due to the fact that in this case, the request-to-vehicle ratio is very large such that the pure information about possibly available requests, which we obtain by sampling, suffices to achieve a good algorithmic performance.
This is particularly possible since any request---independent of its duration or cost---adds the same contribution to the objective when solely maximizing the number of satisfied requests. Note that we observe a completely opposite behavior as soon as we change the objective towards a request-sensitive target when optimizing the total profit (cf. Figure~\ref{subfig:vehicleSizesProfitRelative}): in this case, the \gls*{abk:SB} policy yields an improvement of \resultVehFiveHundertSBGreedy\% compared to a greedy policy as the learning allows to estimate the value of potential future requests.

\begin{result}
Our \gls*{abk:SB} policy outperforms all other online policies across different objectives and various scenarios that reflect real-world applications. Our \gls*{abk:CB} policy yields lower improvements than our \gls*{abk:SB} policy, but shows a more robust performance for artificial corner cases.
\end{result}

\paragraph{Conclusion:}
Our results show that our learning-based policies outperform a greedy and a sampling policy across varying objectives and scenarios with different fleet sizes and request densities. In particular, our \gls*{abk:SB} policy yields improvements up to \resultOverallMAXSBGreedy\% compared to a greedy policy for scenarios that reflect real-world applications. In these cases, our \gls*{abk:CB} policy improves upon a greedy policy up to \resultOverallMAXCBGreedy\%, while a sampling policy shows an unstable performance that sometimes yields improvements, but oftentimes performs even worse than a greedy policy. In artificial corner cases, our \gls*{abk:SB} policy may also perform worse than a greedy policy as it by design inherits some of the structural difficulties of any sampling-based approach. Here, our \gls*{abk:CB} policy still shows a robust performance that improves upon a greedy policy.

\subsection{Extended analysis}
\label{sec:AppliedDiscussion}

In the following, we deepen our analyses of the differences between the proposed benchmark policies. To this end, we analyze structural differences between the different policies and discuss their impact on external effects, as well as their efficiency with respect to computational times.

\begin{figure}[tbp]
    \centering
    \begin{minipage}[c]{.453\linewidth}
      \includegraphics[width=\linewidth]{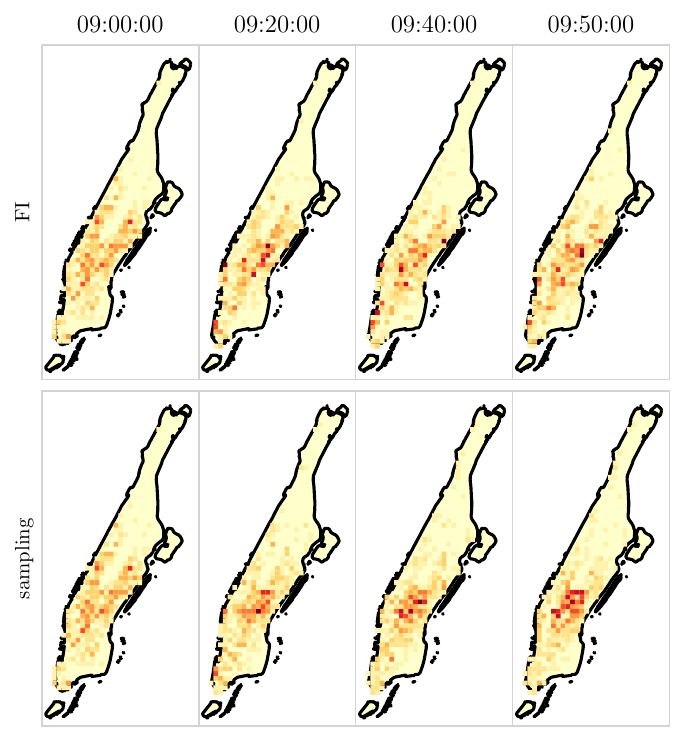}
    \end{minipage}
    \begin{minipage}[c]{.537\linewidth}
      \includegraphics[width=\linewidth]{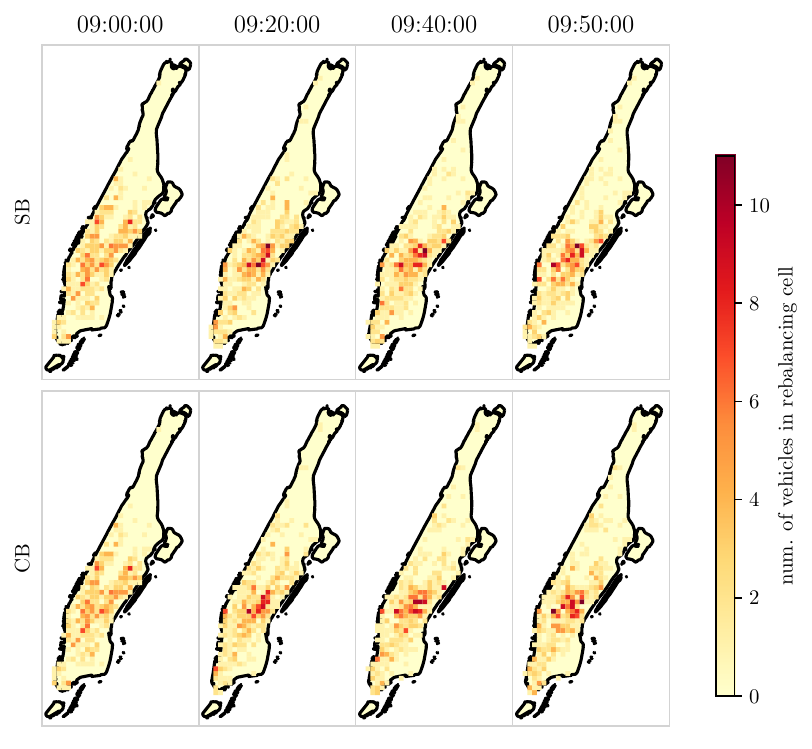}
    \end{minipage}
    \caption{Vehicle distribution for different benchmark policies over different points in time for a representative instance.
    }
    \label{fig:vehicle_distribution}
\end{figure}

\paragraph{Structural Policy Properties:}
Figure~\ref{fig:vehicle_distribution} shows the evolution of the vehicle distribution for a representative instance over different points in time for the different benchmark policies. 
We interpret this evolution of vehicle distributions as an indicator of how the respective benchmark policies differ in their dispatching and rebalancing behavior, and refer to our git repository (\url{https://github.com/tumBAIS/ML-CO-pipeline-AMoD-control}) for an interactive visualization of idling, dispatching, and rebalancing characteristics of all policies. Compared to the full-information bound, the sampling policy leads to a compressed vehicle distribution in the city center, which occurs due to a biased rebalancing behavior that overvalues high-request-density areas over low-request-density areas. This leads to a significant amount of vehicles idling in high-request-density areas, which leads to missed requests in mid to low-density areas and worsens the policy's performance. The \gls*{abk:SB} policy also shows a vehicle distribution that revolves around the city center and its high-request-density areas, which results from the fact that both policies leverage the identical historical request distribution for sampling future ride requests. However, the \gls*{abk:SL} element allows anticipating better the value of relocations, which leads to vehicle distributions that still revolve around high-request-density areas but distribute a significantly larger share of vehicles to mid and low-density areas, which yields the performance improvement of our \gls*{abk:SB} policy.
The \gls*{abk:CB} policy leads to a slightly more uniform vehicle distribution over the operating area as, by design, it restricts the maximum number of vehicles rebalancing to a rebalancing cell via the number of capacity vertices in each rebalancing cell. With these rebalancing actions, the \gls*{abk:CB} policy yields a performance that is always better than greedy but remains below the performance of the \gls*{abk:SB} policy.
\begin{result}
    While sample-based policies tend to rebalance too many vehicles to high-request-density areas, \gls*{abk:CB} policies lead to a slightly more distributed rebalancing.
\end{result}

\paragraph{Externalities:}
Figure~\ref{fig:results_metric} shows the externalities of our online policies compared to a greedy policy for different fleet sizes when optimizing the total profit. Here, we particularly focus on the number of served requests, the average number of kilometers driven per request, and the average number of kilometers driven per vehicle. 
\begingroup
\renewcommand\baselinestretch{0.5}
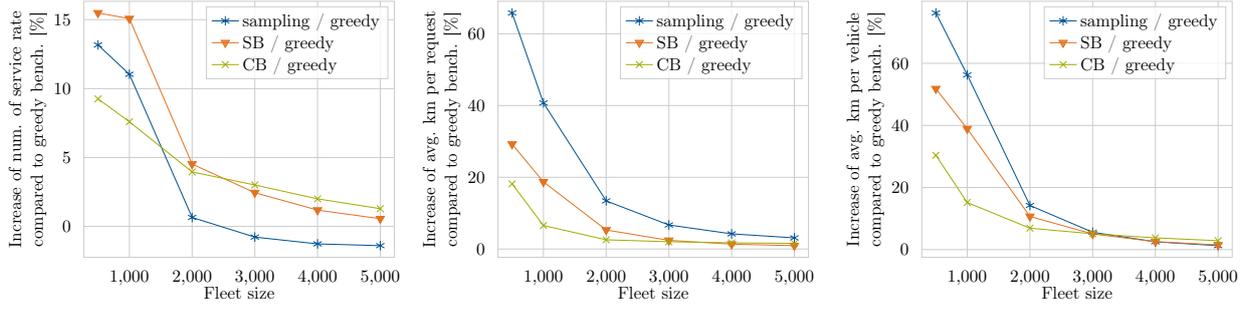
\begin{figure}[t]
    \begin{minipage}[c]{0.33\textwidth}
      \centering
      \resizebox{\textwidth}{!}{ 
      \input{./figures/evaluate_vehicleTest_num_vehicles_service_ratio_ratio.tex}}
    \end{minipage}%
    \begin{minipage}[c]{.33\textwidth}
      \centering
      \resizebox{\textwidth}{!}{
      \input{./figures/evaluate_vehicleTest_num_vehicles_km_per_request_ratio.tex}}
    \end{minipage}
    \begin{minipage}[c]{.33\textwidth}
      \centering
      \resizebox{\textwidth}{!}{
      \input{./figures/evaluate_vehicleTest_num_vehicles_avg_distance_trav_ratio.tex}}
    \end{minipage}
    \caption{Externalities for all online policies compared to a greedy policy over various fleet sizes when optimizing the total profit.}
    \label{fig:results_metric}
\end{figure}
\endgroup 
As can be seen, the increase in the number of served requests is proportional to the profit increase for each policy. Accordingly, the \gls*{abk:SB} policy yields the highest improvement for the number of served customers, and thus contributes to increase customer satisfaction. Remarkably, the \gls*{abk:SB} policy yields a lower increase of average kilometers per request and per vehicle than the sampling policy, while at the same time yielding a higher number of passengers served and a higher profit. 
This highlights the efficacy of the \gls*{abk:SB} policy, even for key performance indicators that are anticorrelated. We note that the \gls*{abk:SB} policy still yields higher average kilometers driven per request and vehicle than the \gls*{abk:CB} and greedy policies for small fleet sizes. This is due to the fact that these policies show a significantly lower performance on the anticorrelated key performance indicators, such as the number of satisfied requests and the total profit. For scenarios with large fleet sizes, when the \gls*{abk:SB} and \gls*{abk:CB} policies lead to a similar performance with regard to total profit, also the KPIs for average kilometers driven per request and vehicle are similar.

\begin{result}
    The \gls*{abk:SB} policy yields a good trade-off between anticorrelated key performance indicators and outperforms the sampling policy on both the number of satisfied customers and the total profit as well as on the distance driven per request and vehicle.
\end{result}

\paragraph{Computational time:}

\begin{table}[t]
    \centering
    \def\arraystretch{0.5}
    \small
    \begin{tabular}{c c c c c c c}
        \toprule
        & \multicolumn{6}{c}{Fleet size} \\
        & 500& 1000& 2000& 3000& 4000& 5000\\
        \midrule
        FI & 2.314 & 4.396 & 10.155 & 14.949 & 27.59 & 28.792\\
        greedy & 0.012 & 0.05 & 0.208 & 0.485 & 0.863 & 1.384\\
        sampling & 0.041 & 0.169 & 0.588 & 1.201 & 2.003 & 3.01\\
        SB & 0.069 & 0.179 & 0.576 & 1.186 & 2.025 & 3.038\\
        CB & 0.093 & 0.312 & 0.905 & 1.798 & 2.966 & 4.167\\
        \bottomrule
    \end{tabular}
    \caption{Average computational time to calculate dispatching and rebalancing actions [sec].}
    \label{tab:my_label}
\end{table}

Table~\ref{tab:my_label} presents the average computational time to compute dispatching and rebalancing actions for all benchmark policies.
We consider a system time period of~\systemTimePeriod~minute for all online policies and do not incorporate the time to construct the digraph.
As can be seen, all policies benefit from the polynomial algorithm introduced in Section~\ref{sec:kdspp}. While the greedy policy shows lower computational times compared to the sampling and learning-based policies, all online policies show computational times below a few seconds, which allow for application in practice when operating an \gls*{abk:amod} system with a one-minute system time period.

\paragraph{Deep learning ML-layer:}

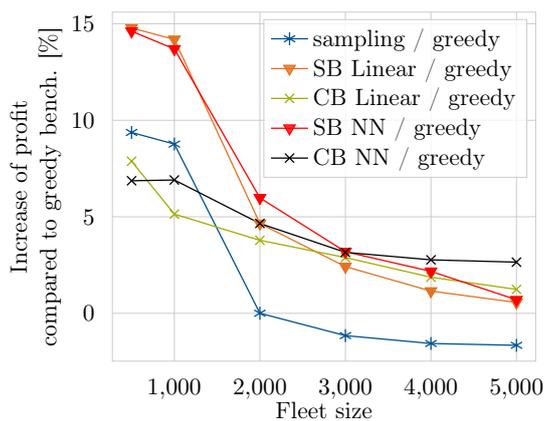
\begin{figure}[b]
      \centering
      \resizebox{0.45\textwidth}{!}{
      \input{./figures/evaluateNN_vehicleTestHyperparameter_num_vehicles_profit_ratio.tex}}
    \caption{Comparison of performance with linear predictor and \gls{abk:FNN} predictor.}
    \label{fig:comparisonNN}
\end{figure}

In the previous results, we used a linear predictor in the \gls*{abk:ML}-layer to predict the digraph weights $\theta$, according to equation~\eqref{eq:MLpredictor}. Instead of using a linear predictor, we can also use a differentiable deep learning predictor of the general form $\varphi_w: \phi(a, x_i) \rightarrow \theta_a \quad \forall a \in A$. For learning a deep learning predictor, some adaptions are in order in comparison to learning a linear predictor. First, instead of using the \gls*{abk:BFGS} method, we use a \gls*{abk:SGD} method to optimize the learning problem~\eqref{eq:learning_problem}. Second, the \gls*{abk:SGD} method updates the weights~$w$ for each training instance instead of updating the weights with a mean gradient over all training instances. Third, for deriving a subgradient~\eqref{eq:gradient} we add a perturbation~$Z$ to~$\theta$ instead of perturbing the weights~$w$.

Figure \ref{fig:comparisonNN} shows the performances of the \gls*{abk:CO}-enriched \gls*{abk:ML} pipeline for different vehicle fleet sizes when using an \acrfull{abk:FNN} predictor, and a linear predictor. We tuned the hyperparameters of the \gls{abk:FNN}, e.g., the learning rate, and the network architecture individually in each scenario. Surprisingly, the \gls*{abk:CO}-enriched \gls*{abk:ML} pipeline including the \gls{abk:FNN} does not outperform the pipeline including the linear predictor in all scenarios: In scenarios with a small vehicle fleet size the linear predictor leads to better performance. This is due to an overfitting issue of the \gls{abk:FNN} on the training data, which we detail in Appendix~\ref{appendix:LinearNNOverfitting}. In all other scenarios, the \gls{abk:FNN} leads to a slightly better performance.

\begin{result}
    The \gls{abk:FNN} predictor in the \gls*{abk:ML}-layer improves the performance of the \gls*{abk:CO}-enriched \gls*{abk:ML} pipeline in most scenarios. However, the \gls*{abk:CO}-enriched \gls*{abk:ML} pipeline including the \gls{abk:FNN} does not generalize well for all scenarios.
\end{result}

%% file: figures/evaluate_vehicleTest_num_vehicles_profit_absolute.tex
\begin{tikzpicture}

\definecolor{chocolate22711434}{RGB}{227,114,34}
\definecolor{darkgoldenrod1621730}{RGB}{162,173,0}
\definecolor{darkslategray38}{RGB}{38,38,38}
\definecolor{lightgray204}{RGB}{204,204,204}
\definecolor{midnightblue05189}{RGB}{0,51,89}
\definecolor{skyblue152198234}{RGB}{152,198,234}
\definecolor{teal082147}{RGB}{0,82,147}

\begin{axis}[
axis line style={lightgray204},
legend cell align={left},
legend style={
  fill opacity=0.8,
  draw opacity=1,
  text opacity=1,
  at={(0.97,0.03)},
  anchor=south east,
  draw=lightgray204
},
scaled y ticks=true,
tick align=outside,
x grid style={lightgray204},
xlabel={Fleet size},
xmajorgrids,
xmajorticks=false,
xmajorticks=true,
xmin=275, xmax=5225,
xtick style={color=darkslategray38},
xtick style={draw=none},
y grid style={lightgray204},
ylabel style={align=center},
ylabel={Profit [\$]},
ymajorgrids,
ymajorticks=false,
ymajorticks=true,
ymin=12401.4857744726, ymax=67504.6911535965,
ytick style={color=darkslategray38},
ytick style={draw=none}
]
\draw[draw=black,thick,dashed] (axis cs:1800,31000) rectangle (axis cs:2200,65000);
\addplot [semithick, midnightblue05189, mark=+, mark size=3, mark options={solid}]
table {%
500 26647.331518853
1000 46175.954212998
2000 63598.1727681903
3000 63922.9981923327
4000 64005.3900575796
5000 64050.4506844398
};
\addlegendentry{FI}
\addplot [semithick, skyblue152198234, mark=*, mark size=3, mark options={solid}]
table {%
500 14906.1769280691
1000 24895.0171666284
2000 37614.1097186284
3000 44671.8555231926
4000 49228.4353871138
5000 52225.5574648163
};
\addlegendentry{greedy}
\addplot [semithick, teal082147, mark=asterisk, mark size=3, mark options={solid}]
table {%
500 16302.8120806746
1000 27077.7627138474
2000 37615.5946559527
3000 44156.0412654607
4000 48458.9997015816
5000 51360.0176189518
};
\addlegendentry{sampling}
\addplot [semithick, chocolate22711434, mark=triangle*, mark size=3, mark options={solid,rotate=180}]
table {%
500 17109.534687486
1000 28427.728899312
2000 39340.3256848816
3000 45752.8260296339
4000 49795.5420191313
5000 52519.3012946185
};
\addlegendentry{SB}
\addplot [semithick, darkgoldenrod1621730, mark=x, mark size=3, mark options={solid}]
table {%
500 16066.8182582072
1000 26213.3960456125
2000 39053.9162914522
3000 45956.7302600405
4000 50167.5742994165
5000 52875.4953677586
};
\addlegendentry{CB}
\draw (axis cs:1700,27000) node[
  scale=0.7,
  anchor=base west,
  text=darkslategray38,
  rotate=0.0
]{\itshape Base scenario};
\end{axis}

\end{tikzpicture}

%% file: figures/evaluate_vehicleTest_num_vehicles_profit_ratio.tex
\begin{tikzpicture}

\definecolor{chocolate22711434}{RGB}{227,114,34}
\definecolor{darkgoldenrod1621730}{RGB}{162,173,0}
\definecolor{darkslategray38}{RGB}{38,38,38}
\definecolor{lightgray204}{RGB}{204,204,204}
\definecolor{teal082147}{RGB}{0,82,147}

\begin{axis}[
axis line style={lightgray204},
legend cell align={left},
legend style={fill opacity=0.8, draw opacity=1, text opacity=1, draw=lightgray204},
scaled y ticks=false,
tick align=outside,
x grid style={lightgray204},
xlabel={Fleet size},
xmajorgrids,
xmajorticks=false,
xmajorticks=true,
xmin=275, xmax=5225,
xtick style={color=darkslategray38},
xtick style={draw=none},
y grid style={lightgray204},
y tick label style={/pgf/number format/fixed},
ylabel style={align=center},
ylabel={Increase of profit compared \\ to greedy benchmark [\%]},
ymajorgrids,
ymajorticks=false,
ymajorticks=true,
ymin=-2.47925183186852, ymax=15.6034491151588,
ytick style={color=darkslategray38},
ytick style={draw=none}
]
\draw[draw=black,thick,dashed] (axis cs:1800,-1.5) rectangle (axis cs:2200,8.5);
\addplot [semithick, teal082147, mark=asterisk, mark size=3, mark options={solid}]
table {%
500 9.36950607352279
1000 8.76780093224818
2000 0.00394781994159887
3000 -1.15467390304398
4000 -1.56299033166839
5000 -1.65731087973091
};
\addlegendentry{sampling / greedy}
\addplot [semithick, chocolate22711434, mark=triangle*, mark size=3, mark options={solid,rotate=180}]
table {%
500 14.7815081630212
1000 14.190437022149
2000 4.58927774488385
3000 2.41980211876394
4000 1.15198995775104
5000 0.562452262955945
};
\addlegendentry{SB / greedy}
\addplot [semithick, darkgoldenrod1621730, mark=x, mark size=3, mark options={solid}]
table {%
500 7.78631124357956
1000 5.29575404652186
2000 3.82783637202696
3000 2.87625110217521
4000 1.90771635319642
5000 1.24448246125493
};
\addlegendentry{CB / greedy}
\draw (axis cs:2300,7) node[
  scale=0.7,
  anchor=base west,
  text=darkslategray38,
  rotate=0.0
]{\itshape Base scenario};
\end{axis}

\end{tikzpicture}

%% file: figures/evaluate_sparsityTest_sparsity_factor_profit_absolute.tex
\begin{tikzpicture}

\definecolor{chocolate22711434}{RGB}{227,114,34}
\definecolor{darkgoldenrod1621730}{RGB}{162,173,0}
\definecolor{darkslategray38}{RGB}{38,38,38}
\definecolor{lightgray204}{RGB}{204,204,204}
\definecolor{midnightblue05189}{RGB}{0,51,89}
\definecolor{skyblue152198234}{RGB}{152,198,234}
\definecolor{teal082147}{RGB}{0,82,147}

\begin{axis}[
axis line style={lightgray204},
legend cell align={left},
legend style={
  fill opacity=0.8,
  draw opacity=1,
  text opacity=1,
  at={(0.97,0.03)},
  anchor=south east,
  draw=lightgray204
},
scaled y ticks=true,
tick align=outside,
x grid style={lightgray204},
xlabel={Request density},
xmajorgrids,
xmajorticks=false,
xmajorticks=true,
xmin=0.055, xmax=1.045,
xtick style={color=darkslategray38},
xtick style={draw=none},
y grid style={lightgray204},
ylabel style={align=center},
ylabel={Profit [\$]},
ymajorgrids,
ymajorticks=false,
ymajorticks=true,
ymin=11796.3604217276, ymax=113290.490926078,
ytick style={color=darkslategray38},
ytick style={draw=none}
]
\draw[draw=black,thick,dashed] (axis cs:0.26,31000) rectangle (axis cs:0.34,76000);
\addplot [semithick, midnightblue05189, mark=+, mark size=3, mark options={solid}]
table {%
0.1 21126.1833904053
0.2 42471.050024362
0.3 63598.1727681903
0.4 80585.4946881512
0.5 88667.7758126678
0.6 95539.9150147555
0.7 98778.1372454884
1 108677.121357698
};
\addlegendentry{FI}
\addplot [semithick, skyblue152198234, mark=*, mark size=3, mark options={solid}]
table {%
0.1 16682.6335953898
0.2 28123.0287107224
0.3 37614.1097186284
0.4 44919.5958739198
0.5 50949.4974924231
0.6 56523.1726972258
0.7 60780.7679185404
1 71459.2582738762
};
\addlegendentry{greedy}
\addplot [semithick, teal082147, mark=asterisk, mark size=3, mark options={solid}]
table {%
0.1 16409.7299901071
0.2 27825.4176102944
0.3 37615.5946559527
0.4 46457.9078956635
0.5 54695.2748532166
0.6 63943.450121774
0.7 69754.4383461335
1 80687.5664236275
};
\addlegendentry{sampling}
\addplot [semithick, chocolate22711434, mark=triangle*, mark size=3, mark options={solid,rotate=180}]
table {%
0.1 16804.256874062
0.2 29009.3665032908
0.3 39371.6797825297
0.4 48219.2758051349
0.5 56317.7566338739
0.6 65473.9486890193
0.7 71181.327185108
1 82042.0088094851
};
\addlegendentry{SB}
\addplot [semithick, darkgoldenrod1621730, mark=x, mark size=3, mark options={solid}]
table {%
0.1 16858.031973267
0.2 29012.2727440785
0.3 39053.9162914522
0.4 47116.7865475167
0.5 52392.8581627927
0.6 58405.6103392521
0.7 62538.8301658452
1 72904.701438667
};
\addlegendentry{CB}
\draw (axis cs:0.25,25000) node[
  scale=0.7,
  anchor=base west,
  text=darkslategray38,
  rotate=0.0
]{\itshape Base scenario};
\end{axis}

\end{tikzpicture}

%% file: figures/evaluate_sparsityTest_sparsity_factor_profit_ratio.tex
\begin{tikzpicture}

\definecolor{chocolate22711434}{RGB}{227,114,34}
\definecolor{darkgoldenrod1621730}{RGB}{162,173,0}
\definecolor{darkslategray38}{RGB}{38,38,38}
\definecolor{lightgray204}{RGB}{204,204,204}
\definecolor{teal082147}{RGB}{0,82,147}

\begin{axis}[
axis line style={lightgray204},
legend cell align={left},
legend style={
  fill opacity=0.8,
  draw opacity=1,
  text opacity=1,
  at={(0.03,0.97)},
  anchor=north west,
  draw=lightgray204
},
scaled y ticks=false,
tick align=outside,
x grid style={lightgray204},
xlabel={Request density},
xmajorgrids,
xmajorticks=false,
xmajorticks=true,
xmin=0.055, xmax=1.045,
xtick style={color=darkslategray38},
xtick style={draw=none},
y grid style={lightgray204},
y tick label style={/pgf/number format/fixed},
ylabel style={align=center},
ylabel={Increase of profit compared \\ to greedy benchmark [\%]},
ymajorgrids,
ymajorticks=false,
ymajorticks=true,
ymin=-2.95557978475253, ymax=18.0671754798031,
ytick style={color=darkslategray38},
ytick style={draw=none}
]
\draw[draw=black,thick,dashed] (axis cs:0.26,-2) rectangle (axis cs:0.34,7.5);
\addplot [semithick, teal082147, mark=asterisk, mark size=3, mark options={solid}]
table {%
0.1 -1.63585445740502
0.2 -1.05824697435419
0.3 0.00394781994159887
0.4 3.42459007436618
0.5 7.35194171709066
0.6 13.1278501726998
0.7 14.7639964661516
1 12.9140833149747
};
\addlegendentry{sampling / greedy}
\addplot [semithick, chocolate22711434, mark=triangle*, mark size=3, mark options={solid,rotate=180}]
table {%
0.1 0.729041239063065
0.2 3.15164416210438
0.3 4.67263502193401
0.4 7.34574714446808
0.5 10.536431968243
0.6 15.8355866535298
0.7 17.1115956950505
1 14.809488359156
};
\addlegendentry{SB / greedy}
\addplot [semithick, darkgoldenrod1621730, mark=x, mark size=3, mark options={solid}]
table {%
0.1 1.05138302579295
0.2 3.16197818699745
0.3 3.82783637202696
0.4 4.89138566554317
0.5 2.83292425128279
0.6 3.33038212860021
0.7 2.89246468498235
1 2.02275142466631
};
\addlegendentry{CB / greedy}
\draw (axis cs:0.14,8.3) node[
  scale=0.7,
  anchor=base west,
  text=darkslategray38,
  rotate=0.0
]{\itshape Base scenario};
\end{axis}

\end{tikzpicture}

%% file: figures/evaluate_vehicleTest_num_vehicles_amount_satisfied_customers_absolute.tex
\begin{tikzpicture}

\definecolor{chocolate22711434}{RGB}{227,114,34}
\definecolor{darkgoldenrod1621730}{RGB}{162,173,0}
\definecolor{darkslategray38}{RGB}{38,38,38}
\definecolor{lightgray204}{RGB}{204,204,204}
\definecolor{midnightblue05189}{RGB}{0,51,89}
\definecolor{skyblue152198234}{RGB}{152,198,234}
\definecolor{teal082147}{RGB}{0,82,147}

\begin{axis}[
axis line style={lightgray204},
legend cell align={left},
legend style={
  fill opacity=0.8,
  draw opacity=1,
  text opacity=1,
  at={(0.97,0.03)},
  anchor=south east,
  draw=lightgray204
},
scaled y ticks=true,
tick align=outside,
x grid style={lightgray204},
xlabel={Fleet size},
xmajorgrids,
xmajorticks=false,
xmajorticks=true,
xmin=275, xmax=5225,
xtick style={color=darkslategray38},
xtick style={draw=none},
y grid style={lightgray204},
ylabel style={align=center},
ylabel={Number of satisfied ride requests \\ in 1000},
ymajorgrids,
ymajorticks=false,
ymajorticks=true,
ymin=1.023735, ymax=6.236965,
ytick style={color=darkslategray38},
ytick style={draw=none}
]
\draw[draw=black,thick,dashed] (axis cs:1800,2.8) rectangle (axis cs:2200,6);
\addplot [semithick, midnightblue05189, mark=+, mark size=3, mark options={solid}]
table {%
500 2.8465
1000 4.609
2000 5.7603
3000 5.77345
4000 5.7788
5000 5.78205
};
\addlegendentry{FI}
\addplot [semithick, skyblue152198234, mark=*, mark size=3, mark options={solid}]
table {%
500 1.3337
1000 2.28615
2000 3.4085
3000 4.03775
4000 4.4335
5000 4.6997
};
\addlegendentry{greedy}
\addplot [semithick, teal082147, mark=asterisk, mark size=3, mark options={solid}]
table {%
500 1.56775
1000 2.58615
2000 3.4598
3000 4.021
4000 4.38665
5000 4.6393
};
\addlegendentry{sampling}
\addplot [semithick, chocolate22711434, mark=triangle*, mark size=3, mark options={solid,rotate=180}]
table {%
500 1.2607
1000 2.66965
2000 3.5738
3000 4.1449
4000 4.49375
5000 4.73245
};
\addlegendentry{SB}
\addplot [semithick, darkgoldenrod1621730, mark=x, mark size=3, mark options={solid}]
table {%
500 1.45495
1000 2.44695
2000 3.5545
3000 4.17945
4000 4.53835
5000 4.78225
};
\addlegendentry{CB}
\draw (axis cs:1700,2.5) node[
  scale=0.7,
  anchor=base west,
  text=darkslategray38,
  rotate=0.0
]{\itshape Base scenario};
\end{axis}

\end{tikzpicture}

%% file: figures/evaluate_vehicleTest_num_vehicles_amount_satisfied_customers_ratio.tex
\begin{tikzpicture}

\definecolor{chocolate22711434}{RGB}{227,114,34}
\definecolor{darkgoldenrod1621730}{RGB}{162,173,0}
\definecolor{darkslategray38}{RGB}{38,38,38}
\definecolor{lightgray204}{RGB}{204,204,204}
\definecolor{teal082147}{RGB}{0,82,147}

\begin{axis}[
axis line style={lightgray204},
legend cell align={left},
legend style={fill opacity=0.8, draw opacity=1, text opacity=1, draw=lightgray204},
scaled y ticks=false,
tick align=outside,
x grid style={lightgray204},
xlabel={Fleet size},
xmajorgrids,
xmajorticks=false,
xmajorticks=true,
xmin=275, xmax=5225,
xtick style={color=darkslategray38},
xtick style={draw=none},
y grid style={lightgray204},
y tick label style={/pgf/number format/fixed},
ylabel style={align=center},
ylabel={Increase of num. sat. ride req. \\ compared to greedy bench. [\%]},
ymajorgrids,
ymajorticks=false,
ymajorticks=true,
ymin=-6.62461573067406, ymax=18.7000449876284,
ytick style={color=darkslategray38},
ytick style={draw=none}
]
\draw[draw=black,thick,dashed] (axis cs:1800,-1) rectangle (axis cs:2200,9);
\addplot [semithick, teal082147, mark=asterisk, mark size=3, mark options={solid}]
table {%
500 17.5489240458874
1000 13.1224985237189
2000 1.50506087721874
3000 -0.414834994737169
4000 -1.05672719070712
5000 -1.2851884162819
};
\addlegendentry{sampling / greedy}
\addplot [semithick, chocolate22711434, mark=triangle*, mark size=3, mark options={solid,rotate=180}]
table {%
500 -5.47349478893304
1000 16.7749272794874
2000 4.84964060437145
3000 2.65370565290075
4000 1.35897146723807
5000 0.696852990616386
};
\addlegendentry{SB / greedy}
\addplot [semithick, darkgoldenrod1621730, mark=x, mark size=3, mark options={solid}]
table {%
500 9.0912499062758
1000 7.03365920871334
2000 4.28340912424823
3000 3.50938022413473
4000 2.3649486861396
5000 1.75649509543161
};
\addlegendentry{CB / greedy}
\draw (axis cs:1700,-2.5) node[
  scale=0.7,
  anchor=base west,
  text=darkslategray38,
  rotate=0.0
]{\itshape Base scenario};
\end{axis}

\end{tikzpicture}

%% file: figures/evaluate_vehicleTest_num_vehicles_service_ratio_ratio.tex
\begin{tikzpicture}

\definecolor{chocolate22711434}{RGB}{227,114,34}
\definecolor{darkgoldenrod1621730}{RGB}{162,173,0}
\definecolor{darkslategray38}{RGB}{38,38,38}
\definecolor{lightgray204}{RGB}{204,204,204}
\definecolor{teal082147}{RGB}{0,82,147}

\begin{axis}[
axis line style={lightgray204},
legend cell align={left},
legend style={fill opacity=0.8, draw opacity=1, text opacity=1, draw=lightgray204},
scaled y ticks=false,
tick align=outside,
x grid style={lightgray204},
xlabel={Fleet size},
xmajorgrids,
xmajorticks=false,
xmajorticks=true,
xmin=275, xmax=5225,
xtick style={color=darkslategray38},
xtick style={draw=none},
y grid style={lightgray204},
y tick label style={/pgf/number format/fixed},
ylabel style={align=center},
ylabel={Increase of num. of service rate \\ compared to greedy bench. [\%]},
ymajorgrids,
ymajorticks=false,
ymajorticks=true,
ymin=-2.23927392741386, ymax=16.3328781039414,
ytick style={color=darkslategray38},
ytick style={draw=none}
]
\addplot [semithick, teal082147, mark=asterisk, mark size=3, mark options={solid}]
table {%
500 13.1679874317008
1000 11.0394990390833
2000 0.650585662446833
3000 -0.774426324890088
4000 -1.26714039633075
5000 -1.39508519871589
};
\addlegendentry{sampling / greedy}
\addplot [semithick, chocolate22711434, mark=triangle*, mark size=3, mark options={solid,rotate=180}]
table {%
500 15.4886893752434
1000 15.0752957229666
2000 4.5224125421661
3000 2.44868545998468
4000 1.18473105811832
5000 0.560584317513196
};
\addlegendentry{SB / greedy}
\addplot [semithick, darkgoldenrod1621730, mark=x, mark size=3, mark options={solid}]
table {%
500 9.25434064530888
1000 7.59599160712904
2000 3.95213912417295
3000 3.00704264839701
4000 2.00047380242642
5000 1.28734358818275
};
\addlegendentry{CB / greedy}
\end{axis}

\end{tikzpicture}

%% file: figures/evaluate_vehicleTest_num_vehicles_km_per_request_ratio.tex
\begin{tikzpicture}

\definecolor{chocolate22711434}{RGB}{227,114,34}
\definecolor{darkgoldenrod1621730}{RGB}{162,173,0}
\definecolor{darkslategray38}{RGB}{38,38,38}
\definecolor{lightgray204}{RGB}{204,204,204}
\definecolor{teal082147}{RGB}{0,82,147}

\begin{axis}[
axis line style={lightgray204},
legend cell align={left},
legend style={fill opacity=0.8, draw opacity=1, text opacity=1, draw=lightgray204},
scaled y ticks=false,
tick align=outside,
x grid style={lightgray204},
xlabel={Fleet size},
xmajorgrids,
xmajorticks=false,
xmajorticks=true,
xmin=275, xmax=5225,
xtick style={color=darkslategray38},
xtick style={draw=none},
y grid style={lightgray204},
y tick label style={/pgf/number format/fixed},
ylabel style={align=center},
ylabel={Increase of avg. km per request \\ compared to greedy bench. [\%]},
ymajorgrids,
ymajorticks=false,
ymajorticks=true,
ymin=-2.29878764449313, ymax=69.0471099208354,
ytick style={color=darkslategray38},
ytick style={draw=none}
]
\addplot [semithick, teal082147, mark=asterisk, mark size=3, mark options={solid}]
table {%
500 65.8041145769568
1000 40.7262275878067
2000 13.4209325534839
3000 6.70790642329918
4000 4.23253592376575
5000 3.11618268070979
};
\addlegendentry{sampling / greedy}
\addplot [semithick, chocolate22711434, mark=triangle*, mark size=3, mark options={solid,rotate=180}]
table {%
500 29.233176610495
1000 18.7245864513213
2000 5.31203632650365
3000 2.4203775524522
4000 1.36162805992988
5000 0.94420769938544
};
\addlegendentry{SB / greedy}
\addplot [semithick, darkgoldenrod1621730, mark=x, mark size=3, mark options={solid}]
table {%
500 18.1615237964699
1000 6.55432067109168
2000 2.57311691891626
3000 2.0385499801291
4000 1.72852574991458
5000 1.57139784134728
};
\addlegendentry{CB / greedy}
\end{axis}

\end{tikzpicture}

%% file: figures/evaluate_vehicleTest_num_vehicles_avg_distance_trav_ratio.tex
\begin{tikzpicture}

\definecolor{chocolate22711434}{RGB}{227,114,34}
\definecolor{darkgoldenrod1621730}{RGB}{162,173,0}
\definecolor{darkslategray38}{RGB}{38,38,38}
\definecolor{lightgray204}{RGB}{204,204,204}
\definecolor{teal082147}{RGB}{0,82,147}

\begin{axis}[
axis line style={lightgray204},
legend cell align={left},
legend style={fill opacity=0.8, draw opacity=1, text opacity=1, draw=lightgray204},
scaled y ticks=false,
tick align=outside,
x grid style={lightgray204},
xlabel={Fleet size},
xmajorgrids,
xmajorticks=false,
xmajorticks=true,
xmin=275, xmax=5225,
xtick style={color=darkslategray38},
xtick style={draw=none},
y grid style={lightgray204},
y tick label style={/pgf/number format/fixed},
ylabel style={align=center},
ylabel={Increase of avg. km per vehicle \\ compared to greedy bench. [\%]},
ymajorgrids,
ymajorticks=false,
ymajorticks=true,
ymin=-2.47889659680291, ymax=79.9685863575651,
ytick style={color=darkslategray38},
ytick style={draw=none}
]
\addplot [semithick, teal082147, mark=asterisk, mark size=3, mark options={solid}]
table {%
500 76.2209734960029
1000 56.251414357088
2000 14.21500950262
3000 5.58678259388668
4000 2.50802905903017
5000 1.26871626475927
};
\addlegendentry{sampling / greedy}
\addplot [semithick, chocolate22711434, mark=triangle*, mark size=3, mark options={solid,rotate=180}]
table {%
500 51.8212210506707
1000 38.9137653184172
2000 10.6782854484474
3000 5.03596208184207
4000 2.55303613295333
5000 1.50195676008165
};
\addlegendentry{SB / greedy}
\addplot [semithick, darkgoldenrod1621730, mark=x, mark size=3, mark options={solid}]
table {%
500 30.3726571231504
1000 15.16445700981
2000 6.89127494885663
3000 5.0966323067447
4000 3.77342201690317
5000 2.86532547646171
};
\addlegendentry{CB / greedy}
\end{axis}

\end{tikzpicture}

%% file: figures/evaluateNN_vehicleTestHyperparameter_num_vehicles_profit_ratio.tex
\begin{tikzpicture}

\definecolor{chocolate22711434}{RGB}{227,114,34}
\definecolor{darkgoldenrod1621730}{RGB}{162,173,0}
\definecolor{darkslategray38}{RGB}{38,38,38}
\definecolor{lightgray204}{RGB}{204,204,204}
\definecolor{teal082147}{RGB}{0,82,147}

\begin{axis}[
axis line style={lightgray204},
legend cell align={left},
legend style={fill opacity=0.8, draw opacity=1, text opacity=1, draw=lightgray204},
scaled y ticks=false,
tick align=outside,
x grid style={lightgray204},
xlabel={Fleet size},
xmajorgrids,
xmajorticks=false,
xmajorticks=true,
xmin=275, xmax=5225,
xtick style={color=darkslategray38},
xtick style={draw=none},
y grid style={lightgray204},
y tick label style={/pgf/number format/fixed},
ylabel style={align=center},
ylabel={Increase of profit \\ compared to greedy bench. [\%]},
ymajorgrids,
ymajorticks=false,
ymajorticks=true,
ymin=-2.47925183186852, ymax=15.6034491151587,
ytick style={color=darkslategray38},
ytick style={draw=none}
]
\addplot [semithick, teal082147, mark=asterisk, mark size=3, mark options={solid}]
table {%
500 9.36950607352279
1000 8.76780093224818
2000 0.00394781994159887
3000 -1.15467390304398
4000 -1.56299033166839
5000 -1.65731087973091
};
\addlegendentry{sampling / greedy}
\addplot [semithick, chocolate22711434, mark=triangle*, mark size=3, mark options={solid,rotate=180}]
table {%
500 14.7815081630211
1000 14.190437022149
2000 4.65225239029337
3000 2.41980211876392
4000 1.15198995775105
5000 0.562452262955945
};
\addlegendentry{SB~Linear / greedy}
\addplot [semithick, darkgoldenrod1621730, mark=x, mark size=3, mark options={solid}]
table {%
500 7.87940730996693
1000 5.1413022762499
2000 3.78305784334606
3000 2.88950041594662
4000 1.86993586804017
5000 1.23969399516029
};
\addlegendentry{CB~Linear / greedy}
\addplot [semithick, red, mark=triangle*, mark size=3, mark options={solid,rotate=180}]
table {%
500 14.6137368354404
1000 13.6952962481376
2000 5.98154579718315
3000 3.19136546034829
4000 2.17399547351373
5000 0.702381120219171
};
\addlegendentry{SB~NN / greedy}
\addplot [semithick, black, mark=x, mark size=3, mark options={solid}]
table {%
500 6.8669683799042
1000 6.90771489299845
2000 4.65493077278515
3000 3.15906781733002
4000 2.77106732399719
5000 2.6509374786656
};
\addlegendentry{CB~NN / greedy}
\end{axis}

\end{tikzpicture}

%% file: contents/conclusion.tex
\section{Conclusion}
\label{Conclusion}

In this paper, we introduced a new family of online control policies for dispatching and rebalancing decisions in \gls*{abk:amod} fleets. To derive such policies, we presented a novel \gls*{abk:CO}-enriched \gls*{abk:ML} pipeline that learns prescriptive dispatching and rebalancing decisions from full-information solutions. To establish this pipeline for our particular application case, we showed how to solve the underlying \gls*{abk:CO} problem in polynomial time, how to leverage \gls*{abk:SL} to learn a parametrization for the \gls*{abk:CO} problem such that it allows for prescriptive dispatching and rebalancing actions, and how to derive corresponding learning-based control policies.

We presented a profound numerical study based on a real-world case study for Manhattan to show the efficacy of our learning-based policies and benchmark them against current state-of-the-art approaches: a greedy policy and a sampling-based receding horizon approach \citep[cf.][]{Alonso-Mora2017a}. Our results show that the sampling policy has an unstable performance and oftentimes performs even worse than the greedy policy.
The performance of the sampling policy is in line with the performance of other benchmarks, e.g., \gls*{abk:DRL} algorithms, that reach an improvement of up to 5\% over myopic policies in similar settings. Our learning-based policies show a significant improvement over such approaches, yielding improvements over a greedy policy of up to \resultReqMAXSBGreedy\%. Moreover, our learning-based policies are robust and always outperform the greedy policy. Finally, our learning-based policies yield a good trade-off between anticorrelated key performance indicators, e.g., when maximizing the total profit or the number of served customers, we still obtain a lower distance driven per vehicle compared to policies that show worse performance on the aforementioned quantities.

\section*{Acknowledgments}
This work was supported by the German Research Foundation (DFG) under grant no. 449261765.

%% file: contents/appendix.tex

\appendix

\section{Fenchel-Young loss}
\label{appendix:fenchel-young-loss}
When summarizing the results from \cite{Berthet2020},
\begin{equation}
    F = \mathbb{E}\Big[\max_{y \in \mathcal{Y}(x_i)}(\theta + Z)^T y \Big]
\end{equation}
is strictly convex, twice differentiable, and has a gradient of
\begin{equation}
    \nabla_{\theta} F = \mathbb{E}\Big[\argmax_{y \in \mathcal{Y}(x_i)}(\theta + Z)^T y \Big].   
\end{equation}
When~$F$ is strictly convex, we can argue that~$L(\theta, x_i, y_i^*)$ is also strictly convex as it is the sum of convex functions with at least one function being strictly convex. Then we can derive the gradient from 
\begin{equation}
    L(\theta, x_i, y_i^*) = \mathbb{E}\Big[\max_{y \in \mathcal{Y}(x_i)}(\theta + Z)^T y \Big] - \theta^T y_i^*
\end{equation}
as
\begin{equation}
    \nabla_{\theta} L(\theta,x_i, y_i^*):= \mathbb{E}\Big[\argmax_{y \in \mathcal{Y}(x_i)}(\theta + Z)^T y \Big] - y_i^*.
\end{equation}

\section{Learning the ML-predictor~$\varphi_w$ for the case study}
\label{appendix:learning}
In this section, we describe how to learn the ML-predictor~$\varphi_{w}$, which predicts the weights~$\theta$ of digraph~$D$ in the \gls*{abk:ML}-layer for the real-world case study. 
In our training set, we consider the weekdays~$7^{th}$~till~$13^{th}$~of~January.
For each day, we consider the time between 9:00 AM and 10:00 AM. 
To derive the target solutions for this training set, we apply the full-information approach to the problem comprising the surrounding time period from 8:30 AM to 10:30 AM. We use this surrounding time period to provide some time for balancing the system.
Thereafter we extract the optimal full-information solution every \instanceExtractionEvery~minutes between 9:00 AM and 10:00 AM and use the solution to rebuild the corresponding \gls*{abk:CO}-layer-digraph solution. Then, the \gls*{abk:CO}-layer-digraph object~$x_i$ and the \gls*{abk:CO}-layer-digraph solution~$y^*_i$ form a training instance~$(x_i, y^*_i)$.
As we consider~5~days and a~1~hour interval with instances every~\instanceExtractionEvery~minutes we receive~\numTrainingInstances~training instances.
Moreover, we draw \numPerturbations~vectors with samples from a normal distribution~$Z \sim \mathcal{N}(0,1)$ to perturb the weight vector~$\theta$ in the loss function~\eqref{eq:loss_perturbed} which we do by perturbing parameter vector~$w$. Accordingly, evaluating the learning problem~\eqref{eq:learning_problem} in one learning iteration requires solving~\numTrainingInstances~$\times$~\numPerturbations~instances of \gls*{abk:k-dSPP}s.
With the training data set, we solve the sample average approximation of the learning problem~\eqref{eq:learning_problem} using a \gls*{abk:BFGS} algorithm.
We detail the numerical values of parameter vector~$w$ that we used in \gls*{abk:SB} and \gls*{abk:CB} policies in Appendix~\ref{appendix:Features}.

\section{The feature vector}
\label{appendix:Features}
In this section, we detail the feature vector~$\phi(\cdot)$.
We manually selected the feature vector~$\phi(\cdot)$ according to the best performance on a validation data set.
We did not apply any further feature selection process as the perturbation of the parameter~$w$ during training acts as regularization and therefore automatically sets the parameter~$w$ of irrelevant features equal to zero.
In Table~\ref{table:parameters_SB}, we present the learned parameters~$w$ for the \gls*{abk:SB} policy, and in Table~\ref{table:parameters_CB} we present the learned parameters~$w$ for the \gls*{abk:CB} policy for its respective features~$\phi(\cdot)$. Note that the respective values are only partially interpretable due to different feature scales. For the \gls*{abk:SB} policy we did not normalize features before training as we achieved the best results without normalization, but for the \gls*{abk:CB} policy we obtained the best results when normalizing features with its standard deviation.

\begingroup 
\setlength{\tabcolsep}{6pt} 
\renewcommand{\arraystretch}{0.6} 
\begin{table}[h!] 
\centering 
\tiny 
\resizebox{1\textwidth}{!}{ 
\begin{tabular}{c c c | c c c } 
\toprule
\multicolumn{6}{c}{The following parameters describe the rebalancing cell where the vehicle starts.} \\ 
& parameter description & $w_k$ & & parameter description & $w_k$\\ 
\midrule
0 & num. of vehicles & 1.091 & 4 & est. num. of future starting requests & -0.229\\
1 & num. of requests & 0.0 & 5 & est. num. of future arriving requests & 0.854\\
2 & num. of requests / num. of vehicles & -0.005 & 6 & est. num. of future starting vehicles & 0.022\\
3 & num. of vehicles / num. of requests & 0.101 & 7 & est. num. of future arriving vehicles & 0.019\\
\toprule
\multicolumn{6}{c}{The following parameters describe the rebalancing cell where the ride request starts.} \\ 
& parameter description & $w_k$ & & parameter description & $w_k$\\ 
\midrule
0 & num. of vehicles & 0.032 & 4 & est. num. of future starting requests & 0.021\\
1 & num. of requests & 0.005 & 5 & est. num. of future arriving requests & 0.02\\
2 & num. of requests / num. of vehicles & 0.001 & 6 & est. num. of future starting vehicles & 0.007\\
3 & num. of vehicles / num. of requests & 0.03 & 7 & est. num. of future arriving vehicles & 0.005\\
\toprule
\multicolumn{6}{c}{The following parameters describe the ride request.} \\ 
& parameter description & $w_k$ & & parameter description & $w_k$\\ 
\midrule
0 & duration of tour & 1.844 & 17 & est. fut. rew. at d.off loc. +[16,18)min & 0.099\\
1 & reward of tour / duration of tour & 0.007 & 18 & est. fut. rew. at d.off loc. +[18,20)min & 0.082\\
2 & distance of tour & 0.007 & 19 & est. fut. rew. at d.off loc. +[20,22)min & 0.071\\
3 & reward of tour / distance of tour & 0.883 & 20 & est. fut. rew. at d.off loc. +[22,24)min & 0.067\\
4 & reward of tour / time till pick-up & 0.295 & 21 & est. fut. rew. at d.off loc. +[24,26)min & 0.06\\
5 & reward of tour / time till drop-off & 0.006 & 22 & est. fut. rew. at d.off loc. +[26,28)min & 0.054\\
6 & cost for tour / duration of tour & 0.001 & 23 & est. fut. rew. at d.off loc. +[28,30)min & 0.05\\
7 & costs for tour / time till drop-off & 0.001 & 24 & est. fut. rew. at d.off loc. +[30,32)min & 0.047\\
8 & reward for tour & 3.115 & 25 & est. fut. rew. at d.off loc. +[32,34)min & 0.044\\
9 & est. fut. rew. at d.off loc. +[0,2)min & 0.611 & 26 & est. fut. rew. at d.off loc. +[34,36)min & 0.042\\
10 & est. fut. rew. at d.off loc. +[2,4)min & 0.853 & 27 & est. fut. rew. at d.off loc. +[36,38)min & 0.042\\
11 & est. fut. rew. at d.off loc. +[4,6)min & 0.59 & 28 & est. fut. rew. at d.off loc. +[38,40)min & 0.042\\
12 & est. fut. rew. at d.off loc. +[6,8)min & 0.301 & 29 & est. fut. rew. at d.off loc. +[40,42)min & 0.042\\
13 & est. fut. rew. at d.off loc. +[8,10)min & 0.204 & 30 & est. fut. rew. at d.off loc. +[42,44)min & 0.041\\
14 & est. fut. rew. at d.off loc. +[10,12)min & 0.17 & 31 & est. fut. rew. at d.off loc. +[44,46)min & 0.042\\
15 & est. fut. rew. at d.off loc. +[12,14)min & 0.141 & 32 & est. fut. rew. at d.off loc. +[46,48)min & 0.042\\
16 & est. fut. rew. at d.off loc. +[14,16)min & 0.116 & 33 & est. fut. rew. at d.off loc. +[48,50)min & 0.042\\
\toprule
\multicolumn{6}{c}{The following parameters describe the tour from the vehicle location to the ride request location.} \\ 
& parameter description & $w_k$ & & parameter description & $w_k$\\ 
\midrule
0 & distance to request & -0.009 & 15 & est. fut. cost at d.off loc. +[20,22)min & 0.116\\
1 & duration to request & -0.505 & 16 & est. fut. cost at d.off loc. +[22,24)min & -0.7\\
2 & cost for tour / duration of tour & 0.0 & 17 & est. fut. cost at d.off loc. +[24,26)min & -0.676\\
3 & cost for tour / time till drop-off & 0.0 & 18 & est. fut. cost at d.off loc. +[26,28)min & -0.63\\
4 & cost for tour / distance to request & 0.256 & 19 & est. fut. cost at d.off loc. +[28,30)min & -0.472\\
5 & est. fut. cost at d.off loc. +[0,2)min & 0.439 & 20 & est. fut. cost at d.off loc. +[30,32)min & -0.47\\
6 & est. fut. cost at d.off loc. +[2,4)min & -0.203 & 21 & est. fut. cost at d.off loc. +[32,34)min & -0.352\\
7 & est. fut. cost at d.off loc. +[4,6)min & -0.309 & 22 & est. fut. cost at d.off loc. +[34,36)min & -0.31\\
8 & est. fut. cost at d.off loc. +[6,8)min & -1.419 & 23 & est. fut. cost at d.off loc. +[36,38)min & -0.352\\
9 & est. fut. cost at d.off loc. +[8,10)min & -1.81 & 24 & est. fut. cost at d.off loc. +[38,40)min & -0.484\\
10 & est. fut. cost at d.off loc. +[10,12)min & -1.048 & 25 & est. fut. cost at d.off loc. +[40,42)min & -0.483\\
11 & est. fut. cost at d.off loc. +[12,14)min & -0.754 & 26 & est. fut. cost at d.off loc. +[42,44)min & -0.484\\
12 & est. fut. cost at d.off loc. +[14,16)min & -0.286 & 27 & est. fut. cost at d.off loc. +[44,46)min & -0.548\\
13 & est. fut. cost at d.off loc. +[16,18)min & 0.106 & 28 & est. fut. cost at d.off loc. +[46,48)min & -0.578\\
14 & est. fut. cost at d.off loc. +[18,20)min & 0.396 & 29 & est. fut. cost at d.off loc. +[48,50)min & -0.582\\
\toprule
\end{tabular}} 
\caption{All parameters with its respective $w$ values listed for the \gls*{abk:SB} policy when optimizing over profit; fleet size: 2000; request density: 0.3.} 
\label{table:parameters_SB} 
\end{table} 
\endgroup
\clearpage

\begingroup 
\setlength{\tabcolsep}{6pt} 
\renewcommand{\arraystretch}{0.6} 
\begin{table}[h!] 
\centering 
\tiny 
\resizebox{1\textwidth}{!}{ 
\begin{tabular}{c c c | c c c } 
\toprule
\multicolumn{6}{c}{The following parameters describe the rebalancing cell where the vehicle starts.} \\ 
& parameter description & $w_k$ & & parameter description & $w_k$\\ 
\midrule
0 & num. of vehicles & -0.048 & 4 & est. num. of future starting requests & -0.859\\
1 & num. of requests & 0.302 & 5 & est. num. of future arriving requests & 1.018\\
2 & num. of requests / num. of vehicles & -0.277 & 6 & est. num. of future starting vehicles & 0.2\\
3 & num. of vehicles / num. of requests & -0.481 & 7 & est. num. of future arriving vehicles & 0.742\\
\toprule
\multicolumn{6}{c}{The following parameters describe the rebalancing cell where the ride request starts.} \\ 
& parameter description & $w_k$ & & parameter description & $w_k$\\ 
\midrule
0 & num. of vehicles & 2.461 & 4 & est. num. of future starting requests & 3.047\\
1 & num. of requests & 3.462 & 5 & est. num. of future arriving requests & 2.629\\
2 & num. of requests / num. of vehicles & 5.924 & 6 & est. num. of future starting vehicles & 3.394\\
3 & num. of vehicles / num. of requests & 2.477 & 7 & est. num. of future arriving vehicles & 3.459\\
\toprule
\multicolumn{6}{c}{The following parameters describe the ride request.} \\ 
& parameter description & $w_k$ & & parameter description & $w_k$\\ 
\midrule
0 & duration of tour & 3.922 & 5 & reward of tour / time till drop-off & 3.577\\
1 & reward of tour / duration of tour & 3.545 & 6 & cost for tour / duration of tour & 3.784\\
2 & distance of tour & 4.117 & 7 & costs for tour / time till drop-off & 3.815\\
3 & reward of tour / distance of tour & 3.501 & 8 & reward for tour & 3.867\\
4 & reward of tour / time till pick-up & 4.956 & & & \\
\toprule
\multicolumn{6}{c}{The following parameters describe the tour from the vehicle location to the ride request location.} \\ 
& parameter description & $w_k$ & & parameter description & $w_k$\\ 
\midrule
0 & distance to request & -0.3 & 3 & cost for tour / time till drop-off & 0.022\\
1 & duration to request & -1.216 & 4 & cost for tour / distance to request & 3.473\\
2 & cost for tour / duration of tour & -0.26 & & & \\
\toprule
\multicolumn{6}{c}{The following parameters describe a rebalancing location} \\ 
& parameter description & $w_k$ & & parameter description & $w_k$\\ 
\midrule
0 & num. of vehicles & -0.077 & 7 & est. num. of future arriving vehicles & 0.378\\
1 & num. of requests & 0.088 & 8 & duration of tour & 0.511\\
2 & num. of requests / num. of vehicles & -0.362 & 9 & reward of tour / duration of tour & -1.121\\
3 & num. of vehicles / num. of requests & 0.042 & 10 & distance of tour & -0.618\\
4 & est. num. of future starting requests & 1.658 & 11 & reward of tour / distance of tour & 0.863\\
5 & est. num. of future arriving requests & -0.017 & 12 & cost for tour / duration of tour & 0.756\\
6 & est. num. of future starting vehicles & 0.133 & 13 & reward for tour & -0.468\\
\toprule
\multicolumn{6}{c}{The following parameters describe the tour to a rebalancing cell.} \\ 
& parameter description & $w_k$ & & parameter description & $w_k$\\ 
\midrule
0 & distance to location & -0.886 & 1 & duration to location & -0.347\\
\toprule
\multicolumn{6}{c}{The following parameters describe the capacity vertices and its rebalancing cell.} \\ 
& parameter description & $w_k$ & & parameter description & $w_k$\\ 
\midrule
0 & num. of vehicles & -0.24 & 8 & duration of tour & 0.706\\
1 & num. of requests & 0.243 & 9 & reward of tour / duration of tour & -1.407\\
2 & num. of requests / num. of vehicles & -0.619 & 10 & distance of tour & -0.761\\
3 & num. of vehicles / num. of requests & 0.076 & 11 & reward of tour / distance of tour & 1.208\\
4 & est. num. of future starting requests & 4.469 & 12 & cost for tour / duration of tour & 0.904\\
5 & est. num. of future arriving requests & -0.051 & 13 & reward for tour & -0.613\\
6 & est. num. of future starting vehicles & 0.244 & 14 & idx of capacity vertex in reb. cell & -50.121\\
7 & est. num. of future arriving vehicles & 0.674 & & & \\
\toprule
\end{tabular}} 
\caption{All parameters with its respective $w$ values listed for the \gls*{abk:CB} policy when optimizing over profit; fleet size: 2000; request density: 0.3.} 
\label{table:parameters_CB} 
\end{table} 
\endgroup 
\clearpage

\section{Sparsification of the full-information digraph}
\label{appendix:Sparsification}

\begin{figure}[H]
    \begin{minipage}[t]{0.48\linewidth}
      \centering
      \resizebox{\textwidth}{!}{
      \input{./figures/evaluate_heuristic_distance_range_profit_absolute.tex}}
      \caption{Spatial sparsification when optimizing the full-information bound over profit; fleet size: 2000; request density: 0.3.}
    \label{fig:heuristics_spatial}
    \end{minipage}
    \hspace{0.1cm}
    \begin{minipage}[t]{0.48\linewidth}
      \centering
      \resizebox{\textwidth}{!}{
      \input{./figures/evaluate_heuristic_time_range_profit_absolute.tex}}
      \caption{Temporal sparsification when optimizing the full-information bound over profit; fleet size: 2000; request density: 0.3.}
    \label{fig:heuristics_temporal}
    \end{minipage} 
\end{figure}
The \gls*{abk:CO}-layer-problem is solvable in polynomial time as a \gls*{abk:k-dSPP}. Still, the computation time in real-world scenarios is crucial as the central controller only has limited time to make dispatching and rebalancing decisions in each system time period. Therefore, we sparsify the \gls*{abk:CO}-layer-digraph~$D$ using temporal and spatial sparsification.

For temporal sparsification, we exclude arcs from the \gls*{abk:CO}-layer-digraph~$D$ if a vehicle has a downtime greater than~$t^{\text{max}}$ seconds, which is the case when a previous request~$r_i$ and the subsequent request~$r_j$ are more than~$t^{\text{max}}$ apart from each other,~$(s_{r_j} - a_{r_i} > t^{\text{max}})$, or if a vehicle~$v$ has a downtime larger than~$t^{max}$ to the first request~$r_i$,~$(s_{r_i} - t_{v} > t^{\text{max}})$. We call this a temporal cut.

For spatial sparsification, we exclude arcs from the \gls*{abk:CO}-layer-digraph~$D$ if a vehicle has to drive a distance further than~$d^{\text{max}}$ to the next request, which is the case when a previous request~$r_i$ and the subsequent request~$r_j$ have a distance greater than~$d^{\text{max}}$,~$\tau (d_{r_i}, o_{r_j}) > d^{\text{max}}$, or if a vehicle~$v$ has to drive a distance greater than~$d^{\text{max}}$ to the first request~$r_i$,~$\tau (l_v, o_{r_i}) > d^{\text{max}}$. We call this a spatial cut.

Figure~\ref{fig:heuristics_spatial} and Figure~\ref{fig:heuristics_temporal} show that the sparsification of the graph with a spatial cut at \distanceSparsification~kilometers and a temporal cut at \temporalSparsification~seconds has no impact on the performance of the algorithm. We use these spatial and temporal cut values in our case study.
\clearpage

\section{Hyperparameter selection for the sampling policy}
\label{appendix:hyper_samplingBench}
\begingroup
\renewcommand\baselinestretch{0.5}
\begin{figure}[H] 
    \begin{minipage}[t]{0.48\linewidth}
        \centering
        \resizebox{\textwidth}{!}{ 
        \input{./figures/evaluate_extended_horizon_profit_absolute.tex}}
        \caption{Performance of the sampling policy for different extended horizons when optimizing over profit. Fleet size: 2000; request density: 0.3.}
        \label{fig:hyper_samplingBench_profit} 
    \end{minipage} 
    \hspace{0.1cm}
    \begin{minipage}[t]{0.48\linewidth}
        \centering
        \resizebox{\textwidth}{!}{
        \input{./figures/evaluate_reduce_weight_of_future_requests_factor_profit_absolute.tex}}
        \caption{Performance of the sampling policy for different discount factors when optimizing over profit. Fleet size: 2000; request density: 0.3.}
        \label{fig:hyper_samplingReduceWeights_profit} 
    \end{minipage}  
\end{figure}
\endgroup
The performance of the sampling policy depends on the sampled requests in the prediction horizon. Therefore we investigate different lengths of the prediction horizon to determine the best prediction horizon length (see Figure~\ref{fig:hyper_samplingBench_profit}). We choose a prediction horizon of \predictionHorizonSampling~minutes as we achieved the best results with this setting on a validation data set. We recall that the prediction horizon is the extended horizon minus the system time period of \systemTimePeriod~minute.

Moreover, the sampling policy improves performance when preferring requests in the system time period over requests in the prediction horizon according to \cite{Alonso-Mora2017a}. Therefore we test the sampling policy over different discount factors for sampled requests (Figure~\ref{fig:hyper_samplingReduceWeights_profit}). The discount factor is a multiplication factor that we use to reduce weights on arcs to sampled requests in the prediction horizon. We choose a discount factor of~\reductionFactorSampling~as we received the best results for this setting on a validation data set.

\clearpage

\section{Hyperparameter selection for the variants SP and CP of the CO-enriched ML pipeline}
\label{appendix:hyper_idea}
\begin{figure}[H]
\renewcommand\baselinestretch{0.5}
    \begin{minipage}[t]{0.32\linewidth}
        \centering
        \resizebox{\textwidth}{!}{
        \input{./figures/evaluate_hyperparameterSelection_extended_horizon_profit_absolute_policy_SB.tex}}
        \caption{\gls*{abk:SB} policy for different extended horizons when optimizing over profit. Fleet size: 2000; request density: 0.3.}
    \label{fig:hyper_idea_SP}
    \end{minipage}
    \hfill
    \begin{minipage}[t]{0.32\linewidth}
        \centering
        \resizebox{\textwidth}{!}{
        \input{./figures/evaluate_hyperparameterSelection_extended_horizon_profit_absolute_policy_CB.tex}}
        \caption{\gls*{abk:CB} policy for different extended horizons when optimizing over profit. Fleet size: 2000; request density: 0.3.}
    \label{fig:hyper_idea_CP}
    \end{minipage}
    \hfill
    \begin{minipage}[t]{0.32\linewidth}
        \centering
        \resizebox{\textwidth}{!}{
        \input{./figures/evaluate_fix_capacity_profit_absolute.tex}}
        \caption{\gls*{abk:CB} policy for different number of capacity vertices in rebalancing cells. Fleet size: 2000; request density: 0.3.}
    \label{fig:maxcap_CP}
    \end{minipage}
\end{figure}
In this section, we determine the hyperparameters for the \gls*{abk:CO}-enriched \gls*{abk:ML} pipeline variants, namely the \gls*{abk:SB} and \gls*{abk:CB} policies. We tested the \gls*{abk:SB} policy (Figure~\ref{fig:hyper_idea_SP}) and the \gls*{abk:CB} policy (Figure~\ref{fig:hyper_idea_CP}) with different prediction horizons. For the \gls*{abk:SB} and \gls*{abk:CB} policy a prediction horizon of \predictionHorizonSB~minutes performs best when testing on a validation data set.
We recall that the prediction horizon is the extended horizon minus the system time period of \systemTimePeriod~minute.
We tested the \gls*{abk:CB} policy with different amounts of capacity vertices for each rebalancing cell (Figure~\ref{fig:maxcap_CP}). We retrieved the best performance with \maxCapacityCB~capacity vertex for each rebalancing cell when testing on a validation data set.

\section{Comparing the performance of the deep learning ML-layer and the linear ML-layer}
\label{appendix:LinearNNOverfitting}
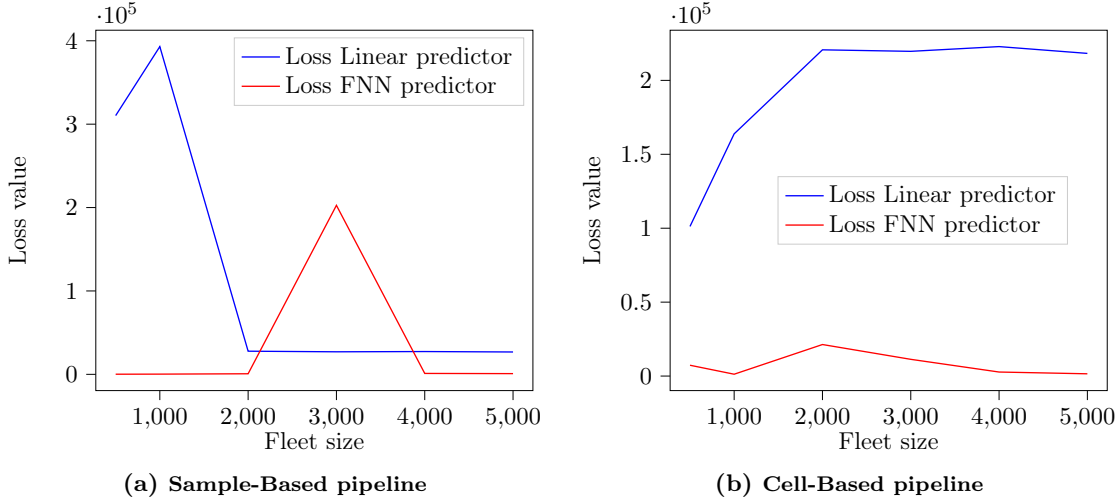
\begin{figure}[H]
\begin{subfigure}[t]{0.45\textwidth}
      \centering
      \resizebox{1\textwidth}{!}{
      \input{./figures/comparisonLossLinearNNpolicy_SB.tex}}
      \caption{Sample-Based pipeline}
\end{subfigure}
\begin{subfigure}[t]{0.45\textwidth}
      \centering
      \resizebox{1\textwidth}{!}{
      \input{./figures/comparisonLossLinearNNpolicy_CB.tex}}
      \caption{Cell-Based pipeline}
\end{subfigure}
    \caption{Comparison of loss value with \gls{abk:FNN} predictor and linear predictor for different flee sizes.}
    \label{fig:comparisonLoss}
\end{figure}

In Figure~\ref{fig:comparisonLoss} we compare the Fenchel-Young loss values of the \gls*{abk:CO}-enriched \gls*{abk:ML} pipeline at the end of training for an \gls{abk:FNN} predictor and a linear predictor in the \gls*{abk:ML}-layer.
Apart from the outlier in the \gls*{abk:SB} setting with a fleet size of~3000 for which the \gls{abk:FNN} has not converged after~70~hours of training, we can see that the loss is always smaller for the \gls{abk:FNN} predictor.
The lower loss value when using the \gls{abk:FNN} predictor shows that, as expected, the \gls{abk:FNN} better approximates the true cost function of the \gls*{abk:CO}-layer-digraph and therefore leads to better results on almost all scenarios (Figure \ref{fig:comparisonNN}).
 
\begin{figure}[H]
\begin{subfigure}[t]{0.45\textwidth}
      \centering
      \resizebox{1\textwidth}{!}{
      \input{./figures/comparisonEqualityLinearNNpolicy_SB.tex}}
      \caption{Sample-Based pipeline}
\end{subfigure}
\begin{subfigure}[t]{0.45\textwidth}
      \centering
      \resizebox{1\textwidth}{!}{
      \input{./figures/comparisonEqualityLinearNNpolicy_CB.tex}}
      \caption{Cell-Based pipeline}
\end{subfigure}
    \caption{Comparison of equality between $\hat{y_i}$ and $y_i^{\star}$ on all instances~$i$ in training set~$(x_1,y_1^*),\ldots,(x_n,y_n^*)$.}
    \label{fig:comparisonEquality}
\end{figure}
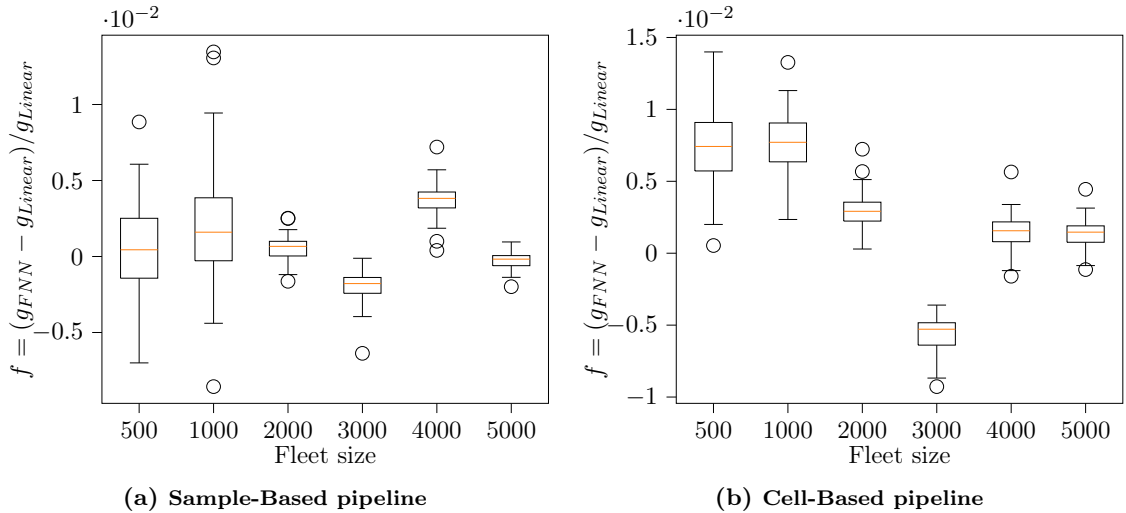

However, the better loss value does not necessarily lead to better generalization: in some scenarios, the \gls*{abk:CO}-enriched \gls*{abk:ML} pipeline using an \gls{abk:FNN} predictor does not outperform the pipeline using the linear predictor~(Figure \ref{fig:comparisonNN}). Figure~\ref{fig:comparisonEquality} shows why this is the case by comparing the difference between $\hat{y_i}$ and $y_i^{\star}$ on the training set~$(x_1,y_1^*),\ldots,(x_n,y_n^*)$. Here, $g_{FNN} = \sum(\hat{y}_{i,FNN} == y_i^{\star})$ and $g_{Linear} = \sum(\hat{y}_{i,Linear} == y_i^{\star})$ indicate the number of matching edges in the predicted solution~$\hat{y}_i$ and the true solution~$y_i^{\star}$, using the \gls{abk:FNN} predictor and the linear predictor, respectively. In scenarios where the \gls*{abk:CO}-enriched \gls*{abk:ML} pipeline using the \gls{abk:FNN} predictor performs worse than the pipeline using the linear predictor (Figure~\ref{fig:comparisonNN}) during testing, the \gls{abk:FNN} predictor sill outperforms the linear model ($g_{FNN} > g_{Linear}$) on the training set. This indicates that the \gls{abk:FNN} predictor is overfitting in these scenarios.

%% file: figures/evaluate_heuristic_distance_range_profit_absolute.tex
\begin{tikzpicture}

\definecolor{darkorange25512714}{RGB}{255,127,14}
\definecolor{darkslategray38}{RGB}{38,38,38}
\definecolor{lightgray204}{RGB}{204,204,204}
\definecolor{steelblue31119180}{RGB}{31,119,180}

\begin{axis}[
axis line style={lightgray204},
legend style={
  fill opacity=0.8,
  draw opacity=1,
  text opacity=1,
  at={(0.03,0.97)},
  anchor=north west,
  draw=lightgray204
},
scaled y ticks=true,
tick align=outside,
x grid style={lightgray204},
xlabel style={below=6mm},
xlabel={Maximum distance to next request [km]},
xmajorgrids,
xmajorticks=false,
xmajorticks=true,
xmin=-0.5, xmax=7.5,
xtick style={color=darkslategray38},
xtick style={draw=none},
xtick={0,1,2,3,4,5,6,7},
xticklabel style={rotate=90.0},
xticklabels={0.3,0.5,0.7,0.9,1.1,1.3,1.5,1.7},
y grid style={lightgray204},
ylabel style={align=center},
ylabel={Profit [\$]},
ymajorgrids,
ymajorticks=false,
ymajorticks=true,
ymin=30054.9255696242, ymax=75601.3279167477,
ytick style={color=darkslategray38},
ytick style={draw=none}
]
\path [draw=black, fill=steelblue31119180]
(axis cs:-0.25,52207.5192598186)
--(axis cs:0.25,52207.5192598186)
--(axis cs:0.25,56755.9597049168)
--(axis cs:-0.25,56755.9597049168)
--(axis cs:-0.25,52207.5192598186)
--cycle;
\addplot [black, forget plot]
table {%
0 52207.5192598186
0 52207.5192598186
};
\addplot [black, forget plot]
table {%
0 56755.9597049168
0 57828.4388877406
};
\addplot [black, forget plot]
table {%
-0.125 52207.5192598186
0.125 52207.5192598186
};
\addplot [black, forget plot]
table {%
-0.125 57828.4388877406
0.125 57828.4388877406
};
\addplot [black, mark=o, mark size=3, mark options={solid,fill opacity=0}, only marks, forget plot]
table {%
0 32125.2165854026
0 32125.2165854026
0 32125.2165854026
0 32125.2165854026
0 32125.2165854026
0 32125.2165854026
0 32125.2165854026
0 32125.2165854026
};
\path [draw=black, fill=steelblue31119180]
(axis cs:0.75,61431.8910744907)
--(axis cs:1.25,61431.8910744907)
--(axis cs:1.25,65857.6474690116)
--(axis cs:0.75,65857.6474690116)
--(axis cs:0.75,61431.8910744907)
--cycle;
\addplot [black, forget plot]
table {%
1 61431.8910744907
1 61431.8910744907
};
\addplot [black, forget plot]
table {%
1 65857.6474690116
1 67626.7651257717
};
\addplot [black, forget plot]
table {%
0.875 61431.8910744907
1.125 61431.8910744907
};
\addplot [black, forget plot]
table {%
0.875 67626.7651257717
1.125 67626.7651257717
};
\addplot [black, mark=o, mark size=3, mark options={solid,fill opacity=0}, only marks, forget plot]
table {%
1 35752.4818487034
1 35752.4818487034
1 35752.4818487034
1 35752.4818487034
1 35752.4818487034
1 35752.4818487034
1 35752.4818487034
1 35752.4818487034
};
\path [draw=black, fill=steelblue31119180]
(axis cs:1.75,65145.6245877649)
--(axis cs:2.25,65145.6245877649)
--(axis cs:2.25,69157.3050238322)
--(axis cs:1.75,69157.3050238322)
--(axis cs:1.75,65145.6245877649)
--cycle;
\addplot [black, forget plot]
table {%
2 65145.6245877649
2 65145.6245877649
};
\addplot [black, forget plot]
table {%
2 69157.3050238322
2 71631.5360484503
};
\addplot [black, forget plot]
table {%
1.875 65145.6245877649
2.125 65145.6245877649
};
\addplot [black, forget plot]
table {%
1.875 71631.5360484503
2.125 71631.5360484503
};
\addplot [black, mark=o, mark size=3, mark options={solid,fill opacity=0}, only marks, forget plot]
table {%
2 36101.2256390214
2 36101.2256390214
2 36101.2256390214
2 36101.2256390214
2 36101.2256390214
2 36101.2256390214
2 36101.2256390214
2 36101.2256390214
};
\path [draw=black, fill=steelblue31119180]
(axis cs:2.75,66196.2979919029)
--(axis cs:3.25,66196.2979919029)
--(axis cs:3.25,70385.4119085169)
--(axis cs:2.75,70385.4119085169)
--(axis cs:2.75,66196.2979919029)
--cycle;
\addplot [black, forget plot]
table {%
3 66196.2979919029
3 66196.2979919029
};
\addplot [black, forget plot]
table {%
3 70385.4119085169
3 73070.4636204806
};
\addplot [black, forget plot]
table {%
2.875 66196.2979919029
3.125 66196.2979919029
};
\addplot [black, forget plot]
table {%
2.875 73070.4636204806
3.125 73070.4636204806
};
\addplot [black, mark=o, mark size=3, mark options={solid,fill opacity=0}, only marks, forget plot]
table {%
3 36139.5869836973
3 36139.5869836973
3 36139.5869836973
3 36139.5869836973
3 36139.5869836973
3 36139.5869836973
3 36139.5869836973
3 36139.5869836973
};
\path [draw=black, fill=steelblue31119180]
(axis cs:3.75,66353.9195416992)
--(axis cs:4.25,66353.9195416992)
--(axis cs:4.25,70714.2731650138)
--(axis cs:3.75,70714.2731650138)
--(axis cs:3.75,66353.9195416992)
--cycle;
\addplot [black, forget plot]
table {%
4 66353.9195416992
4 66353.9195416992
};
\addplot [black, forget plot]
table {%
4 70714.2731650138
4 73457.0152936829
};
\addplot [black, forget plot]
table {%
3.875 66353.9195416992
4.125 66353.9195416992
};
\addplot [black, forget plot]
table {%
3.875 73457.0152936829
4.125 73457.0152936829
};
\addplot [black, mark=o, mark size=3, mark options={solid,fill opacity=0}, only marks, forget plot]
table {%
4 36143.4786013127
4 36143.4786013127
4 36143.4786013127
4 36143.4786013127
4 36143.4786013127
4 36143.4786013127
4 36143.4786013127
4 36143.4786013127
};
\path [draw=black, fill=steelblue31119180]
(axis cs:4.75,66370.0363311033)
--(axis cs:5.25,66370.0363311033)
--(axis cs:5.25,70764.8608313665)
--(axis cs:4.75,70764.8608313665)
--(axis cs:4.75,66370.0363311033)
--cycle;
\addplot [black, forget plot]
table {%
5 66370.0363311033
5 66370.0363311033
};
\addplot [black, forget plot]
table {%
5 70764.8608313665
5 73479.3107524759
};
\addplot [black, forget plot]
table {%
4.875 66370.0363311033
5.125 66370.0363311033
};
\addplot [black, forget plot]
table {%
4.875 73479.3107524759
5.125 73479.3107524759
};
\addplot [black, mark=o, mark size=3, mark options={solid,fill opacity=0}, only marks, forget plot]
table {%
5 36149.5620377344
5 36149.5620377344
5 36149.5620377344
5 36149.5620377344
5 36149.5620377344
5 36149.5620377344
5 36149.5620377344
5 36149.5620377344
};
\path [draw=black, fill=steelblue31119180]
(axis cs:5.75,66392.8959064822)
--(axis cs:6.25,66392.8959064822)
--(axis cs:6.25,70768.9185790466)
--(axis cs:5.75,70768.9185790466)
--(axis cs:5.75,66392.8959064822)
--cycle;
\addplot [black, forget plot]
table {%
6 66392.8959064822
6 66392.8959064822
};
\addplot [black, forget plot]
table {%
6 70768.9185790466
6 73508.3859868621
};
\addplot [black, forget plot]
table {%
5.875 66392.8959064822
6.125 66392.8959064822
};
\addplot [black, forget plot]
table {%
5.875 73508.3859868621
6.125 73508.3859868621
};
\addplot [black, mark=o, mark size=3, mark options={solid,fill opacity=0}, only marks, forget plot]
table {%
6 36149.9816426179
6 36149.9816426179
6 36149.9816426179
6 36149.9816426179
6 36149.9816426179
6 36149.9816426179
6 36149.9816426179
6 36149.9816426179
};
\path [draw=black, fill=steelblue31119180]
(axis cs:6.75,66413.6040664385)
--(axis cs:7.25,66413.6040664385)
--(axis cs:7.25,70770.8697506807)
--(axis cs:6.75,70770.8697506807)
--(axis cs:6.75,66413.6040664385)
--cycle;
\addplot [black, forget plot]
table {%
7 66413.6040664385
7 66413.6040664385
};
\addplot [black, forget plot]
table {%
7 70770.8697506807
7 73531.0369009694
};
\addplot [black, forget plot]
table {%
6.875 66413.6040664385
7.125 66413.6040664385
};
\addplot [black, forget plot]
table {%
6.875 73531.0369009694
7.125 73531.0369009694
};
\addplot [black, mark=o, mark size=3, mark options={solid,fill opacity=0}, only marks, forget plot]
table {%
7 36150.0764804787
7 36150.0764804787
7 36150.0764804787
7 36150.0764804787
7 36150.0764804787
7 36150.0764804787
7 36150.0764804787
7 36150.0764804787
};
\addplot [darkorange25512714, forget plot]
table {%
-0.25 54988.4217218814
0.25 54988.4217218814
};
\addplot [darkorange25512714, forget plot]
table {%
0.75 63895.6330510536
1.25 63895.6330510536
};
\addplot [darkorange25512714, forget plot]
table {%
1.75 66301.3629379395
2.25 66301.3629379395
};
\addplot [darkorange25512714, forget plot]
table {%
2.75 67526.8039538576
3.25 67526.8039538576
};
\addplot [darkorange25512714, forget plot]
table {%
3.75 68249.251939201
4.25 68249.251939201
};
\addplot [darkorange25512714, forget plot]
table {%
4.75 68439.4692597291
5.25 68439.4692597291
};
\addplot [darkorange25512714, forget plot]
table {%
5.75 68448.5121459573
6.25 68448.5121459573
};
\addplot [darkorange25512714, forget plot]
table {%
6.75 68453.8026420818
7.25 68453.8026420818
};
\end{axis}

\end{tikzpicture}

%% file: figures/evaluate_heuristic_time_range_profit_absolute.tex
\begin{tikzpicture}

\definecolor{darkorange25512714}{RGB}{255,127,14}
\definecolor{darkslategray38}{RGB}{38,38,38}
\definecolor{lightgray204}{RGB}{204,204,204}
\definecolor{steelblue31119180}{RGB}{31,119,180}

\begin{axis}[
axis line style={lightgray204},
legend style={fill opacity=0.8, draw opacity=1, text opacity=1, draw=lightgray204},
scaled y ticks=true,
tick align=outside,
x grid style={lightgray204},
xlabel style={below=6mm},
xlabel={Maximum time to next request [sec]},
xmajorgrids,
xmajorticks=false,
xmajorticks=true,
xmin=-0.5, xmax=5.5,
xtick style={color=darkslategray38},
xtick style={draw=none},
xtick={0,1,2,3,4,5},
xticklabel style={rotate=90.0},
xticklabels={120.0,180.0,360.0,540.0,720.0,900.0},
y grid style={lightgray204},
ylabel style={align=center},
ylabel={Profit [\$]},
ymajorgrids,
ymajorticks=false,
ymajorticks=true,
ymin=13595.3863558897, ymax=76443.7670832316,
ytick style={color=darkslategray38},
ytick style={draw=none}
]
\path [draw=black, fill=steelblue31119180]
(axis cs:-0.25,37252.8965276867)
--(axis cs:0.25,37252.8965276867)
--(axis cs:0.25,41523.0424373222)
--(axis cs:-0.25,41523.0424373222)
--(axis cs:-0.25,37252.8965276867)
--cycle;
\addplot [black, forget plot]
table {%
0 37252.8965276867
0 37252.8965276867
};
\addplot [black, forget plot]
table {%
0 41523.0424373222
0 42781.153829374
};
\addplot [black, forget plot]
table {%
-0.125 37252.8965276867
0.125 37252.8965276867
};
\addplot [black, forget plot]
table {%
-0.125 42781.153829374
0.125 42781.153829374
};
\addplot [black, mark=o, mark size=3, mark options={solid,fill opacity=0}, only marks, forget plot]
table {%
0 16452.1309344052
0 16452.1309344052
0 16452.1309344052
0 16452.1309344052
0 16452.1309344052
0 16452.1309344052
};
\path [draw=black, fill=steelblue31119180]
(axis cs:0.75,58902.0829567758)
--(axis cs:1.25,58902.0829567758)
--(axis cs:1.25,62310.0600371106)
--(axis cs:0.75,62310.0600371106)
--(axis cs:0.75,58902.0829567758)
--cycle;
\addplot [black, forget plot]
table {%
1 58902.0829567758
1 58902.0829567758
};
\addplot [black, forget plot]
table {%
1 62310.0600371106
1 62359.7240589399
};
\addplot [black, forget plot]
table {%
0.875 58902.0829567758
1.125 58902.0829567758
};
\addplot [black, forget plot]
table {%
0.875 62359.7240589399
1.125 62359.7240589399
};
\addplot [black, mark=o, mark size=3, mark options={solid,fill opacity=0}, only marks, forget plot]
table {%
1 29826.8042933334
1 29826.8042933334
1 29826.8042933334
1 29826.8042933334
1 29826.8042933334
1 29826.8042933334
};
\path [draw=black, fill=steelblue31119180]
(axis cs:1.75,65915.9008803178)
--(axis cs:2.25,65915.9008803178)
--(axis cs:2.25,70306.2520233265)
--(axis cs:1.75,70306.2520233265)
--(axis cs:1.75,65915.9008803178)
--cycle;
\addplot [black, forget plot]
table {%
2 65915.9008803178
2 65915.9008803178
};
\addplot [black, forget plot]
table {%
2 70306.2520233265
2 72736.6939531337
};
\addplot [black, forget plot]
table {%
1.875 65915.9008803178
2.125 65915.9008803178
};
\addplot [black, forget plot]
table {%
1.875 72736.6939531337
2.125 72736.6939531337
};
\addplot [black, mark=o, mark size=3, mark options={solid,fill opacity=0}, only marks, forget plot]
table {%
2 35581.7825953578
2 35581.7825953578
2 35581.7825953578
2 35581.7825953578
2 35581.7825953578
2 35581.7825953578
};
\path [draw=black, fill=steelblue31119180]
(axis cs:2.75,66257.4230130377)
--(axis cs:3.25,66257.4230130377)
--(axis cs:3.25,70614.5662154435)
--(axis cs:2.75,70614.5662154435)
--(axis cs:2.75,66257.4230130377)
--cycle;
\addplot [black, forget plot]
table {%
3 66257.4230130377
3 66257.4230130377
};
\addplot [black, forget plot]
table {%
3 70614.5662154435
3 73415.6091609932
};
\addplot [black, forget plot]
table {%
2.875 66257.4230130377
3.125 66257.4230130377
};
\addplot [black, forget plot]
table {%
2.875 73415.6091609932
3.125 73415.6091609932
};
\addplot [black, mark=o, mark size=3, mark options={solid,fill opacity=0}, only marks, forget plot]
table {%
3 35808.3855116191
3 35808.3855116191
3 35808.3855116191
3 35808.3855116191
3 35808.3855116191
3 35808.3855116191
};
\path [draw=black, fill=steelblue31119180]
(axis cs:3.75,66376.7017478538)
--(axis cs:4.25,66376.7017478538)
--(axis cs:4.25,70720.9550403843)
--(axis cs:3.75,70720.9550403843)
--(axis cs:3.75,66376.7017478538)
--cycle;
\addplot [black, forget plot]
table {%
4 66376.7017478538
4 66376.7017478538
};
\addplot [black, forget plot]
table {%
4 70720.9550403843
4 73547.7028786351
};
\addplot [black, forget plot]
table {%
3.875 66376.7017478538
4.125 66376.7017478538
};
\addplot [black, forget plot]
table {%
3.875 73547.7028786351
4.125 73547.7028786351
};
\addplot [black, mark=o, mark size=3, mark options={solid,fill opacity=0}, only marks, forget plot]
table {%
4 35936.013933633
4 35936.013933633
4 35936.013933633
4 35936.013933633
4 35936.013933633
4 35936.013933633
};
\path [draw=black, fill=steelblue31119180]
(axis cs:4.75,66410.5251036266)
--(axis cs:5.25,66410.5251036266)
--(axis cs:5.25,70751.1400095557)
--(axis cs:4.75,70751.1400095557)
--(axis cs:4.75,66410.5251036266)
--cycle;
\addplot [black, forget plot]
table {%
5 66410.5251036266
5 66410.5251036266
};
\addplot [black, forget plot]
table {%
5 70751.1400095557
5 73587.022504716
};
\addplot [black, forget plot]
table {%
4.875 66410.5251036266
5.125 66410.5251036266
};
\addplot [black, forget plot]
table {%
4.875 73587.022504716
5.125 73587.022504716
};
\addplot [black, mark=o, mark size=3, mark options={solid,fill opacity=0}, only marks, forget plot]
table {%
5 36018.2213686719
5 36018.2213686719
5 36018.2213686719
5 36018.2213686719
5 36018.2213686719
5 36018.2213686719
};
\addplot [darkorange25512714, forget plot]
table {%
-0.25 39128.2206275796
0.25 39128.2206275796
};
\addplot [darkorange25512714, forget plot]
table {%
0.75 59370.6088718419
1.25 59370.6088718419
};
\addplot [darkorange25512714, forget plot]
table {%
1.75 67891.8360448031
2.25 67891.8360448031
};
\addplot [darkorange25512714, forget plot]
table {%
2.75 68248.567276309
3.25 68248.567276309
};
\addplot [darkorange25512714, forget plot]
table {%
3.75 68342.6942011278
4.25 68342.6942011278
};
\addplot [darkorange25512714, forget plot]
table {%
4.75 68372.0303663439
5.25 68372.0303663439
};
\end{axis}

\end{tikzpicture}

%% file: figures/evaluate_extended_horizon_profit_absolute.tex
\begin{tikzpicture}

\definecolor{darkorange25512714}{RGB}{255,127,14}
\definecolor{darkslategray38}{RGB}{38,38,38}
\definecolor{lightgray204}{RGB}{204,204,204}
\definecolor{steelblue31119180}{RGB}{31,119,180}

\begin{axis}[
axis line style={lightgray204},
legend style={fill opacity=0.8, draw opacity=1, text opacity=1, draw=lightgray204},
scaled y ticks=true,
tick align=outside,
x grid style={lightgray204},
xlabel style={below=12mm},
xlabel={Extended horizon [hour:min:sec]},
xmajorgrids,
xmajorticks=false,
xmajorticks=true,
xmin=-0.5, xmax=5.5,
xtick style={color=darkslategray38},
xtick style={draw=none},
xtick={0,1,2,3,4,5},
xticklabel style={rotate=90.0},
xticklabels={00:05:00,00:10:00,00:15:00,00:20:00,00:25:00,00:30:00},
y grid style={lightgray204},
ylabel style={align=center},
ylabel={Profit [\$]},
ymajorgrids,
ymajorticks=false,
ymajorticks=true,
ymin=18640.9484589367, ymax=43379.9801868768,
ytick style={color=darkslategray38},
ytick style={draw=none}
]
\path [draw=black, fill=steelblue31119180]
(axis cs:-0.25,38783.5472112402)
--(axis cs:0.25,38783.5472112402)
--(axis cs:0.25,41998.3510159261)
--(axis cs:-0.25,41998.3510159261)
--(axis cs:-0.25,38783.5472112402)
--cycle;
\addplot [black, forget plot]
table {%
0 38783.5472112402
0 38783.5472112402
};
\addplot [black, forget plot]
table {%
0 41998.3510159261
0 42255.4787446977
};
\addplot [black, forget plot]
table {%
-0.125 38783.5472112402
0.125 38783.5472112402
};
\addplot [black, forget plot]
table {%
-0.125 42255.4787446977
0.125 42255.4787446977
};
\addplot [black, mark=o, mark size=3, mark options={solid,fill opacity=0}, only marks, forget plot]
table {%
0 23225.0952357815
};
\path [draw=black, fill=steelblue31119180]
(axis cs:0.75,37166.5581241198)
--(axis cs:1.25,37166.5581241198)
--(axis cs:1.25,40160.1145527217)
--(axis cs:0.75,40160.1145527217)
--(axis cs:0.75,37166.5581241198)
--cycle;
\addplot [black, forget plot]
table {%
1 37166.5581241198
1 37166.5581241198
};
\addplot [black, forget plot]
table {%
1 40160.1145527217
1 40231.5880680933
};
\addplot [black, forget plot]
table {%
0.875 37166.5581241198
1.125 37166.5581241198
};
\addplot [black, forget plot]
table {%
0.875 40231.5880680933
1.125 40231.5880680933
};
\addplot [black, mark=o, mark size=3, mark options={solid,fill opacity=0}, only marks, forget plot]
table {%
1 21200.1876107291
};
\path [draw=black, fill=steelblue31119180]
(axis cs:1.75,35840.2961676439)
--(axis cs:2.25,35840.2961676439)
--(axis cs:2.25,39030.7068256501)
--(axis cs:1.75,39030.7068256501)
--(axis cs:1.75,35840.2961676439)
--cycle;
\addplot [black, forget plot]
table {%
2 35840.2961676439
2 35840.2961676439
};
\addplot [black, forget plot]
table {%
2 39030.7068256501
2 39043.8373769016
};
\addplot [black, forget plot]
table {%
1.875 35840.2961676439
2.125 35840.2961676439
};
\addplot [black, forget plot]
table {%
1.875 39043.8373769016
2.125 39043.8373769016
};
\addplot [black, mark=o, mark size=3, mark options={solid,fill opacity=0}, only marks, forget plot]
table {%
2 20654.7108001617
};
\path [draw=black, fill=steelblue31119180]
(axis cs:2.75,35214.1076863684)
--(axis cs:3.25,35214.1076863684)
--(axis cs:3.25,38458.846067493)
--(axis cs:2.75,38458.846067493)
--(axis cs:2.75,35214.1076863684)
--cycle;
\addplot [black, forget plot]
table {%
3 35214.1076863684
3 35214.1076863684
};
\addplot [black, forget plot]
table {%
3 38458.846067493
3 38692.389088338
};
\addplot [black, forget plot]
table {%
2.875 35214.1076863684
3.125 35214.1076863684
};
\addplot [black, forget plot]
table {%
2.875 38692.389088338
3.125 38692.389088338
};
\addplot [black, mark=o, mark size=3, mark options={solid,fill opacity=0}, only marks, forget plot]
table {%
3 20078.6569655261
};
\path [draw=black, fill=steelblue31119180]
(axis cs:3.75,35456.8132671539)
--(axis cs:4.25,35456.8132671539)
--(axis cs:4.25,38028.0175072808)
--(axis cs:3.75,38028.0175072808)
--(axis cs:3.75,35456.8132671539)
--cycle;
\addplot [black, forget plot]
table {%
4 35456.8132671539
4 35456.8132671539
};
\addplot [black, forget plot]
table {%
4 38028.0175072808
4 38110.8432047728
};
\addplot [black, forget plot]
table {%
3.875 35456.8132671539
4.125 35456.8132671539
};
\addplot [black, forget plot]
table {%
3.875 38110.8432047728
4.125 38110.8432047728
};
\addplot [black, mark=o, mark size=3, mark options={solid,fill opacity=0}, only marks, forget plot]
table {%
4 20048.2738925235
};
\path [draw=black, fill=steelblue31119180]
(axis cs:4.75,34856.8666526468)
--(axis cs:5.25,34856.8666526468)
--(axis cs:5.25,37894.7654029047)
--(axis cs:4.75,37894.7654029047)
--(axis cs:4.75,34856.8666526468)
--cycle;
\addplot [black, forget plot]
table {%
5 34856.8666526468
5 34856.8666526468
};
\addplot [black, forget plot]
table {%
5 37894.7654029047
5 38384.4556252295
};
\addplot [black, forget plot]
table {%
4.875 34856.8666526468
5.125 34856.8666526468
};
\addplot [black, forget plot]
table {%
4.875 38384.4556252295
5.125 38384.4556252295
};
\addplot [black, mark=o, mark size=3, mark options={solid,fill opacity=0}, only marks, forget plot]
table {%
5 19765.4499011158
};
\addplot [darkorange25512714, forget plot]
table {%
-0.25 39651.2157652232
0.25 39651.2157652232
};
\addplot [darkorange25512714, forget plot]
table {%
0.75 38121.2876697168
1.25 38121.2876697168
};
\addplot [darkorange25512714, forget plot]
table {%
1.75 37031.5957771412
2.25 37031.5957771412
};
\addplot [darkorange25512714, forget plot]
table {%
2.75 36504.8990486021
3.25 36504.8990486021
};
\addplot [darkorange25512714, forget plot]
table {%
3.75 36246.8706748777
4.25 36246.8706748777
};
\addplot [darkorange25512714, forget plot]
table {%
4.75 36349.2166932611
5.25 36349.2166932611
};
\end{axis}

\end{tikzpicture}

%% file: figures/evaluate_reduce_weight_of_future_requests_factor_profit_absolute.tex
\begin{tikzpicture}

\definecolor{darkorange25512714}{RGB}{255,127,14}
\definecolor{darkslategray38}{RGB}{38,38,38}
\definecolor{lightgray204}{RGB}{204,204,204}
\definecolor{steelblue31119180}{RGB}{31,119,180}

\begin{axis}[
axis line style={lightgray204},
legend style={
  fill opacity=0.8,
  draw opacity=1,
  text opacity=1,
  at={(0.91,0.5)},
  anchor=east,
  draw=lightgray204
},
scaled y ticks=true,
tick align=outside,
x grid style={lightgray204},
xlabel style={align=center, below=7.7mm},
xlabel={Factor reducing weights \\ to requests in prediction horizon},
xmajorgrids,
xmajorticks=false,
xmajorticks=true,
xmin=-0.5, xmax=9.5,
xtick style={color=darkslategray38},
xtick style={draw=none},
xtick={0,1,2,3,4,5,6,7,8,9},
xticklabel style={rotate=90.0},
xticklabels={0.1,0.2,0.3,0.4,0.5,0.6,0.7,0.8,0.9,1.0},
y grid style={lightgray204},
ylabel style={align=center},
ylabel={Profit [\$]},
ymajorgrids,
ymajorticks=false,
ymajorticks=true,
ymin=22181.1136189223, ymax=43283.6961475403,
ytick style={color=darkslategray38},
ytick style={draw=none}
]
\path [draw=black, fill=steelblue31119180]
(axis cs:-0.25,38706.301588859)
--(axis cs:0.25,38706.301588859)
--(axis cs:0.25,41957.691377905)
--(axis cs:-0.25,41957.691377905)
--(axis cs:-0.25,38706.301588859)
--cycle;
\addplot [black, forget plot]
table {%
0 38706.301588859
0 38706.301588859
};
\addplot [black, forget plot]
table {%
0 41957.691377905
0 42234.5542278657
};
\addplot [black, forget plot]
table {%
-0.125 38706.301588859
0.125 38706.301588859
};
\addplot [black, forget plot]
table {%
-0.125 42234.5542278657
0.125 42234.5542278657
};
\addplot [black, mark=o, mark size=3, mark options={solid,fill opacity=0}, only marks, forget plot]
table {%
0 23219.9683217656
};
\path [draw=black, fill=steelblue31119180]
(axis cs:0.75,38783.5472112402)
--(axis cs:1.25,38783.5472112402)
--(axis cs:1.25,41998.3510159261)
--(axis cs:0.75,41998.3510159261)
--(axis cs:0.75,38783.5472112402)
--cycle;
\addplot [black, forget plot]
table {%
1 38783.5472112402
1 38783.5472112402
};
\addplot [black, forget plot]
table {%
1 41998.3510159261
1 42255.4787446977
};
\addplot [black, forget plot]
table {%
0.875 38783.5472112402
1.125 38783.5472112402
};
\addplot [black, forget plot]
table {%
0.875 42255.4787446977
1.125 42255.4787446977
};
\addplot [black, mark=o, mark size=3, mark options={solid,fill opacity=0}, only marks, forget plot]
table {%
1 23225.0952357815
};
\path [draw=black, fill=steelblue31119180]
(axis cs:1.75,38709.9050636285)
--(axis cs:2.25,38709.9050636285)
--(axis cs:2.25,42060.3126434388)
--(axis cs:1.75,42060.3126434388)
--(axis cs:1.75,38709.9050636285)
--cycle;
\addplot [black, forget plot]
table {%
2 38709.9050636285
2 38709.9050636285
};
\addplot [black, forget plot]
table {%
2 42060.3126434388
2 42257.9864358213
};
\addplot [black, forget plot]
table {%
1.875 38709.9050636285
2.125 38709.9050636285
};
\addplot [black, forget plot]
table {%
1.875 42257.9864358213
2.125 42257.9864358213
};
\addplot [black, mark=o, mark size=3, mark options={solid,fill opacity=0}, only marks, forget plot]
table {%
2 23140.3219156777
};
\path [draw=black, fill=steelblue31119180]
(axis cs:2.75,38747.9324873342)
--(axis cs:3.25,38747.9324873342)
--(axis cs:3.25,41758.5800812158)
--(axis cs:2.75,41758.5800812158)
--(axis cs:2.75,38747.9324873342)
--cycle;
\addplot [black, forget plot]
table {%
3 38747.9324873342
3 38747.9324873342
};
\addplot [black, forget plot]
table {%
3 41758.5800812158
3 42235.3478161762
};
\addplot [black, forget plot]
table {%
2.875 38747.9324873342
3.125 38747.9324873342
};
\addplot [black, forget plot]
table {%
2.875 42235.3478161762
3.125 42235.3478161762
};
\addplot [black, mark=o, mark size=3, mark options={solid,fill opacity=0}, only marks, forget plot]
table {%
3 23189.7621260598
};
\path [draw=black, fill=steelblue31119180]
(axis cs:3.75,38859.2677372875)
--(axis cs:4.25,38859.2677372875)
--(axis cs:4.25,41756.3601601775)
--(axis cs:3.75,41756.3601601775)
--(axis cs:3.75,38859.2677372875)
--cycle;
\addplot [black, forget plot]
table {%
4 38859.2677372875
4 38859.2677372875
};
\addplot [black, forget plot]
table {%
4 41756.3601601775
4 42324.4878507849
};
\addplot [black, forget plot]
table {%
3.875 38859.2677372875
4.125 38859.2677372875
};
\addplot [black, forget plot]
table {%
3.875 42324.4878507849
4.125 42324.4878507849
};
\addplot [black, mark=o, mark size=3, mark options={solid,fill opacity=0}, only marks, forget plot]
table {%
4 23153.1302600446
};
\path [draw=black, fill=steelblue31119180]
(axis cs:4.75,38852.4736187963)
--(axis cs:5.25,38852.4736187963)
--(axis cs:5.25,41815.5433015615)
--(axis cs:4.75,41815.5433015615)
--(axis cs:4.75,38852.4736187963)
--cycle;
\addplot [black, forget plot]
table {%
5 38852.4736187963
5 38852.4736187963
};
\addplot [black, forget plot]
table {%
5 41815.5433015615
5 42251.2696682654
};
\addplot [black, forget plot]
table {%
4.875 38852.4736187963
5.125 38852.4736187963
};
\addplot [black, forget plot]
table {%
4.875 42251.2696682654
5.125 42251.2696682654
};
\addplot [black, mark=o, mark size=3, mark options={solid,fill opacity=0}, only marks, forget plot]
table {%
5 23172.242884841
};
\path [draw=black, fill=steelblue31119180]
(axis cs:5.75,38956.366558096)
--(axis cs:6.25,38956.366558096)
--(axis cs:6.25,41961.2249857196)
--(axis cs:5.75,41961.2249857196)
--(axis cs:5.75,38956.366558096)
--cycle;
\addplot [black, forget plot]
table {%
6 38956.366558096
6 38956.366558096
};
\addplot [black, forget plot]
table {%
6 41961.2249857196
6 42094.5316623377
};
\addplot [black, forget plot]
table {%
5.875 38956.366558096
6.125 38956.366558096
};
\addplot [black, forget plot]
table {%
5.875 42094.5316623377
6.125 42094.5316623377
};
\addplot [black, mark=o, mark size=3, mark options={solid,fill opacity=0}, only marks, forget plot]
table {%
6 23189.5543826725
};
\path [draw=black, fill=steelblue31119180]
(axis cs:6.75,38857.8679632635)
--(axis cs:7.25,38857.8679632635)
--(axis cs:7.25,41844.3748230238)
--(axis cs:6.75,41844.3748230238)
--(axis cs:6.75,38857.8679632635)
--cycle;
\addplot [black, forget plot]
table {%
7 38857.8679632635
7 38857.8679632635
};
\addplot [black, forget plot]
table {%
7 41844.3748230238
7 42122.021272072
};
\addplot [black, forget plot]
table {%
6.875 38857.8679632635
7.125 38857.8679632635
};
\addplot [black, forget plot]
table {%
6.875 42122.021272072
7.125 42122.021272072
};
\addplot [black, mark=o, mark size=3, mark options={solid,fill opacity=0}, only marks, forget plot]
table {%
7 23209.0289297899
7 23209.0289297899
};
\path [draw=black, fill=steelblue31119180]
(axis cs:7.75,38905.5000582277)
--(axis cs:8.25,38905.5000582277)
--(axis cs:8.25,41939.6804836496)
--(axis cs:7.75,41939.6804836496)
--(axis cs:7.75,38905.5000582277)
--cycle;
\addplot [black, forget plot]
table {%
8 38905.5000582277
8 38905.5000582277
};
\addplot [black, forget plot]
table {%
8 41939.6804836496
8 41958.2904332653
};
\addplot [black, forget plot]
table {%
7.875 38905.5000582277
8.125 38905.5000582277
};
\addplot [black, forget plot]
table {%
7.875 41958.2904332653
8.125 41958.2904332653
};
\addplot [black, mark=o, mark size=3, mark options={solid,fill opacity=0}, only marks, forget plot]
table {%
8 23227.6441217697
8 23227.6441217697
};
\path [draw=black, fill=steelblue31119180]
(axis cs:8.75,38807.7939327315)
--(axis cs:9.25,38807.7939327315)
--(axis cs:9.25,41839.3479928796)
--(axis cs:8.75,41839.3479928796)
--(axis cs:8.75,38807.7939327315)
--cycle;
\addplot [black, forget plot]
table {%
9 38807.7939327315
9 38807.7939327315
};
\addplot [black, forget plot]
table {%
9 41839.3479928796
9 42130.3475819492
};
\addplot [black, forget plot]
table {%
8.875 38807.7939327315
9.125 38807.7939327315
};
\addplot [black, forget plot]
table {%
8.875 42130.3475819492
9.125 42130.3475819492
};
\addplot [black, mark=o, mark size=3, mark options={solid,fill opacity=0}, only marks, forget plot]
table {%
9 23168.7967843926
};
\addplot [darkorange25512714, forget plot]
table {%
-0.25 39720.0426511069
0.25 39720.0426511069
};
\addplot [darkorange25512714, forget plot]
table {%
0.75 39651.2157652232
1.25 39651.2157652232
};
\addplot [darkorange25512714, forget plot]
table {%
1.75 39612.3307146386
2.25 39612.3307146386
};
\addplot [darkorange25512714, forget plot]
table {%
2.75 39674.190549741
3.25 39674.190549741
};
\addplot [darkorange25512714, forget plot]
table {%
3.75 39665.8329138655
4.25 39665.8329138655
};
\addplot [darkorange25512714, forget plot]
table {%
4.75 39619.6325999346
5.25 39619.6325999346
};
\addplot [darkorange25512714, forget plot]
table {%
5.75 39644.6224014136
6.25 39644.6224014136
};
\addplot [darkorange25512714, forget plot]
table {%
6.75 39693.4247726213
7.25 39693.4247726213
};
\addplot [darkorange25512714, forget plot]
table {%
7.75 39701.7476806548
8.25 39701.7476806548
};
\addplot [darkorange25512714, forget plot]
table {%
8.75 39762.9873265188
9.25 39762.9873265188
};
\end{axis}

\end{tikzpicture}

%% file: figures/evaluate_hyperparameterSelection_extended_horizon_profit_absolute_policy_SB.tex
\begin{tikzpicture}

\definecolor{darkorange25512714}{RGB}{255,127,14}
\definecolor{darkslategray38}{RGB}{38,38,38}
\definecolor{lightgray204}{RGB}{204,204,204}
\definecolor{steelblue31119180}{RGB}{31,119,180}

\begin{axis}[
axis line style={lightgray204},
legend style={fill opacity=0.8, draw opacity=1, text opacity=1, draw=lightgray204},
scaled y ticks=true,
tick align=outside,
x grid style={lightgray204},
xlabel style={below=11mm},
xlabel={Extended horizon [hour:min:sec]},
xmajorgrids,
xmajorticks=false,
xmajorticks=true,
xmin=-0.5, xmax=5.5,
xtick style={color=darkslategray38},
xtick style={draw=none},
xtick={0,1,2,3,4,5},
xticklabel style={rotate=90.0},
xticklabels={00:05:00,00:10:00,00:15:00,00:20:00,00:25:00,00:30:00},
y grid style={lightgray204},
ylabel style={align=center},
ylabel={Profit [\$]},
ymajorgrids,
ymajorticks=false,
ymajorticks=true,
ymin=20771.7997640182, ymax=45223.7607517915,
ytick style={color=darkslategray38},
ytick style={draw=none}
]
\path [draw=black, fill=steelblue31119180]
(axis cs:-0.25,39954.0404728343)
--(axis cs:0.25,39954.0404728343)
--(axis cs:0.25,44046.2103281767)
--(axis cs:-0.25,44046.2103281767)
--(axis cs:-0.25,39954.0404728343)
--cycle;
\addplot [black, forget plot]
table {%
0 39954.0404728343
0 39954.0404728343
};
\addplot [black, forget plot]
table {%
0 44046.2103281767
0 44112.30797962
};
\addplot [black, forget plot]
table {%
-0.125 39954.0404728343
0.125 39954.0404728343
};
\addplot [black, forget plot]
table {%
-0.125 44112.30797962
0.125 44112.30797962
};
\addplot [black, mark=o, mark size=3, mark options={solid,fill opacity=0}, only marks, forget plot]
table {%
0 24543.233017493
};
\path [draw=black, fill=steelblue31119180]
(axis cs:0.75,38516.353481078)
--(axis cs:1.25,38516.353481078)
--(axis cs:1.25,41937.4352235067)
--(axis cs:0.75,41937.4352235067)
--(axis cs:0.75,38516.353481078)
--cycle;
\addplot [black, forget plot]
table {%
1 38516.353481078
1 38516.353481078
};
\addplot [black, forget plot]
table {%
1 41937.4352235067
1 41948.2609816883
};
\addplot [black, forget plot]
table {%
0.875 38516.353481078
1.125 38516.353481078
};
\addplot [black, forget plot]
table {%
0.875 41948.2609816883
1.125 41948.2609816883
};
\addplot [black, mark=o, mark size=3, mark options={solid,fill opacity=0}, only marks, forget plot]
table {%
1 23345.6271565581
};
\path [draw=black, fill=steelblue31119180]
(axis cs:1.75,36711.5423972301)
--(axis cs:2.25,36711.5423972301)
--(axis cs:2.25,39762.8710503252)
--(axis cs:1.75,39762.8710503252)
--(axis cs:1.75,36711.5423972301)
--cycle;
\addplot [black, forget plot]
table {%
2 36711.5423972301
2 36711.5423972301
};
\addplot [black, forget plot]
table {%
2 39762.8710503252
2 39782.3201339274
};
\addplot [black, forget plot]
table {%
1.875 36711.5423972301
2.125 36711.5423972301
};
\addplot [black, forget plot]
table {%
1.875 39782.3201339274
2.125 39782.3201339274
};
\addplot [black, mark=o, mark size=3, mark options={solid,fill opacity=0}, only marks, forget plot]
table {%
2 21883.2525361897
};
\path [draw=black, fill=steelblue31119180]
(axis cs:2.75,36154.2589766238)
--(axis cs:3.25,36154.2589766238)
--(axis cs:3.25,38743.0879261908)
--(axis cs:2.75,38743.0879261908)
--(axis cs:2.75,36154.2589766238)
--cycle;
\addplot [black, forget plot]
table {%
3 36154.2589766238
3 36154.2589766238
};
\addplot [black, forget plot]
table {%
3 38743.0879261908
3 39572.0517133768
};
\addplot [black, forget plot]
table {%
2.875 36154.2589766238
3.125 36154.2589766238
};
\addplot [black, forget plot]
table {%
2.875 39572.0517133768
3.125 39572.0517133768
};
\addplot [black, mark=o, mark size=3, mark options={solid,fill opacity=0}, only marks, forget plot]
table {%
3 22254.6890761458
};
\path [draw=black, fill=steelblue31119180]
(axis cs:3.75,36316.2126606656)
--(axis cs:4.25,36316.2126606656)
--(axis cs:4.25,39185.3813978673)
--(axis cs:3.75,39185.3813978673)
--(axis cs:3.75,36316.2126606656)
--cycle;
\addplot [black, forget plot]
table {%
4 36316.2126606656
4 36316.2126606656
};
\addplot [black, forget plot]
table {%
4 39185.3813978673
4 39596.0118442899
};
\addplot [black, forget plot]
table {%
3.875 36316.2126606656
4.125 36316.2126606656
};
\addplot [black, forget plot]
table {%
3.875 39596.0118442899
4.125 39596.0118442899
};
\addplot [black, mark=o, mark size=3, mark options={solid,fill opacity=0}, only marks, forget plot]
table {%
4 22096.558955722
};
\path [draw=black, fill=steelblue31119180]
(axis cs:4.75,36503.6406423196)
--(axis cs:5.25,36503.6406423196)
--(axis cs:5.25,39505.5181677966)
--(axis cs:4.75,39505.5181677966)
--(axis cs:4.75,36503.6406423196)
--cycle;
\addplot [black, forget plot]
table {%
5 36503.6406423196
5 36503.6406423196
};
\addplot [black, forget plot]
table {%
5 39505.5181677966
5 39577.4505269692
};
\addplot [black, forget plot]
table {%
4.875 36503.6406423196
5.125 36503.6406423196
};
\addplot [black, forget plot]
table {%
4.875 39577.4505269692
5.125 39577.4505269692
};
\addplot [black, mark=o, mark size=3, mark options={solid,fill opacity=0}, only marks, forget plot]
table {%
5 22121.1505170432
};
\addplot [darkorange25512714, forget plot]
table {%
-0.25 41674.2393668681
0.25 41674.2393668681
};
\addplot [darkorange25512714, forget plot]
table {%
0.75 39702.810233994
1.25 39702.810233994
};
\addplot [darkorange25512714, forget plot]
table {%
1.75 37818.4572199072
2.25 37818.4572199072
};
\addplot [darkorange25512714, forget plot]
table {%
2.75 37897.7112711995
3.25 37897.7112711995
};
\addplot [darkorange25512714, forget plot]
table {%
3.75 37676.5923149048
4.25 37676.5923149048
};
\addplot [darkorange25512714, forget plot]
table {%
4.75 37872.0705204379
5.25 37872.0705204379
};
\end{axis}

\end{tikzpicture}

%% file: figures/evaluate_hyperparameterSelection_extended_horizon_profit_absolute_policy_CB.tex
\begin{tikzpicture}

\definecolor{darkorange25512714}{RGB}{255,127,14}
\definecolor{darkslategray38}{RGB}{38,38,38}
\definecolor{lightgray204}{RGB}{204,204,204}
\definecolor{steelblue31119180}{RGB}{31,119,180}

\begin{axis}[
axis line style={lightgray204},
legend style={fill opacity=0.8, draw opacity=1, text opacity=1, draw=lightgray204},
scaled y ticks=true,
tick align=outside,
x grid style={lightgray204},
xlabel style={below=11mm},
xlabel={Extended horizon [hour:min:sec]},
xmajorgrids,
xmajorticks=false,
xmajorticks=true,
xmin=-0.5, xmax=5.5,
xtick style={color=darkslategray38},
xtick style={draw=none},
xtick={0,1,2,3,4,5},
xticklabel style={rotate=90.0},
xticklabels={00:05:00,00:10:00,00:15:00,00:20:00,00:25:00,00:30:00},
y grid style={lightgray204},
ylabel style={align=center},
ylabel={Profit [\$]},
ymajorgrids,
ymajorticks=false,
ymajorticks=true,
ymin=23869.7674486626, ymax=45107.4121536059,
ytick style={color=darkslategray38},
ytick style={draw=none}
]
\path [draw=black, fill=steelblue31119180]
(axis cs:-0.25,39765.1908491305)
--(axis cs:0.25,39765.1908491305)
--(axis cs:0.25,43568.6001336639)
--(axis cs:-0.25,43568.6001336639)
--(axis cs:-0.25,39765.1908491305)
--cycle;
\addplot [black, forget plot]
table {%
0 39765.1908491305
0 39765.1908491305
};
\addplot [black, forget plot]
table {%
0 43568.6001336639
0 44142.0646670176
};
\addplot [black, forget plot]
table {%
-0.125 39765.1908491305
0.125 39765.1908491305
};
\addplot [black, forget plot]
table {%
-0.125 44142.0646670176
0.125 44142.0646670176
};
\addplot [black, mark=o, mark size=3, mark options={solid,fill opacity=0}, only marks, forget plot]
table {%
0 25183.8060637194
};
\path [draw=black, fill=steelblue31119180]
(axis cs:0.75,39134.898622891)
--(axis cs:1.25,39134.898622891)
--(axis cs:1.25,43152.0236107122)
--(axis cs:0.75,43152.0236107122)
--(axis cs:0.75,39134.898622891)
--cycle;
\addplot [black, forget plot]
table {%
1 39134.898622891
1 39134.898622891
};
\addplot [black, forget plot]
table {%
1 43152.0236107122
1 43293.9451783183
};
\addplot [black, forget plot]
table {%
0.875 39134.898622891
1.125 39134.898622891
};
\addplot [black, forget plot]
table {%
0.875 43293.9451783183
1.125 43293.9451783183
};
\addplot [black, mark=o, mark size=3, mark options={solid,fill opacity=0}, only marks, forget plot]
table {%
1 24838.6746800561
};
\path [draw=black, fill=steelblue31119180]
(axis cs:1.75,38261.182876198)
--(axis cs:2.25,38261.182876198)
--(axis cs:2.25,42389.0567243709)
--(axis cs:1.75,42389.0567243709)
--(axis cs:1.75,38261.182876198)
--cycle;
\addplot [black, forget plot]
table {%
2 38261.182876198
2 38261.182876198
};
\addplot [black, forget plot]
table {%
2 42389.0567243709
2 42621.8648831788
};
\addplot [black, forget plot]
table {%
1.875 38261.182876198
2.125 38261.182876198
};
\addplot [black, forget plot]
table {%
1.875 42621.8648831788
2.125 42621.8648831788
};
\addplot [black, mark=o, mark size=3, mark options={solid,fill opacity=0}, only marks, forget plot]
table {%
2 24935.4641566316
};
\path [draw=black, fill=steelblue31119180]
(axis cs:2.75,38402.6199521687)
--(axis cs:3.25,38402.6199521687)
--(axis cs:3.25,42423.8237290912)
--(axis cs:2.75,42423.8237290912)
--(axis cs:2.75,38402.6199521687)
--cycle;
\addplot [black, forget plot]
table {%
3 38402.6199521687
3 38402.6199521687
};
\addplot [black, forget plot]
table {%
3 42423.8237290912
3 42647.9557420067
};
\addplot [black, forget plot]
table {%
2.875 38402.6199521687
3.125 38402.6199521687
};
\addplot [black, forget plot]
table {%
2.875 42647.9557420067
3.125 42647.9557420067
};
\addplot [black, mark=o, mark size=3, mark options={solid,fill opacity=0}, only marks, forget plot]
table {%
3 24835.1149352509
};
\path [draw=black, fill=steelblue31119180]
(axis cs:3.75,38320.9134582053)
--(axis cs:4.25,38320.9134582053)
--(axis cs:4.25,42380.0981848445)
--(axis cs:3.75,42380.0981848445)
--(axis cs:3.75,38320.9134582053)
--cycle;
\addplot [black, forget plot]
table {%
4 38320.9134582053
4 38320.9134582053
};
\addplot [black, forget plot]
table {%
4 42380.0981848445
4 42620.9482140373
};
\addplot [black, forget plot]
table {%
3.875 38320.9134582053
4.125 38320.9134582053
};
\addplot [black, forget plot]
table {%
3.875 42620.9482140373
4.125 42620.9482140373
};
\addplot [black, mark=o, mark size=3, mark options={solid,fill opacity=0}, only marks, forget plot]
table {%
4 24923.0160957724
};
\path [draw=black, fill=steelblue31119180]
(axis cs:4.75,38229.4684007854)
--(axis cs:5.25,38229.4684007854)
--(axis cs:5.25,42456.0917659359)
--(axis cs:4.75,42456.0917659359)
--(axis cs:4.75,38229.4684007854)
--cycle;
\addplot [black, forget plot]
table {%
5 38229.4684007854
5 38229.4684007854
};
\addplot [black, forget plot]
table {%
5 42456.0917659359
5 42591.5722259594
};
\addplot [black, forget plot]
table {%
4.875 38229.4684007854
5.125 38229.4684007854
};
\addplot [black, forget plot]
table {%
4.875 42591.5722259594
5.125 42591.5722259594
};
\addplot [black, mark=o, mark size=3, mark options={solid,fill opacity=0}, only marks, forget plot]
table {%
5 24875.2019506649
};
\addplot [darkorange25512714, forget plot]
table {%
-0.25 41121.5555192677
0.25 41121.5555192677
};
\addplot [darkorange25512714, forget plot]
table {%
0.75 40842.1140348613
1.25 40842.1140348613
};
\addplot [darkorange25512714, forget plot]
table {%
1.75 40390.8015610913
2.25 40390.8015610913
};
\addplot [darkorange25512714, forget plot]
table {%
2.75 40313.8721561623
3.25 40313.8721561623
};
\addplot [darkorange25512714, forget plot]
table {%
3.75 40344.3461351448
4.25 40344.3461351448
};
\addplot [darkorange25512714, forget plot]
table {%
4.75 40463.1580098545
5.25 40463.1580098545
};
\end{axis}

\end{tikzpicture}

%% file: figures/evaluate_fix_capacity_profit_absolute.tex
\begin{tikzpicture}

\definecolor{darkorange25512714}{RGB}{255,127,14}
\definecolor{darkslategray38}{RGB}{38,38,38}
\definecolor{lightgray204}{RGB}{204,204,204}
\definecolor{steelblue31119180}{RGB}{31,119,180}

\begin{axis}[
axis line style={lightgray204},
legend style={
  fill opacity=0.8,
  draw opacity=1,
  text opacity=1,
  at={(0.91,0.5)},
  anchor=east,
  draw=lightgray204
},
scaled y ticks=true,
tick align=outside,
x grid style={lightgray204},
xlabel style={align=center, below=5.8mm},
xlabel={Num. of capacity vertices \\ in rebalancing cells},
xmajorgrids,
xmajorticks=false,
xmajorticks=true,
xmin=-0.5, xmax=9.5,
xtick style={color=darkslategray38},
xtick style={draw=none},
xtick={0,1,2,3,4,5,6,7,8,9},
xticklabel style={rotate=90.0},
xticklabels={1.0,2.0,3.0,4.0,5.0,6.0,7.0,8.0,9.0,10.0},
y grid style={lightgray204},
ylabel style={align=center},
ylabel={Profit [\$]},
ymajorgrids,
ymajorticks=false,
ymajorticks=true,
ymin=23678.5101278782, ymax=44530.1301433265,
ytick style={color=darkslategray38},
ytick style={draw=none}
]
\path [draw=black, fill=steelblue31119180]
(axis cs:-0.25,39489.825292102)
--(axis cs:0.25,39489.825292102)
--(axis cs:0.25,43229.4472659629)
--(axis cs:-0.25,43229.4472659629)
--(axis cs:-0.25,39489.825292102)
--cycle;
\addplot [black, forget plot]
table {%
0 39489.825292102
0 39489.825292102
};
\addplot [black, forget plot]
table {%
0 43229.4472659629
0 43582.3292335334
};
\addplot [black, forget plot]
table {%
-0.125 39489.825292102
0.125 39489.825292102
};
\addplot [black, forget plot]
table {%
-0.125 43582.3292335334
0.125 43582.3292335334
};
\addplot [black, mark=o, mark size=3, mark options={solid,fill opacity=0}, only marks, forget plot]
table {%
0 25151.9176293159
};
\path [draw=black, fill=steelblue31119180]
(axis cs:0.75,39298.6463248644)
--(axis cs:1.25,39298.6463248644)
--(axis cs:1.25,43198.2489937904)
--(axis cs:0.75,43198.2489937904)
--(axis cs:0.75,39298.6463248644)
--cycle;
\addplot [black, forget plot]
table {%
1 39298.6463248644
1 39298.6463248644
};
\addplot [black, forget plot]
table {%
1 43198.2489937904
1 43236.9208225487
};
\addplot [black, forget plot]
table {%
0.875 39298.6463248644
1.125 39298.6463248644
};
\addplot [black, forget plot]
table {%
0.875 43236.9208225487
1.125 43236.9208225487
};
\addplot [black, mark=o, mark size=3, mark options={solid,fill opacity=0}, only marks, forget plot]
table {%
1 24812.1361910198
};
\path [draw=black, fill=steelblue31119180]
(axis cs:1.75,39426.074826634)
--(axis cs:2.25,39426.074826634)
--(axis cs:2.25,43108.3976181312)
--(axis cs:1.75,43108.3976181312)
--(axis cs:1.75,39426.074826634)
--cycle;
\addplot [black, forget plot]
table {%
2 39426.074826634
2 39426.074826634
};
\addplot [black, forget plot]
table {%
2 43108.3976181312
2 43168.0818305295
};
\addplot [black, forget plot]
table {%
1.875 39426.074826634
2.125 39426.074826634
};
\addplot [black, forget plot]
table {%
1.875 43168.0818305295
2.125 43168.0818305295
};
\addplot [black, mark=o, mark size=3, mark options={solid,fill opacity=0}, only marks, forget plot]
table {%
2 24723.7303652211
};
\path [draw=black, fill=steelblue31119180]
(axis cs:2.75,39203.926315721)
--(axis cs:3.25,39203.926315721)
--(axis cs:3.25,43139.2984191791)
--(axis cs:2.75,43139.2984191791)
--(axis cs:2.75,39203.926315721)
--cycle;
\addplot [black, forget plot]
table {%
3 39203.926315721
3 39203.926315721
};
\addplot [black, forget plot]
table {%
3 43139.2984191791
3 43161.4765642622
};
\addplot [black, forget plot]
table {%
2.875 39203.926315721
3.125 39203.926315721
};
\addplot [black, forget plot]
table {%
2.875 43161.4765642622
3.125 43161.4765642622
};
\addplot [black, mark=o, mark size=3, mark options={solid,fill opacity=0}, only marks, forget plot]
table {%
3 24675.6327486772
};
\path [draw=black, fill=steelblue31119180]
(axis cs:3.75,39145.383365635)
--(axis cs:4.25,39145.383365635)
--(axis cs:4.25,43026.4180086463)
--(axis cs:3.75,43026.4180086463)
--(axis cs:3.75,39145.383365635)
--cycle;
\addplot [black, forget plot]
table {%
4 39145.383365635
4 39145.383365635
};
\addplot [black, forget plot]
table {%
4 43026.4180086463
4 43154.8088089188
};
\addplot [black, forget plot]
table {%
3.875 39145.383365635
4.125 39145.383365635
};
\addplot [black, forget plot]
table {%
3.875 43154.8088089188
4.125 43154.8088089188
};
\addplot [black, mark=o, mark size=3, mark options={solid,fill opacity=0}, only marks, forget plot]
table {%
4 24662.3711425711
};
\path [draw=black, fill=steelblue31119180]
(axis cs:4.75,39120.0456579621)
--(axis cs:5.25,39120.0456579621)
--(axis cs:5.25,42934.0769606003)
--(axis cs:4.75,42934.0769606003)
--(axis cs:4.75,39120.0456579621)
--cycle;
\addplot [black, forget plot]
table {%
5 39120.0456579621
5 39120.0456579621
};
\addplot [black, forget plot]
table {%
5 42934.0769606003
5 42945.6172824248
};
\addplot [black, forget plot]
table {%
4.875 39120.0456579621
5.125 39120.0456579621
};
\addplot [black, forget plot]
table {%
4.875 42945.6172824248
5.125 42945.6172824248
};
\addplot [black, mark=o, mark size=3, mark options={solid,fill opacity=0}, only marks, forget plot]
table {%
5 24692.8439116451
};
\path [draw=black, fill=steelblue31119180]
(axis cs:5.75,39047.5493102499)
--(axis cs:6.25,39047.5493102499)
--(axis cs:6.25,42995.1920577414)
--(axis cs:5.75,42995.1920577414)
--(axis cs:5.75,39047.5493102499)
--cycle;
\addplot [black, forget plot]
table {%
6 39047.5493102499
6 39047.5493102499
};
\addplot [black, forget plot]
table {%
6 42995.1920577414
6 43056.6774635515
};
\addplot [black, forget plot]
table {%
5.875 39047.5493102499
6.125 39047.5493102499
};
\addplot [black, forget plot]
table {%
5.875 43056.6774635515
6.125 43056.6774635515
};
\addplot [black, mark=o, mark size=3, mark options={solid,fill opacity=0}, only marks, forget plot]
table {%
6 24626.3110376713
};
\path [draw=black, fill=steelblue31119180]
(axis cs:6.75,39173.0580066535)
--(axis cs:7.25,39173.0580066535)
--(axis cs:7.25,42933.7592350992)
--(axis cs:6.75,42933.7592350992)
--(axis cs:6.75,39173.0580066535)
--cycle;
\addplot [black, forget plot]
table {%
7 39173.0580066535
7 39173.0580066535
};
\addplot [black, forget plot]
table {%
7 42933.7592350992
7 42967.744589816
};
\addplot [black, forget plot]
table {%
6.875 39173.0580066535
7.125 39173.0580066535
};
\addplot [black, forget plot]
table {%
6.875 42967.744589816
7.125 42967.744589816
};
\addplot [black, mark=o, mark size=3, mark options={solid,fill opacity=0}, only marks, forget plot]
table {%
7 24671.4733870043
};
\path [draw=black, fill=steelblue31119180]
(axis cs:7.75,39089.9040822837)
--(axis cs:8.25,39089.9040822837)
--(axis cs:8.25,42950.5272761989)
--(axis cs:7.75,42950.5272761989)
--(axis cs:7.75,39089.9040822837)
--cycle;
\addplot [black, forget plot]
table {%
8 39089.9040822837
8 39089.9040822837
};
\addplot [black, forget plot]
table {%
8 42950.5272761989
8 43018.3511729734
};
\addplot [black, forget plot]
table {%
7.875 39089.9040822837
8.125 39089.9040822837
};
\addplot [black, forget plot]
table {%
7.875 43018.3511729734
8.125 43018.3511729734
};
\addplot [black, mark=o, mark size=3, mark options={solid,fill opacity=0}, only marks, forget plot]
table {%
8 24671.4661751173
};
\path [draw=black, fill=steelblue31119180]
(axis cs:8.75,39155.4087479198)
--(axis cs:9.25,39155.4087479198)
--(axis cs:9.25,42933.2112978891)
--(axis cs:8.75,42933.2112978891)
--(axis cs:8.75,39155.4087479198)
--cycle;
\addplot [black, forget plot]
table {%
9 39155.4087479198
9 39155.4087479198
};
\addplot [black, forget plot]
table {%
9 42933.2112978891
9 42953.6618880534
};
\addplot [black, forget plot]
table {%
8.875 39155.4087479198
9.125 39155.4087479198
};
\addplot [black, forget plot]
table {%
8.875 42953.6618880534
9.125 42953.6618880534
};
\addplot [black, mark=o, mark size=3, mark options={solid,fill opacity=0}, only marks, forget plot]
table {%
9 24655.6192893025
};
\addplot [darkorange25512714, forget plot]
table {%
-0.25 40962.4424348062
0.25 40962.4424348062
};
\addplot [darkorange25512714, forget plot]
table {%
0.75 41154.0416698584
1.25 41154.0416698584
};
\addplot [darkorange25512714, forget plot]
table {%
1.75 40874.5032206692
2.25 40874.5032206692
};
\addplot [darkorange25512714, forget plot]
table {%
2.75 40980.4548581995
3.25 40980.4548581995
};
\addplot [darkorange25512714, forget plot]
table {%
3.75 40848.7892998484
4.25 40848.7892998484
};
\addplot [darkorange25512714, forget plot]
table {%
4.75 40893.5231243157
5.25 40893.5231243157
};
\addplot [darkorange25512714, forget plot]
table {%
5.75 40900.5067830352
6.25 40900.5067830352
};
\addplot [darkorange25512714, forget plot]
table {%
6.75 40907.8905474628
7.25 40907.8905474628
};
\addplot [darkorange25512714, forget plot]
table {%
7.75 40876.3195437724
8.25 40876.3195437724
};
\addplot [darkorange25512714, forget plot]
table {%
8.75 40907.3579119684
9.25 40907.3579119684
};
\end{axis}

\end{tikzpicture}

%% file: figures/comparisonLossLinearNNpolicy_SB.tex
\begin{tikzpicture}

\definecolor{darkgray176}{RGB}{176,176,176}
\definecolor{lightgray204}{RGB}{204,204,204}

\begin{axis}[
legend cell align={left},
legend style={fill opacity=0.8, draw opacity=1, text opacity=1, draw=lightgray204},
tick align=outside,
tick pos=left,
x grid style={darkgray176},
xlabel={Fleet size},
xmin=275, xmax=5225,
xtick style={color=black},
y grid style={darkgray176},
ylabel={Loss value},
ymin=-19493.0662718108, ymax=412665.730763065,
ytick style={color=black}
]
\addplot [semithick, blue]
table {%
500 310274.513522515
1000 393022.149079662
2000 27748.3670513977
3000 27023.1862374981
4000 27326.2627585553
5000 26804.7837020397
};
\addlegendentry{Loss~Linear~predictor}
\addplot [semithick, red]
table {%
500 150.515411592691
1000 250.596321228213
2000 687.3280245445
3000 202650.449760221
4000 1081.24512125929
5000 850.720202001004
};
\addlegendentry{Loss~FNN~predictor}
\end{axis}

\end{tikzpicture}

%% file: figures/comparisonLossLinearNNpolicy_CB.tex
\begin{tikzpicture}

\definecolor{darkgray176}{RGB}{176,176,176}
\definecolor{lightgray204}{RGB}{204,204,204}

\begin{axis}[
legend cell align={left},
legend style={
  fill opacity=0.8,
  draw opacity=1,
  text opacity=1,
  at={(0.91,0.5)},
  anchor=east,
  draw=lightgray204
},
tick align=outside,
tick pos=left,
x grid style={darkgray176},
xlabel={Fleet size},
xmin=275, xmax=5225,
xtick style={color=black},
y grid style={darkgray176},
ylabel={Loss value},
ymin=-9827.29498329475, ymax=233922.20808877,
ytick style={color=black}
]
\addplot [semithick, blue]
table {%
500 101214.590085728
1000 163828.165592908
2000 220657.631486652
3000 219636.730853338
4000 222842.685221858
5000 218245.697489549
};
\addlegendentry{Loss~Linear~predictor}
\addplot [semithick, red]
table {%
500 7318.9943910905
1000 1252.2278836173
2000 21329.2682745505
3000 11319.9238196824
4000 2726.03830349432
5000 1529.14288897699
};
\addlegendentry{Loss~FNN~predictor}
\end{axis}

\end{tikzpicture}

%% file: figures/comparisonEqualityLinearNNpolicy_SB.tex
\begin{tikzpicture}

\definecolor{darkgray176}{RGB}{176,176,176}
\definecolor{darkorange25512714}{RGB}{255,127,14}
\definecolor{lightgray204}{RGB}{204,204,204}

\begin{axis}[
legend style={fill opacity=0.8, draw opacity=1, text opacity=1, draw=lightgray204},
tick align=outside,
tick pos=left,
x grid style={darkgray176},
xlabel={Fleet size},
xmin=0.5, xmax=6.5,
xtick style={color=black},
xtick={1,2,3,4,5,6},
xticklabels={500,1000,2000,3000,4000,5000},
y grid style={darkgray176},
ylabel={\(\displaystyle f = (g_{FNN} - g_{Linear}) / g_{Linear}\)},
ymin=-0.00966089469268139, ymax=0.0145593606559635,
ytick style={color=black}
]
\addplot [black, forget plot]
table {%
0.75 -0.00142700477298154
1.25 -0.00142700477298154
1.25 0.00251317115649964
0.75 0.00251317115649964
0.75 -0.00142700477298154
};
\addplot [black, forget plot]
table {%
1 -0.00142700477298154
1 -0.00699912510936133
};
\addplot [black, forget plot]
table {%
1 0.00251317115649964
1 0.00607034458132476
};
\addplot [black, forget plot]
table {%
0.875 -0.00699912510936133
1.125 -0.00699912510936133
};
\addplot [black, forget plot]
table {%
0.875 0.00607034458132476
1.125 0.00607034458132476
};
\addplot [black, mark=o, mark size=3, mark options={solid,fill opacity=0}, only marks, forget plot]
table {%
1 0.00885167464114833
};
\addplot [black, forget plot]
table {%
1.75 -0.000280816500399304
2.25 -0.000280816500399304
2.25 0.00386168392653228
1.75 0.00386168392653228
1.75 -0.000280816500399304
};
\addplot [black, forget plot]
table {%
2 -0.000280816500399304
2 -0.00439334455038325
};
\addplot [black, forget plot]
table {%
2 0.00386168392653228
2 0.00943964651536453
};
\addplot [black, forget plot]
table {%
1.875 -0.00439334455038325
2.125 -0.00439334455038325
};
\addplot [black, forget plot]
table {%
1.875 0.00943964651536453
2.125 0.00943964651536453
};
\addplot [black, mark=o, mark size=3, mark options={solid,fill opacity=0}, only marks, forget plot]
table {%
2 -0.00855997399501571
2 0.0130620737486996
2 0.0134584399582978
};
\addplot [black, forget plot]
table {%
2.75 3.48569702321474e-05
3.25 3.48569702321474e-05
3.25 0.000999309462610134
2.75 0.000999309462610134
2.75 3.48569702321474e-05
};
\addplot [black, forget plot]
table {%
3 3.48569702321474e-05
3 -0.00119013271783035
};
\addplot [black, forget plot]
table {%
3 0.000999309462610134
3 0.00176712197917662
};
\addplot [black, forget plot]
table {%
2.875 -0.00119013271783035
3.125 -0.00119013271783035
};
\addplot [black, forget plot]
table {%
2.875 0.00176712197917662
3.125 0.00176712197917662
};
\addplot [black, mark=o, mark size=3, mark options={solid,fill opacity=0}, only marks, forget plot]
table {%
3 -0.00163126919140225
3 0.00252089690858432
3 0.00250072695551032
};
\addplot [black, forget plot]
table {%
3.75 -0.00242026818816975
4.25 -0.00242026818816975
4.25 -0.00137900546528754
3.75 -0.00137900546528754
3.75 -0.00242026818816975
};
\addplot [black, forget plot]
table {%
4 -0.00242026818816975
4 -0.00395341382626349
};
\addplot [black, forget plot]
table {%
4 -0.00137900546528754
4 -0.000116577290743763
};
\addplot [black, forget plot]
table {%
3.875 -0.00395341382626349
4.125 -0.00395341382626349
};
\addplot [black, forget plot]
table {%
3.875 -0.000116577290743763
4.125 -0.000116577290743763
};
\addplot [black, mark=o, mark size=3, mark options={solid,fill opacity=0}, only marks, forget plot]
table {%
4 -0.00637095804784229
};
\addplot [black, forget plot]
table {%
4.75 0.00319617189559113
5.25 0.00319617189559113
5.25 0.00424226490909636
4.75 0.00424226490909636
4.75 0.00319617189559113
};
\addplot [black, forget plot]
table {%
5 0.00319617189559113
5 0.00185905019435525
};
\addplot [black, forget plot]
table {%
5 0.00424226490909636
5 0.00570342205323194
};
\addplot [black, forget plot]
table {%
4.875 0.00185905019435525
5.125 0.00185905019435525
};
\addplot [black, forget plot]
table {%
4.875 0.00570342205323194
5.125 0.00570342205323194
};
\addplot [black, mark=o, mark size=3, mark options={solid,fill opacity=0}, only marks, forget plot]
table {%
5 0.00100032153192097
5 0.000400106695118698
5 0.00720036285293117
};
\addplot [black, forget plot]
table {%
5.75 -0.00060338056976615
6.25 -0.00060338056976615
6.25 5.93384332501099e-05
5.75 5.93384332501099e-05
5.75 -0.00060338056976615
};
\addplot [black, forget plot]
table {%
6 -0.00060338056976615
6 -0.00137102565776588
};
\addplot [black, forget plot]
table {%
6 5.93384332501099e-05
6 0.000956988675634005
};
\addplot [black, forget plot]
table {%
5.875 -0.00137102565776588
6.125 -0.00137102565776588
};
\addplot [black, forget plot]
table {%
5.875 0.000956988675634005
6.125 0.000956988675634005
};
\addplot [black, mark=o, mark size=3, mark options={solid,fill opacity=0}, only marks, forget plot]
table {%
6 -0.00198548116895204
};
\addplot [darkorange25512714, forget plot]
table {%
0.75 0.000442119414435371
1.25 0.000442119414435371
};
\addplot [darkorange25512714, forget plot]
table {%
1.75 0.00159974636725261
2.25 0.00159974636725261
};
\addplot [darkorange25512714, forget plot]
table {%
2.75 0.000666611855987409
3.25 0.000666611855987409
};
\addplot [darkorange25512714, forget plot]
table {%
3.75 -0.00178148399759784
4.25 -0.00178148399759784
};
\addplot [darkorange25512714, forget plot]
table {%
4.75 0.00382409774723375
5.25 0.00382409774723375
};
\addplot [darkorange25512714, forget plot]
table {%
5.75 -0.000173870566808934
6.25 -0.000173870566808934
};
\end{axis}

\end{tikzpicture}

%% file: figures/comparisonEqualityLinearNNpolicy_CB.tex
\begin{tikzpicture}

\definecolor{darkgray176}{RGB}{176,176,176}
\definecolor{darkorange25512714}{RGB}{255,127,14}
\definecolor{lightgray204}{RGB}{204,204,204}

\begin{axis}[
legend style={fill opacity=0.8, draw opacity=1, text opacity=1, draw=lightgray204},
tick align=outside,
tick pos=left,
x grid style={darkgray176},
xlabel={Fleet size},
xmin=0.5, xmax=6.5,
xtick style={color=black},
xtick={1,2,3,4,5,6},
xticklabels={500,1000,2000,3000,4000,5000},
y grid style={darkgray176},
ylabel={\(\displaystyle f = (g_{FNN} - g_{Linear}) / g_{Linear}\)},
ymin=-0.0104383566080676, ymax=0.0151582712314779,
ytick style={color=black}
]
\addplot [black, forget plot]
table {%
0.75 0.00571986536566289
1.25 0.00571986536566289
1.25 0.00909092303129483
0.75 0.00909092303129483
0.75 0.00571986536566289
};
\addplot [black, forget plot]
table {%
1 0.00571986536566289
1 0.00200059276822762
};
\addplot [black, forget plot]
table {%
1 0.00909092303129483
1 0.0139947881478622
};
\addplot [black, forget plot]
table {%
0.875 0.00200059276822762
1.125 0.00200059276822762
};
\addplot [black, forget plot]
table {%
0.875 0.0139947881478622
1.125 0.0139947881478622
};
\addplot [black, mark=o, mark size=3, mark options={solid,fill opacity=0}, only marks, forget plot]
table {%
1 0.000535004585753592
};
\addplot [black, forget plot]
table {%
1.75 0.00635336765621494
2.25 0.00635336765621494
2.25 0.00905433299593861
1.75 0.00905433299593861
1.75 0.00635336765621494
};
\addplot [black, forget plot]
table {%
2 0.00635336765621494
2 0.00234591597355513
};
\addplot [black, forget plot]
table {%
2 0.00905433299593861
2 0.011308259402724
};
\addplot [black, forget plot]
table {%
1.875 0.00234591597355513
2.125 0.00234591597355513
};
\addplot [black, forget plot]
table {%
1.875 0.011308259402724
2.125 0.011308259402724
};
\addplot [black, mark=o, mark size=3, mark options={solid,fill opacity=0}, only marks, forget plot]
table {%
2 0.0132674134801928
};
\addplot [black, forget plot]
table {%
2.75 0.00223448720112305
3.25 0.00223448720112305
3.25 0.00354866097046186
2.75 0.00354866097046186
2.75 0.00223448720112305
};
\addplot [black, forget plot]
table {%
3 0.00223448720112305
3 0.000291766353504114
};
\addplot [black, forget plot]
table {%
3 0.00354866097046186
3 0.00511993446483885
};
\addplot [black, forget plot]
table {%
2.875 0.000291766353504114
3.125 0.000291766353504114
};
\addplot [black, forget plot]
table {%
2.875 0.00511993446483885
3.125 0.00511993446483885
};
\addplot [black, mark=o, mark size=3, mark options={solid,fill opacity=0}, only marks, forget plot]
table {%
3 0.00723120837297812
3 0.00567906364706697
};
\addplot [black, forget plot]
table {%
3.75 -0.0063904379503072
4.25 -0.0063904379503072
4.25 -0.00483417128279404
3.75 -0.00483417128279404
3.75 -0.0063904379503072
};
\addplot [black, forget plot]
table {%
4 -0.0063904379503072
4 -0.00868095821874689
};
\addplot [black, forget plot]
table {%
4 -0.00483417128279404
4 -0.00360668624141678
};
\addplot [black, forget plot]
table {%
3.875 -0.00868095821874689
4.125 -0.00868095821874689
};
\addplot [black, forget plot]
table {%
3.875 -0.00360668624141678
4.125 -0.00360668624141678
};
\addplot [black, mark=o, mark size=3, mark options={solid,fill opacity=0}, only marks, forget plot]
table {%
4 -0.00927487352445194
};
\addplot [black, forget plot]
table {%
4.75 0.000800292520600977
5.25 0.000800292520600977
5.25 0.00217765458502174
4.75 0.00217765458502174
4.75 0.000800292520600977
};
\addplot [black, forget plot]
table {%
5 0.000800292520600977
5 -0.0012085255987695
};
\addplot [black, forget plot]
table {%
5 0.00217765458502174
5 0.00338815189345565
};
\addplot [black, forget plot]
table {%
4.875 -0.0012085255987695
5.125 -0.0012085255987695
};
\addplot [black, forget plot]
table {%
4.875 0.00338815189345565
5.125 0.00338815189345565
};
\addplot [black, mark=o, mark size=3, mark options={solid,fill opacity=0}, only marks, forget plot]
table {%
5 -0.0016001097218095
5 0.00564759036144578
};
\addplot [black, forget plot]
table {%
5.75 0.000761562324389274
6.25 0.000761562324389274
6.25 0.00190175804324666
5.75 0.00190175804324666
5.75 0.000761562324389274
};
\addplot [black, forget plot]
table {%
6 0.000761562324389274
6 -0.000859008131943649
};
\addplot [black, forget plot]
table {%
6 0.00190175804324666
6 0.00313729524005452
};
\addplot [black, forget plot]
table {%
5.875 -0.000859008131943649
6.125 -0.000859008131943649
};
\addplot [black, forget plot]
table {%
5.875 0.00313729524005452
6.125 0.00313729524005452
};
\addplot [black, mark=o, mark size=3, mark options={solid,fill opacity=0}, only marks, forget plot]
table {%
6 -0.00113698380935055
6 0.00444378610576211
};
\addplot [darkorange25512714, forget plot]
table {%
0.75 0.00742285724530754
1.25 0.00742285724530754
};
\addplot [darkorange25512714, forget plot]
table {%
1.75 0.00771200218010142
2.25 0.00771200218010142
};
\addplot [darkorange25512714, forget plot]
table {%
2.75 0.00291447055651374
3.25 0.00291447055651374
};
\addplot [darkorange25512714, forget plot]
table {%
3.75 -0.00528388838393098
4.25 -0.00528388838393098
};
\addplot [darkorange25512714, forget plot]
table {%
4.75 0.00156127604206223
5.25 0.00156127604206223
};
\addplot [darkorange25512714, forget plot]
table {%
5.75 0.00146491844357492
6.25 0.00146491844357492
};
\end{axis}

\end{tikzpicture}